\documentclass[a4paper]{preprint}            
 \usepackage{amssymb}
\usepackage{amsfonts}
\usepackage{mhequ}
\usepackage{times}
\usepackage{theorem}
\usepackage{enumerate}

\newcommand{\A}{{\bf \mathcal A}}

\newcommand{\C}{{\mathcal C}}
\newcommand{\D}{{\rm I \! D}}
\newcommand{\E}{{\mathbf E}} 
\newcommand{\Epsilon}{{\mathcal E}}
\newcommand{\filt}{{\mathcal F}^{x_0}}
\newcommand{\F}{{\mathcal F}}
\newcommand{\g}{{\mathfrak g}}
\newcommand{\G}{{\mathcal G}}
\newcommand{\HH}{{\mathcal H}}
\newcommand{\K}{{\mathcal K}}
\newcommand{\LL}{\mathbb {L}}

\newcommand{\R}{{\mathbf R}}
\newcommand{\pnabla}{{\nabla\!\!\!\!\nabla}}

\newcommand{\1}{{\bf 1}}
\newcommand{\I}{{\mathcal I}}
\newcommand{\W}{\mathbf W}
\newcommand{\X}{{X\!\!\!\!\!\!X}}
\newcommand{\Y}{{Y\!\!\!\!Y}}

\def\paral{/\kern-0.55ex/}
\def\parals_#1{/\kern-0.55ex/_{\!#1}}
\def\n#1{|\kern-0.24em|\kern-0.24em|#1|\kern-0.24em|\kern-0.24em|}
\def\const{\rm {const.}}
\def\cf{\rm {c.f. }}
\def\eg{\textit{e.g. }}
\def\as{\textit{a.s. }}
\def\ie{\textit{i.e. }}
\def\s.t.{\mathop {\rm s.t.}}
\def\esssup{\mathop{\rm ess\; sup}}
\def\Ric{\mathop{\rm Ric}}
\def\div{\mathop{\rm div}}
\def\ker{\mathop{\rm ker}}
\def\Image{\mathop{\rm Image}}
\def\Dom{\mathop{\rm Dom}}
\def\dom{\mathop {\rm {\cal D}om}}
\def\id{\mathop {\rm Id}}
\def\Cyl{\mathop {\rm Cyl}}
\def\Span{\mathop {\rm Span}}
\def\trace{\mathop{\rm trace}}
\def\ev{\mathop {\rm ev}}
\def\constant{\mathop {\rm const.}}
\def\adjointconnection_#1{{\nabla^{1\prime}_{\!#1}}}
\def\adjointconnectionb_#1{{{\breve \nabla}^{1\prime}_{\!#1}}}

\newtheorem{theorem}{Theorem}[section] 
\newtheorem{proposition}[theorem]{Proposition}
\newtheorem{lemma}[theorem]{Lemma}
\newtheorem{corollary}[theorem]{Corollary}
\newtheorem{property}{Property}[section]

{\theorembodyfont{\rmfamily}
\newtheorem{definition}[theorem]{Definition}

\newtheorem{remark}[theorem]{Remark}

}

\newenvironment{enumerateroman}
{\begin{enumerate}[i.]}{\end{enumerate}}
\newenvironment{proof}{
 \noindent\textbf{Proof}\ }{\hspace*{\fill}
  \begin{math}\Box\end{math}\medskip}

\begin{document}

\author{K. D. Elworthy 
\& Xue-Mei Li }
\institute{}
\title{It\^o maps and analysis on path spaces}
\date{}
\maketitle

\begin{abstract}
We consider versions of Malliavin calculus on path spaces of compact manifolds with diffusion measures, defining Gross-Sobolev spaces of differentiable functions and proving their intertwining with solution maps, $\I$, of certain stochastic differential equations. This is shown to shed light on fundamental uniqueness questions for this calculus including uniqueness of the closed derivative operator $d$ and Markov uniqueness of the associated Dirichlet form. A continuity result for the divergence operator by Kree and Kree is extended to this situation. The regularity of conditional expectations of smooth functionals of classical Wiener space, given $\I$, is considered and shown to have strong implications for these questions. A major role is played by the (possibly sub-Riemannian) connections induced by stochastic differential equations: Damped Markovian connections are used for the covariant derivatives.
\end{abstract}


\section{Introduction}

A natural approach to geometric analysis on path spaces, or loop spaces, of manifolds is to base it on continuous paths with Brownian motion (or other) diffusion) measure. It became clear in the 1970's from the early work of L. Gross on analysis on Banach spaces with Gaussian measures, that in such analysis the differentiation should be restricted to differentiation in directions given by a certain Hilbert space, the Cameron-Martin space. These $H$-derivatives formed the basis of the highly successful Malliavin Calculus, see \eg Malliavin \cite{Malliavin78} \cite{Malliavin-book}. Key tools in this were the Sobolev spaces they generated,  see for example the books by Ikeda-Watanabe \cite{Ikeda-Watanabe} and Nualart \cite{Nualart-book}.
For paths on a Riemannian manifold, here based paths for simplicity, with Brownian motion measure, it was realized that the Cameron-Martin space should be replaced by Hilbert spaces of tangent vectors at almost all points of the path space: the so called Bismut tangent spaces. These are described in terms of parallel translation of the usual Cameron-Martin space of finite energy paths in the tangent space to $M$ at
 the base point, Jones-L\'eandre \cite{Jones-Leandre-91}. The parallel translation was 
 that of the Levi-Civita connection. This was extended to more general connections by
  Driver in \cite{Driver92} whose work led to rapid progress in creating a Sobolev calculus (depending on the choice of connections) over the path spaces and loop spaces see \eg Aida\cite{Aida97}, L\'eandre \cite{Leandre93}. When $M$ has curvature a major difficulty in this analysis comes from the non-holonomic nature of the Bismut tangent ``bundle''. There are no known ``local charts'' which adequately preserve the structure. The standard method, as in Driver \cite{Driver92} has been to use the stochastic development which gives a measure theoretic isomorphism of the Wiener space of based paths on Euclidean space with that on the manifold, and classically gives a diffeomorphism between the corresponding spaces of finite energy
   paths. However although this is smooth in the sense of Malliavin calculus, in general its derivative does not map the Cameron-Martin space to the Bismut Hilbert spaces. Moreover it seems clear from X-D Li \cite{X-D-Li-2003} that the composition of it with a differentiable function, say in $\D^{p,1}$, on the path space may not be differentiable on the flat space, \eg not in $\D^{q,1}$ for any $1\le q<\infty$: a loss of differentiability occurs. There is an intertwining formula, Thm 2.6 in Cruzeiro-Malliavin 
  \cite{Cruzeiro-Malliavin-96}, but it is for differentiation given by ``tangent processes'' not by tangent vectors.

There are also fundamental unresolved uniqueness problems in the calculus on these path spaces. The most basic is of the derivative operator itself: a standard approach is to take the closure in $L^p$ of the $H$-derivative defined on some initial domain of manifestly regular functions, \eg smooth cylindrical functions or bounded Fr\'echet differentiable functions with bounded derivatives. In Wiener space the result does not depend on any, reasonable, choice of such initial domain, Sugita \cite{Sugita85}. For paths on $M$ when there is non-zero curvature this is not known. Alhough there is a self-adjoint analogue $d^*d$ of the finite dimensional Laplace-Beltrami operator it is unknown whether it is essentially self-adjoint or whether it, or equivalently the associated Dirichlet form, has Markov uniqueness, taking the space of smooth cylindrical functions as initial domain. The latter concept relates to the uniqueness of a Markov process on the path space which would play the role of a Brownian motion (or Ornstein-Uhlenbeck process), see \eg Eberle \cite{Eberle-book} and \S\ref{section:Markov-Unique} below. Note that Aida has shown that such operators on certain finite co-dimensional submanifolds of Wiener space \cite{Aida93} are essentially self adjoint, and similarly for paths and loops on Lie groups  \cite{Aida95}. An earlier work of Costa has shown the essential self-adjointness for a larger core on such path groups.

Here we continue the approach of Aida-Elworthy \cite{Aida-Elworthy}, Aida \cite{Aida97}, Elworthy-LeJan-Li \cite{Elworthy-LeJan-Li-Tani} and Elworthy-Li \cite{Elworthy-Li-Hodge-1} using It\^o maps, \ie solution maps of stochastic differential equations, as substitutes for charts, and filtering techniques. We work with a fairly general class of, possibly degenerate, diffusion measures, with metric connections to define the Bismut tangent spaces. The stochastic differential equations are those whose solutions form the given diffusion process on $M$. We take $M$ compact and all coefficients smooth. The It\^o maps are then infinitely differentiable in the sense of Malliavin Calculus, but as with the stochastic development their H-derivatives will not in general map into the Bismut tangent spaces, nor can we expect there to be a `chain rule'  to say that  composition with them maps ``differentiable'' functions to `differentiable' functions. However if we restrict to stochastic differential equations whose associated connection, in the sense of Elworthy-LeJan-Li \cite{Elworthy-LeJan-Li-Tani} \cite{Elworthy-LeJan-Li-book}, agrees with the connection defining the Bismut tangent spaces, it turns out that such compositions are well behaved. The aim of this article is to describe this and its possible deficiences, and show what light it sheds on the fundamental uniqueness questions mentioned above. In particular we show that the latter are related to a question on Classical Wiener space concerning the regularity of differentiable functions after conditioning with respect to any of these It\^o maps, see Remark \ref{Remark-7.9} below. 

Since we are working in greater generality than usual, the calculus on path space is developed from scratch, based on the integration by parts formulae of Elworthy-Li \cite{Elworthy-Li-ibp} and Elworthy-LeJan-Li \cite{Elworthy-LeJan-Li-book}, \cite{Elworthy-LeJan-Li-CR}. This is done in section \ref{section-4} for scalars and with a covariant calculus in section \ref{section-covariant} more generally. In section \ref{section-2} the basic setup of Bismut tangent bundles, damped parallel translations, stochastic differential equations and associated connections, are described. In section \ref{section-3} there is the key result, Theorem \ref{th:I-star}, that if the connection associated to the stochastic differential equation is the same as that used to define the Bismut tangent spaces then the It\^o map can be used to pull back any measurable $\HH$-1-form $\phi$ on the path space of $M$ to an $H$-1-form on flat Wiener space.  We also give an explicit expression for the pull back $\I^*\phi$ as a stochastic integral. From this one obtains the pull back theorem, or chain rule, for functions in a Sobolev space $\D^{p,1}(\C_{x_0}M;\R)$, Theorem \ref{theorem-closability-2} and Corollary \ref{co:closed-range}. In particular if $f\in \D^{p,1}(\C_{x_0}M;\R)$ then $f\circ \I \in \D^{p,1}(\Omega;\R)$ and $d(f\circ \I)=\I^*(df)$ given some conditions on the connection.

In section \ref{section-divergence} we consider the divergence operator acting on $\HH$-vector fields and its intertwining by these It\^o maps. We show that $V$ lies in its domain if a certain pull back of $V$ is in the domain of the divergence on flat space, Corollary \ref{divergence-skorohod}. This enable us to extend the flat space result of Kree-Kree \cite{Kree-Kree} and see that $\D^{p,1}$ $\HH$-vector fields lie in the domain of the divergence, Theorem \ref{th:divergence-21}, a crucial result for our discussion of weak differentiability later.

In section \ref{section:Markov-Unique} we introduce weak differentiability and the weak Sobolev spaces $W^{p,1}$. Theorem \ref{th:weak_Sob} extends the chain rule to a precise intertwining:
$$f\in W^{p,1}(\C_{x_0}M;\R) \hbox { iff } f\circ \I\in \D^{p,1}(\C_0\R^m;\R).$$
The question is posed as to whether $W^{p,1}=\D^{p,1}$, as in flat space. In Theorem \ref{th:Markov-unique-2}, following Eberle \cite{Eberle-book}, this is shown for the case of $p=2$ to be equivalent to Markov uniqueness, after demonstrating $W^{2,1}={}^{0}W^{2,1}$, the latter being the weak Sobolev space used in \cite{Eberle-book}. A key step in the proof is to show that smooth cylindrical forms are dense in the space of $\D^{2,1}$ $\HH$-1-forms, Proposition \ref{dense-cylindrical}.

Other uniqueness questions are considered in section \ref{se:weak-uniqueness}. In
particular it is shown that the closure of the differentiation operator is independent of the initial domain if that domain contains $\Cyl$, the set of smooth cylindrical functions and consists of $BC^2$, twice Fr\'echet differentiable functions  whose derivatives are bounded, Corollary \ref{co:weak-unique}. However we are not able to prove the uniqueness of $d$ when its initial domain is allowed to contain general $BC^1$ functions. Some of these results for the special case of Brownian motion measures and Levi-Civita connections are summarised in \cite{Elworthy-Li-MUniqueness} correcting \cite{Elworthy-Li-CR-03}.

  In the case when it is possible to find a stochastic differential
 equation whose It\^o map has no redundant noise, all the main results
 in this article hold without the Condition $(M_0)$ which was often needed in the general situation. See section \ref{se:special-case}. In particular we have $\I^*d=d\I^*$ on $L^p$ for $1<p<\infty$ and the Markov uniqueness.  This applies to paths on Lie groups with left or right invariant connections and to paths on the orthonormal frame bundle of a Riemannian manifold with measure associated to the horizontal Laplacian. In this case our It\^o map is essentially the stochastic development map and our results are an extension of some of the isomorphism results by Fang-Franchi \cite{Fang-Franchi-97} for path spaces on Lie groups.

The culmination of section \ref{section-covariant} is Theorem \ref{higher-order-pull}
on the pull back by composition with $\I$ of higher order Sobolev spaces $\D^{q,k}$, and weak Sobolev spaces $W^{q,k}$, $k=1,2\dots$. As for $k=1$, in the weak case there is a precise intertwining. To differentiate these Sobolev spaces requires a connection on the Bismut tangent `bundle'. We use the `damped Markovian' connection. This was introduced in Cruzeiro-Fang \cite{Cruzeiro-Fang95} for Brownian motion measures with Levi-Civita connections. One key point in this work is how well they fit into this situation, in some sense being induced by the derivative of $\I$, Proposition \ref{pr:connection-3}.

Although we work in considerably greater generality, the main results here, Theorem \ref{th:I-star} and Theorem \ref{higher-order-pull} on intertwining, Theorem \ref{th:continuity-2} on the continuity of the divergence, Theorem 
\ref{th:Markov-unique-2} on Markov uniqueness, and Corollary \ref{co:weak-unique} ($BC^2$ functions are in $\D^{2,1}$), are essentially novel for the more standard case of Brownian motion measures and Levi-Civita connections (thought there is a version of Theorem 
\ref{th:I-star} in Elworthy-Li \cite{Elworthy-Li-Hodge-1}), as are the treatment of the covariant calculus in section \ref{section-covariant}, and the importance shown for the rather general problem of the smoothness of conditional expectations in classical Wiener space, Remark \ref{Remark-7.9}.

\tableofcontents

\section{Basic Assumptions}
\label{section-2}

 Let $M$ be a $C^\infty$ connected manifold of dimension
$n$. For simplicity assume it is compact. Otherwise some bounded
geometry  assumptions on the manifold and bounds on the
coefficients of the stochastic differential equations we consider
will need to be imposed. Let $\A$ be a smooth semi-elliptic second order operator with no zero order term. Assume its symbol $\sigma^\A: T^*M\to TM$ has constant rank $p$ so that its image is a sub-bundle $E$ of $TM$. It has a natural Riemannian metric induced by $\sigma^\A$. Let $\nabla$ be a metric connection
on $E$. Then $\A$ can be written in the following form:
\begin{equation}
\label{Generator-A}
\A f={1\over 2} \hbox{trace}_E \nabla_-(df\vert_{E})+L_{A}(f)
\end{equation}
where $A$ is a smooth vector field. Denote by $\mu_{x_0}$ the law of
 the Markov process $(x_t:~0\le t\le T)$ corresponding to $\A$ with initial value $x_0$ for some point $x_0\in M$ and fixed $T>0$.

Consider $$\C_{x_0}M=\big\{\sigma: [0,T]\to M\;\big\vert \;\sigma(x_0)=0,
 \sigma \hbox{ is } C^0    \big\},$$
the space of continuous paths on $M$ starting from $x_0$ equipped
with the probability measure $\mu_{x_0}$.

\subsection{An SDE which induces the connection $\nabla$}\label{section-1}

 The underlying probability space $\Omega$ will be taken to be the canonical
  space $\C_0\R^m$,  given by
\begin{equation}
\label{Omega}
\Omega=\C_0\R^m=\big\{\omega\colon[0,T]\to\R^m \; \big\vert \,\omega(0)=0,
\hbox{ $\omega$ is continuous} \big\},
\end{equation}
some natural number $m$. It is equipped with the Wiener measure
$\mathbf P$ and its natural filtration $\{\F_*\}$. Let  $\displaystyle{\{B_t: 0\le t\le T\}}$
be  the canonical Brownian motion on $\R^m$, that is
 $B_t(\omega)=\omega(t)$, the evaluation map.

Denote by  $\LL(E; F)$  the space of bounded linear maps between
linear spaces $E$ and $F$ and let  $X: \R^m\times M\to TM$ be $C^\infty$ with
 $X(x)\in \LL(\R^m; T_xM)$ for each $x$.  For each $x\in M$, $\Image[X(x)]$
 inherits an inner product.
 We shall choose $X$ so that $\Image[X(x)]=E_x$ as a Hilbert space.
Let $C^r\Gamma E$ be the space of $C^r$ sections of $E$.  For $e$ in
$\displaystyle{\R^m}$ let $X^e$ be the section of $E$ given by
$\displaystyle{X^e(x)=X(x)(e)}$ and $Y: E\to \R^m$ the adjoint of $X$. Note that $X(x)Y(x)(v)=v$ for all $v \in T_xM$.  
Write $\ker X(x)$ and $[\ker X (x)]^\perp$ respectively for the kernel of the map $X(x)$ 
and its orthogonal complement.
 The result on which, Elworthy-LeJan-Li \cite{Elworthy-LeJan-Li-book},
 and this article are based is the following, c.f. Quillen \cite{Quillen}, Narasimhan-Ramanan \cite{Narasimhan-Ramanan-61}:

\begin{proposition}
\label{theorem-connection}
(Elworthy-LeJan-Li \cite{Elworthy-LeJan-Li-Tani}, \cite{Elworthy-LeJan-Li-book},
 c.f. Quillen\cite{Quillen})
  For each such map $X:\R^m\times M\to TM$ there is a unique connection 
  $\breve \nabla$ on $E$ such that
  \begin{equation}
  \label{the-connection}
  \breve \nabla_v X^e=0, \quad \forall   v\in T_yM, y\in M, e\in[\ker X(y)]^\perp
  \end{equation}
  This connection is metric. In fact 
\begin{equation}
\label{connection}
\breve\nabla_vU
=X(x)d\left (Y(U(\cdot))\right)(v), \quad v\in T_xM,  U\in C^1\Gamma E.
\end{equation}
Furthermore all metric connections on $E$ can be obtained this way for
 some $X$ and some number $m$.
\end{proposition}
\begin{itemize}
\item{\bf Assumption ($X$).} 
By this proposition we can and will suppose from now on that map $X$
induces the Riemannian metric and the connection $\nabla$ on $E$. Hence (\ref{the-connection}) holds for $\nabla$. 
\end{itemize}
For  $A$ as in (\ref{Generator-A}), consider the stochastic differential equation
\begin{equation}\label{sde}
dx_t=X(x_t)\circ dB_t +A(x_t) dt,
\hskip 20pt  0\le t\le T.
\end{equation}\
It induces the diffusion measure $\mu_{x_0}$ on $M$.  
 
\subsubsection{ Examples \cite{Elworthy-LeJan-Li-book} }
\label{sec:Examples}
 {\it Example 1 (\it Gradient S.D.E.).}
 Consider a Riemannian manifold $M$  isometrically immersed in
$\R^m$ with immersion $j$. Set $X^e\equiv X(\cdot)(e)=\nabla\langle
j(\cdot), e\rangle$. Then $X(x): \R^m\to TM$ is
 the orthogonal projection of $\R^m$ to $T_xM$. The stochastic differential
equation  has solutions which are Brownian motions on $M$
and it induces the Levi-Civita connection on $M$ as its associated
connection.

{\it Example 2. (\it Left invariant S.D.E.).}
 Let $M$ be a Lie group with left invariant metric,  identity
$e$ and Lie algebra ${\mathfrak g}:=T_eG$. Let $(B_t)$ be a Brownian motion
 on ${\mathfrak g}$. The connection associated  to the left invariant
 stochastic differential equation
$$dx_t=X(x_t)\circ dB_t$$
is the flat left invariant connection. Here $X(e):\R^n\to \G$  is some
 isometry and $X(g)e=(TL_g)X(e)$.  The solution of the equation is a
process on the Lie group whose filtration is the same as
that of the noise $(B_t)$. It is a Brownian motion if the metric is
bi-invariant.

{\it Example 3 (\it Symmetric space S.D.E.).}
 Let $M=G$ be a Lie group with bi-invariant metric. It has a standard
symmetric space structure: $$G={G\times G\over \{(g,g): g\in G\}}$$
where $G\times G$ acts on $G$ as follows:
$$(g_1,g_2)x=g_1xg_2^{-1}.$$
 Denote by $L_x$ and $R_x$  respectively the left and right
group multiplications. Consider the stochastic differential equation on $G$
$$dx_t={1\over \sqrt 2}TL_{x_t}\circ dB_t-{1\over \sqrt 2}TR_{x_t}\circ dB_t^\prime$$
where $(B_t)$ and $(B_t^\prime)$ are two independent Brownian motions on
$\mathfrak g$.
The corresponding connection is the Levi-Civita connection on $G$.

\bigskip

{\it Example 4 (Canonical SDE on frame bundles).}
Let $N$ be an $m$ dimensional  Riemannian manifold and let $M$ be its
 orthonormal frame bundle, $M=ON$, with $\pi: ON\to N$ the projection.
 Using the  Levi-Civita connection for $N$ consider the canonical
 stochastic differential equation on $ON$. Then $X(u)e=H_u(u(e))$
where $H_u: T_{\pi(u)}N\to T_uON$ denotes the horizontal lift map.
Then  $E$ is the horizontal tangent bundle of $ON$ and $p=m=\dim N$. The connection on $E$ is the flat connection induced by
the trivialization $X$. The solutions to the S.D.E. on $ON$
project to Brownian motions on $N$, and are the horizontal lifts
of those Brownian motions.

\subsection{The Covariant Differentiation operator ${\D\over dt}$ } 
\label{subsection-2.1}
There is an adjoint {\it semi-connection}, $\nabla^\prime$, of $\nabla$. For each
 smooth vector field $V$ on $M$ this gives a derivative
$$\nabla^\prime_u V\in T_y M$$
for each $u\in E_y$, $y\in M$. It is defined by
\begin{equation}
\nabla^\prime_u V= \nabla_v U+[U,V](y)
\end{equation}
for $v=V(y)$ and $U$ any smooth section of $E$ with $U(y)=u$.

 Using $\nabla$ there are parallel translations along smooth paths
 $\sigma$ in $M$
$${\parals_{t}}={\parals_{t}}(\sigma):
E_{\sigma(0)}\to E_{\sigma(t)}$$
and these preserve the inner products.
 Using $\nabla^\prime$ we obtain parallel translation
$$\parals_{t}^\prime=\parals_{t}^\prime(\sigma):
T_{\sigma(0)}M\to T_{\sigma(t)}M  $$
along smooth paths which are `horizontal'. (A path $\sigma$ is horizontal
 if $\dot \sigma(t)$ belongs to  $E_{\sigma(t)}$ for each $t$.)
 There are also the operators
${D\over dt}$ and ${D^\prime \over dt}$
$${D\over dt}U_t={\parals_{t}}{d\over dt}({ \parals_{t}}^{-1}U_t)\in E_{\sigma(t)}$$
$${D^\prime \over dt}V_t
={\parals_{t}^\prime}{d\over dt}({\parals_{t}^\prime}^{-1}V_t) \in
T_{\sigma(t)}M$$ defined for vector fields along $\sigma$, for
$U_t\in E_{\sigma(t)}$ each $t$, and for $\sigma$ horizontal in
the case of ${D^\prime\over dt}$.

\bigskip

We will need `damped' versions of these operations. Let $Z$ be a vector field on $M$.
When  $Z(x)\in E_x$ for each $x$, the
damped parallel translation $W_t^Z\equiv W_t^Z(\sigma):T_{\sigma(0)}M\to T_{\sigma(t)}M$ along a horizontal smooth path $\sigma$ is defined by
\begin{equation}
\label{damped-parallel-1}
\left\{\begin{array}{ll}
{D^\prime\over dt}[W_t^Z(v_0)]
&=-{1\over 2}\Ric^{\#} (W_t^Z(v_0))+ \nabla_{W_t^Z(v_0)}Z,
\hskip 20pt 0\le t\le T\\
W_0^Z(v_0)&=v_0.
\end{array}\right.
\end{equation}  Here
${\Ric}^{\#}: TM\to E$ is defined by the Ricci
curvature ${\Ric}$ corresponding to the connection $\nabla$:
$\langle{\Ric}^{\#}(u), v \rangle_y={\Ric}_y(u,v)={\trace}_E \langle R(u,-)-, v\rangle_y$.
 Under these conditions the corresponding operator
 ${\D \over dt}$ on vector fields along $\sigma$ is given by
$${\D \over dt} V_t= W_t^Z{d\over dt}\left(({W_t^Z})^{-1}V_t\right).$$
Thus
\begin{equation}
\label{damped-covariant}
{\D\over dt}
={D^\prime\over dt}+{1\over 2}{\Ric}^{\#}-\nabla Z.
\end{equation}

In this case the damped (and undamped) parallel translation
is defined almost surely along the  sample paths of the solutions to our
stochastic differential equation (\ref{sde}), with ${D^\prime\over dt}$
and ${\D\over dt}$ being defined correspondingly on suitable vector
fields along the paths.

If $Z$ is not a section of $E$ then the solution
paths are not `horizontal' and it is convenient  to introduce an
auxiliary connection  $\nabla^1$ on $TM$. To obtain
this take a Riemannian metric on $TM$. Let $E^\perp$ be the orthogonal
bundle to $E$ in $TM$ and take any metric connection $\nabla^\perp$
on $E^\perp$. Set
\begin{equation}
\label{nabla-1}
\nabla^1=\nabla +\nabla^\perp
\end{equation}
and let ${\nabla^1}^\prime$ be its adjoint. Now extend the definition
of  ${\D\over dt}$ to define
\begin{equation}
\label{covariant-3}
{\D\over dt}V_t:=
{{D^1}^\prime \over dt}V_t+{1\over 2}{\Ric}^{\#}V_t-
\nabla_{V_t}^1 Z
\end{equation}
for any suitably regular vector field $V$ along the paths of any continuous
semi-martingale on $M$.

In particular ${\D\over dt}$ is defined $\mu_{x_0}$-almost surely for suitably
regular vector fields along the elements of $\C_{x_0}M$. It follows from
Proposition 3.3.9 of Elworthy-LeJan-Li \cite{Elworthy-LeJan-Li-book} and the Girsanov-Maruyama theorem
that as such it depends only on $\nabla$ and $Z$ and not on the choice of $E^\perp$ or $\nabla^\perp$ provided that
$Z-A\in \Gamma(E)$ where $A$ is the drift coefficient of the SDE.
 Consequently  for such $Z$ the solution to
 \begin{equation}
\label{damped-parallel-2}
{\D\over dt}v_t=0
\end{equation}
given $v_0\in T_{x_0}M$ is defined along $\mu_{x_0}$-almost all paths and
 is  independent of the choice of the auxiliary connection on $E^\perp$.
When $Z=A$ it shall be denoted $W_t(v_0)$ to extend that defined by
 (\ref{damped-parallel-1}). With these extensions it remains true that
\begin{equation}
{\D\over dt}v_t=W_t{d\over dt} (W_t^{-1}(v_0)).
\end{equation}

\subsubsection{Condition $(M_0)$}
Some additional conditions will sometimes be imposed on our
connection $\nabla$ and on the stochastic differential equation. 
We are given a Riemannian metric only on $E$
and so will formulate the conditions in terms of that metric,
avoiding using the metric we imposed on $TM$, even though for
a compact manifold $M$ that is not really essential. \hskip 10pt
\begin{itemize}
\item {\bf Condition $(M_0)$ }: The damped parallel translation $W_t$
satisfies:
\begin{enumerate}
\item[(i)]
$\sup_{0\le t\le T}\left\vert ( {W_t^{-1}} \vert_{E_{x_t}})
\right\vert _{\LL(E_{x_t}; T_{x_0}M)}\in L^\infty$
   \hskip 10pt and\\
 \item[(ii)] $\sup_{0\le t\le T}\left\vert \nabla_{W_t(-)}X
\right \vert _{\LL(T_{x_0}M; \LL(\R^m; E_{x_t}))}\in L^\infty$
\end{enumerate}
\item {\bf Condition ($M$) }:
 The adjoint connection $\nabla^\prime$ is metric
for some Riemannian metric on $TM$, (which we will denote by
$\langle \cdot , \cdot\rangle^\prime$).
\end{itemize}
Note that if $E=TM$ condition (M) holds  with $\langle\cdot , \cdot
\rangle^\prime =\langle\cdot , \cdot\rangle$ if and only if $\nabla$
is torsion skew symmetric as described by Driver in \cite{Driver92}.
In particular Condition ($M$) holds for the SDE's in Examples 1-3, section~\ref{sec:Examples}. 
For examples where it does not hold
 see  Elworthy-LeJan-Li \cite{Elworthy-LeJan-Li-book}, in which
there is also the following result (Proposition 3.3.11, p72):

\begin{proposition} \cite{Elworthy-LeJan-Li-book}
\label{norm-TI}
For compact $M$,
$$\sup_{0\le s\le T}\left\vert W_s\right\vert
_{\LL( T_{x_0}M; T_{x_s}M)} \hskip 6pt \hbox{and} \hskip 6pt
\sup_{0\le s\le T}\left\vert W_s^{-1}\right\vert _{\LL( T_{x_s}M;
T_{x_0}M)}$$ lie in $L^p$ for all $1\le p<\infty$. If also
condition ($M$) holds then both are in $L^\infty$. (Here we are
using any Riemannian metric on $M$.)
\end{proposition}
From this we see immediately that condition $(M)$ implies
condition $(M_0)$
 in the compact case under consideration.

\subsection{The $L^2$ tangent bundles $L^2\Epsilon$, $L^2T\C_{x_0}M$ 
and the Bismut tangent bundle $\HH$} \label{se:L2_tangent_spaces}

Recall that $\C_{x_0}M$ is a $C^\infty$ Banach manifold,  Eells \cite{Eells58}, see Eliasson \cite{Eliasson67},
and its tangent space $T_\sigma \C_{x_0}M$ at a path $\sigma$
can be identified with the following space of vector fields along it:
$$T_\sigma \C_{x_0}M
=\big\{v: [0,T]\to TM \; \big \vert \,  v_t\in T_{\sigma(t)}M, v_0=0\big\}.$$
By the $L^2$ tangent space $L^2\Epsilon_\sigma$ at $\sigma$ we mean
 the following  set of measurable vector fields along $\sigma$:
\begin{equation}
\label{L2-tangent}
L^2\Epsilon_\sigma=\left\{ v: [0,T]\to E \;\;
  \big\vert\;\; v_t\in E_{\sigma(t)}, \;
\vert v_\cdot\vert_{L^2\Epsilon}<\infty\right\}.
\end{equation}
where $$\vert v_\cdot\vert_{L^2\Epsilon}:= \left(\int_0^T\vert v_s\vert ^2 \, ds\right)^{1/2}.$$
These form the fibres of a smooth Hilbert bundle $L^2\Epsilon$ over $\C_{x_0}M$. It is associated to
 the principal bundle $$C_{x_0}OE\to \C_{x_0}M$$
where $\pi: OE\to M$ is the  orthonormal frame bundle of $E$ and
$$C_{x_0}OE=\{u: [0,T]\to OE\;  \big |\; u(0)\in \pi^{-1}(x_0)\}.$$
Given the choice of a Riemannian metric on $TM$ extending that of $E$ we also have 
the Hilbert subbundle of $L^2$ tangent vectors $L^2TC_{x_0}M$, obtained as 
$L^2\Epsilon$ but using $TM$ rather than $E$. Then $L^2\Epsilon$ is
a subbundle of $L^2T\C_{x_0}M$.  Let 
$$\Pi: L^2T\C_{x_0}M\to L^2\Epsilon$$ denote the orthogonal projection.

Using the metric connection $\nabla$ which we have imposed on $E$, we can define a family  of subspaces $ \HH_\sigma\subset T_\sigma \C_{x_0}M$:
\begin{equation}\label{HH}
 \HH_\sigma:=\left\{
v\in T_\sigma \C_{x_0}M \left \vert \hskip 5pt
 {\D\over dt}v_t\in E_{\sigma(t)},
\int_0^T\vert {\D v_t\over dt}\vert^2_{\sigma(t)} \; dt<\infty
\right.\right\}\end{equation}
(with the usual convention of absolute continuity after translation back
to $T_{x_0}M$). This is a Hilbert space under the obvious inner product
$$\langle u,v \rangle_\sigma=\int_0^T \left \langle {\D u_t\over dt},  {\D v_t\over dt} \right \rangle_{\sigma(t)}\; dt.$$

Note also that ${\D\over dt}$ determines 
an isometry of $\HH_\sigma\to L^2\Epsilon_\sigma$ for almost all $\sigma$, with inverse
$$\W_\cdot: L^2\Epsilon_\sigma\to \HH_\sigma$$
given by 
\begin{equation}
\label{wt}
\W_t(v)=W_t\int_0^t W_s^{-1} v_s ds.
\end{equation}
Let $\HH=\sqcup_{\sigma}\HH_\sigma$. Then it inherits a vector bundle structure (over a subset
of full measure in $\C_{x_0}M$) from $L^2\Epsilon$ via ${\D\over dt}$ as does its dual $\HH^*=\sqcup_\sigma\HH_\sigma^*$.

In particular an  $L^p$ $H$-form  (or written as $\HH$-form) $\phi$ on $\C_{x_0}M$ is an $L^p$ section  of $\HH^*$,
 i.e. an assignment of  $\phi_\sigma:\HH_\sigma\to \R$, continuous linear,
 for almost all   $\sigma$ in $\C_{x_0}M$, measurable in $\sigma$ in the sense
that $\sigma\mapsto \phi_\sigma \left({\D\over d \cdot}-\right)$
is a measurable section of the dual of the vector bundle $L^2 \Epsilon$
with
$$\n{\phi}_{L^p}:=\int_{\C_{x_0}M}\vert \phi_\sigma\vert^p_{\HH_\sigma^*}\,
d\mu_{x_0}(\sigma)<\infty.$$
Let $L^p\Gamma\HH^*$ be the space of equivalence classes of $L^p$
$H$-forms.

\begin{remark}
Suppose $V\in L^p\Gamma\HH$, the space of $L^p$ $H$-vector fields on $\C_{x_0}M$. Then,
for any inner product on $T_{x_0}M$ and almost all $\sigma\in \C_{x_0}M$
\begin{eqnarray*}
\sup_t|W_t^{-1}V_t(\sigma)|&=&\sup_t \left|\int_0^t W_s^{-1}{\D\over ds}
V_s(\sigma) ds\right|\\
&\le& T^{1\over 2} \sup_t |W_t^{-1}|_{\LL(E_{\sigma(t)}; T_{x_0}M)}\|V\|_{\HH_\sigma}
\end{eqnarray*}
and so $\sup_t|W_t^{-1}V_t(\sigma)|$ is in $L^p$ if condition $(M_0)$ holds.
\end{remark}
\section{Pull backs of $H$-forms by It\^ o maps.}
\label{section-3}

\subsection{The derivative $T\I$ of the It\^o map and $\overline{T\I}$}
 Let
$$\I: \C_0\R^m\to \C_{x_0}M$$
$$\I(\omega)_t=x_t(\omega)$$
 be the It\^o map of (\ref{sde}) for
 $\{\xi_t: 0\le t\le T\}$  the solution flow of (\ref{sde}) and
$x_t(\omega)=\xi_t(x_0,\omega)$. For each $\omega\in \C_0\R^m$, let
$$T_\omega\I: H\equiv L_0^{2,1}\R^m \longrightarrow
 T_{x_\cdot(\omega)} \C_{x_0}M$$
 be its $H$-derivative in the sense of Malliavin calculus. Strictly speaking
 conventional Malliavin calculus just gives a derivative at each time $t$
$$T_\omega\I_t:  H\longrightarrow T_{x_t(\omega)} \C_{x_0}M.$$
However there is the formula, due to Bismut, for $v_t(\omega)=T_\omega\I_t(h)$,
$ h\in H$:
\begin{equation}\label{Bismut}
v_t=T_{x_0}\xi_t\int _0^t(T_{x_0}\xi_s)^{-1} X(x_s) \dot h_s ds,
\hskip 15pt 0\le t\le T
\end{equation}
where $T_{x_0}\xi_t: T_{x_0}M\to T_{x_t}M$  is the derivative  at $x_0$ of $\xi_t$. This shows that we do have, for almost all
 $\omega\in \C_0\R^m$, a continuous linear version
 $T_\omega\I: H\to T_{x_\cdot(\omega)}\C_{x_0}M$.
Moreover $\sup_t |T\I_t|_{\LL(H;T_{x_t}M)} $ lies in $L^p$ for all
$1\le p<\infty$, (c.f. Proposition \ref{norm-TI} below), for any Riemannian
 metric on our compact manifold $M$.

One of the key points in our discussion will be the decomposition of the
`noise' $\{B_t :0\le t\le T\}$ into `redundant' and `relevant' parts
\begin{equation}
\label{noise-decomposition} dB_t=\tilde {\parals_t} d\tilde B_t
+\tilde{\parals_t} d\beta_t, 
\end{equation}
as described in Elworthy-Yor \cite{Elworthy-Yor} for gradient systems and
Elworthy-LeJan-Li  \cite{Elworthy-LeJan-Li-Tani}, \cite{Elworthy-LeJan-Li-book} more generally. Here

\begin{enumerate}
\item[(i)]
$\tilde{\parals_t}(\omega): \R^m\to \R^m$ is an orthogonal transformation of
$\R^m$, mapping $\ker X(x_0)$ to $\ker X(x_t(\omega))$, given by parallel
translation along $\{x_t: 0\le t\le T\}$ using a connection on the trivial
 $\R^m$-bundle over $M$, canonically determined by $X$.

\item[(ii)]
$\tilde B_t:=\int_0^t \tilde \parals_s^{-1} K^\perp (x_s) dB_s$,  for $K^\perp(x)$ the orthogonal projection of $\R^m$ onto $[\ker(X(x))]^\perp$; so $\{\tilde B_t: 0\le t\le T\}$  is a Brownian motion on $[ker X(x_0)]^\perp$.
It has the same filtration as that of $\{x_t: 0\le t\le T\}$.

\item[(iii)]
$\beta_t:=\int_0^t \tilde \parals_s^{-1} K(x_s)dB_s$ with $K(x)=\1-K^\perp(x)$; so $\{\beta_t: 0\le t\le T\}$ is an $\F_*$-Brownian motion on $\ker X(x_0)$,  independent of $\{x_s: 0\le s\le T\}$. From the point of the view of the solution
$\{x_t: 0\le t<\infty\}$ it is the `redundant noise'.

\end{enumerate}

From Elworthy-LeJan-Li \cite{Elworthy-LeJan-Li-CR} or equation (4.16) p79 of
\cite{Elworthy-LeJan-Li-book}, we have the covariant It\^o
equation for $v_t\equiv T\I_t(h)$ any $h\in H$, using the
connection ${\nabla^1}^\prime$ on $TM$, defined via (\ref{nabla-1}),
\begin{equation}
\label{eq-filter} {D^1}^\prime v_t = \nabla_{v_t} X(\tilde
{\parals_t} d\beta_t) -{1\over 2} {\Ric}^{\#}(v_t)dt
+\nabla^1_{v_t} A dt+X(x_t)\dot h_t dt
\end{equation}
which may be written, using notation (\ref{covariant-3}),
\begin{equation}
\label{eq-filter-2} \D v_t= \nabla _{v_t} X(\tilde
{\parals_t} d \beta_t)+X(x_t)\dot h_t dt.
\end{equation}
The equation (\ref{eq-filter}) comes from (\ref{noise-decomposition})
and the defining property, (\ref{the-connection}),  of the connection.
\bigskip

For almost all $\sigma\in C_{x_0}M$ and $h\in H$ define
\begin{equation}
\overline{T\I}_\sigma(h)=\E\left\{T_\omega\I(h)| x_\cdot(\omega)=\sigma\right\}.
\end{equation}
From (\ref{eq-filter-2}), as in \cite{Elworthy-LeJan-Li-book}, we obtain a key property

\begin{property}[ Elworthy-Li \cite{Elworthy-Li-forms-CR}\cite{Elworthy-Li-Hodge-1}]
\label{property-1}
Suppose the connection defined by the SDE (\ref{sde}) is the same as that defining $\HH$.
Then the map $\overline{T\I}_\sigma$ gives a projection
$$\overline{T\I}_\sigma: H\to \HH_\sigma$$
for almost all $\sigma\in C_{x_0}M$. It is given by
$$\overline{T\I}_\sigma(h)_t
=W_t\int_0^t W_s^{-1} X(\sigma_s) \dot h_s ds$$
with isometric right inverse $v\mapsto \int_0^\cdot Y_{\sigma(s)}\left({\D\over ds}v_s\right)ds$.
\end{property}

\subsection{Some useful lemmas}

\begin{lemma}
\label{stopping-times}
  Let $(M_t, 0 \le t\le T)$ be a continuous local martingale 
with respect to some filtration $\G_*$, with values in $\LL (\R^k ; G )$ for some separable Hilbert space $G$. Suppose the tensor quadratic variation of $(M_t)$ has a continuous density with respect to $t$. Then the map
 $$f \mapsto \int_0^T {dM}_s ( f_s ).$$
is continuous in probability as a map
\[ L^0 \left( \Omega,\mathcal{G}_0,P; \LL^2 ( [ 0, T ] ;\R^k ) \right)
   \longrightarrow L^0 \left( \Omega,\G_T,P; G \right) \]
   \end{lemma}

\begin{proof}
  Suppose $\{f_n, n\ge 1 \}$ is a sequence of $\G_0$ measurable
  functions converging in probability to $f$. For $m = 1, 2, 3, \dots$, set
   $\tau_m = \inf_{t>0} \{ \sup_n \{ \int_0^t |f_{n}(s) |^2 ds \} \ge m \}$, 
  giving it the value $T$ if the set is empty.
  Note that by going to a subsequence which converges almost surely we can
  assume that $\sup_n \{ \int_0^T | f_{n}( s) |^2 {ds} \}$ is almost
  surely finite and so these times increase to $T$ almost surely. They are
  also $\G_0$-measurable and so can be used as stopping times. The
  processes$\{\chi_{[0,\tau_m)}(\cdot) f_n  \}_{n=1}^\infty$ are bounded and so converge in $L^p$ to $\chi_{[0,\tau_m)}(\cdot)f$ for each $m$ and $p<\infty$. Their stochastic integrals, after localisation, will then converge in $L^p$ and the result follows.
\end{proof}

  The next proposition extends the main technical tool
used in Aida-Elworthy \cite{Aida-Elworthy}. Compactness of $M$ is not
used though  non-explosion of the underlying diffusion needs to be assumed.
First we record an easy consequence of the Burkholder-Davis-Gundy
inequalities.

\begin{lemma}\label{lemma-BDG}
For $0<p<\infty$ let $c_p, C_p$ be the constants in the
Burkholder-Davis-Gundy inequalities. If $\{Z_t: 0\le t\le T\}$ is a
real-valued continuous local martingale with respect to a filtration
$\{\G_t: 0\le t\le T\}$ with $Z_0=0$, then almost surely
$$c_p\E\left\{\langle Z, Z\rangle_T^{p\over 2}\hskip 4pt \vert \hskip 6pt
\G_0\right\}
\le \E\left\{\sup_{0\le t\le T} \vert Z_t\vert ^p
\hskip 4pt \vert \hskip 6pt
\G_0\right\}
\le C_p\E\left\{\langle Z, Z\rangle_T^{p\over 2}\hskip 4pt \vert \hskip 6pt
\G_0\right\}.$$
\end{lemma}

\begin{proof} 
Let $\lambda$ be non-negative, $\G_0$-measurable, and bounded.
Then $\{\lambda^{1\over p} Z_t: 0\le t\le T\}$ is a $\G_*$-local martingale
to which we can apply the Burkholder-Davis-Gundy inequalities to see
$$c_p\E \lambda \langle Z, Z\rangle_T^{p\over 2}
\le \E \lambda \sup_{0\le t\le T} \vert Z_t\vert ^p
\le C_p\E\lambda \langle Z, Z\rangle_T^{p\over 2} $$
giving the result. \end{proof}

Set $\F^{x_0}=\sigma\{x_s: 0\le s \le T\}$.
\begin{proposition}
\label{proposition-3.2} Assume condition $(M_0)$ holds. Then for
all $1\le p<\infty $ there is a constant $\alpha _p$ with
$$\E \left\{ \sup_{0\le s\le T \,}\left| W_s^{-1}T\I  %
_s(h)\right| _{T_{x_0}M}^p\;\Big\vert \;\F ^{x_0}\right\} \le
\alpha _p\left\| h\right\| _H^p,\hskip 16pt\mathrm{\,all\,}\hskip 5pt h\in
H\,\,\,\,\,\,\,\,\,\,\mathrm{a.s.}$$
\end{proposition}

\begin{proof} Take $h\in H$. We only need to show the inequality
for $p\ge 2$. From (\ref{eq-filter-2}) we have the It\^o equation for
 $u_t:=W_t^{-1}\left( T\I  _t(h)\right) $:

\begin{equation}
\label{eq-TI}
du_t=W_t^{-1}\nabla _{W_t(u_t)}X\left( \tilde{\parals_t}
   d\beta _t\right)    +    W_t^{-1}X(x_t)(\dot{h}_t)dt
\end{equation}
giving
\begin{equs}[3]
\sup_{0\le t\le \tau \,} |u_t|^p
&\le 2^{p-1}\sup_{0\le t\le \tau 
}\left| \int_0^tW_s^{-1}\nabla_{W_s(u_s)}X\left( \tilde{\parals_s}
d\beta_s\right) \right| ^p\\
&+2^{p-1}\sup_{0\le t\le \tau }\left|
\int_0^tW_s^{-1}X(x_s)(\dot{h}_s)ds\right| ^p  \label{TI-norm-3}
\end{equs}
for any $0\le \tau \le T$. Set $\G_t=\F_t\vee\F^{x_0}$. Then $(\beta_\cdot)$ is a $\G_*$-Brownian motion and so we can apply Lemma \ref{lemma-BDG} to give
\begin{eqnarray*}
&&\E\left\{ \sup_{0\le t\le \tau \;}\left| \int_0^tW_s^{-1}\nabla_{W_s(u_s)}X\left( \tilde{\parals_s} d\beta _s\right) \right|
^p \;\Big\vert\; \F ^{x_0} \right\}  \\
&\le &C_p\E \left\{ \left( \int_0^\tau \sum_j\left| W_s^{-1}\nabla_{W_s(u_s)}X^j\right| ^2\,ds\right) ^{\frac p2}\;\Big\vert\; \F ^{x_0}\right\}
\end{eqnarray*}
for $X^j(x)=X(x)(e^j)$ where $e^1,\dots ,e^m$ is an orthonormal
base for $\R^m$. Since $p\ge 2$,  condition $(M_0)$ plus
Jensen's inequality gives
\begin{eqnarray*}
&&\E \Big\{ \Big( \int_0^\tau \sum_j\left| W_s^{-1}\nabla_{W_s(u_s)}X^j\right| ^2\,ds\Big) ^{\frac p2}\Big| \F^{x_0} \Big\}  \\
&& \le \constant \Big(\n {\sup_{0\le s\le \tau}
\left| W_s^{-1}\,_{|_{E_{x_s}}}\right| _{\LL(E_{x_s};T_{x_0}M)}}
 _{L^\infty }\Big)^p \cdot\\
&&\Big( \n {\sup_{0\le s\le \tau}\left| \nabla_{W_s(-)}X\right|
_{\LL(T_{x_0}M;\LL_2(R^m;E_{x_s}))}} _{L^\infty }\Big)^p
  \big(\tau ^{p/2-1}\big)\E\big\{ \int_0^\tau |u_s|^p\,ds\Big| \F ^{x_0} \big\}\\
&&\le \mathrm{const.}\,\int_0^\tau \E \Big\{
\sup_{0\le r\le s}|u_r|^p\;\big\vert \F ^{x_0}\; \Big\}
\,ds .
\end{eqnarray*}
Applying condition $(M_0)$ to the second term on the right hand
side of (\ref {TI-norm-3}) we see that (\ref{TI-norm-3}) leads to
$$\E \left\{ \sup_{0\le t\le \tau }|u_t|^p\;\big| \F %
^{x_0}\; \right\} \le \mathrm{const.}\,\int_0^\tau \E \left\{
\sup_{0\le r\le s}|u_r|^p\;\big\vert \F ^{x_0}\; \right\}
\,ds+  \mathrm{const.}\,\left( \left\| h\right\| _H\right) ^p$$
and the result follows by Gronwall's lemma. 
\end{proof}

\subsection{The pull back map $\I^*$ and the push forward map $\overline{T\I(-)}$}

If $f: \C_{x_0}M\to \R$ is Fr\'echet  $C^1$ with bounded derivative
$df: T\C_{x_0}M\to \R$ we see that $\I^*(df)_\omega:=df\circ T_\omega\I$
is almost surely defined as a continuous linear functional on $H$,
(and by the usual approximation techniques this is
 $\bar d(f\circ \I)_\omega$).  Similarly  we can pull back any
geometric 1-form $\phi: T\C_{x_0}M\to \R$ to obtain an $H$-form
$\I^*(\phi)=\phi(T\I-)$ on $\C_0\R^m$. If $f\in L^2(\C_{x_0}M;\R)$ is
an arbitrary element in $\Dom(\bar d)$ we have now the $\HH$-form
$\bar d f_\sigma: \HH_\sigma\to \R$ for almost all $\sigma \in
\C_{x_0}M$. Set $v_t=T\I_t(h)$ for $h\in H$. From Bismut's formula (\ref{Bismut}) we cannot expect $v_\cdot$ to be in $\HH$ and thus the usual pull back map $\I^*: \Lambda^1\to L^0(\Omega;  H^*)$
 defined by $\I^*(\phi)(h)=\phi(T\I(h))$, 
 on geometric differential 1-forms
does not obviously extend to $\HH$-forms.  In particular it is not at all clear that we can define
$\I^*(\bar df)$. 
As shall be seen below we
 will need to interpret $\phi(v_\cdot)$ as a stochastic integral.

\begin{theorem}
\label{th:I-star} 
Under the standing assumption Assumption ($X$), the map $\phi \mapsto \I  ^{*}\phi :=\phi \circ T%
\I  $ defined on measurable geometric forms on $\C_{x_0}M$ extends to
a continuous linear injective map
$$\I  ^{*}:L^0\Gamma \HH^{*}\quad
\longrightarrow\quad  L^0\left( \C_0\R
^m;H^{*}\right)
$$
from measurable $\HH $-one forms on $\C_{x_0}M$ to measurable $H$-one -forms on $\C_0\R^m$, using the topology of convergence in
probability. The map is given by the It\^{o} stochastic integral
\begin{equation}
\label{I-phi}
\I  ^{*}(\phi )(h)=\int_0^T\left\langle
 \frac{{\D }\phi _s^{\#}}{
d s},\nabla_{T\I  _s(h)}X (\tilde{\parals_s} d\beta
_s)+X(x_s)(\dot h_{s})\,ds \right\rangle_{x_s},      \qquad  h\in H
\end{equation}
using the filtration $\G_t:=\F_t^\beta \vee \F^{x_0}$, 0$\le t\le T$.
Moreover for $1\le p<\infty$,
\begin{enumerate}
\item[(a)]
 the map $\I^*$ restricts to a continuous linear map
\begin{equation}
\I  ^{*}:L^{p+\varepsilon }\Gamma \HH ^{*}\longrightarrow
L^p\left( \C_0\R ^m;H^{*}\right)  
 \label{I-star}
\end{equation}
for any $\varepsilon >0$. If condition ($M_0$) holds,
(\ref{I-star}) holds for $\varepsilon =0$. 
\item[(b)]
 if $\phi\in L^0\Gamma\HH$
\begin{equation}
\E \left\{ \I  ^{*}(\phi(-) )\left| \F ^{x_0}\right.\right\} =\phi(\overline{T\I}_{x_\cdot}-)
\end{equation}
and 
\begin{equation}
\label{I-star-domination}
\|\phi\|_{L^p}\le \|\I^*(\phi)\|_{L^p}.
\end{equation}
\end{enumerate}
Consequently  $\I^*[L^0\Gamma\HH^*]\cap L^p(\C_0\R^m;H)$ is closed in $L^p$ and is contained in
 $\I^*[L^p\Gamma\HH^*]$ with equality if condition ($M_0$) holds.
\end{theorem}

\begin{proof}
(1) For $h\in H$ set $v_t=T\I_t(h)$.  From (\ref{eq-filter-2} ),
\begin{equation}
\label{derivative-Ito-map}
v_t=W_t\int_0^t W_s^{-1} \nabla _{v_s} X(\tilde{\parals_s} d\beta_s)
+\W_t\left(X(x_\cdot)(\dot h_\cdot)\right).
\end{equation}
Consider $\phi(v_\cdot)$ for $\phi$ a geometric 1-form. We can treat $\phi$ as an $\HH$-form by restriction.
 As such it has a dual $\HH$-vector field $\phi^{\#}$,
 so if $u\in \HH_\sigma$ then
$$\phi(u)=\phi_{\sigma}(u)=\int_0^T \left \langle
 {\D\over dt}u_t,
{\D\over dt}\phi_t^{\#}(\sigma)
 \right \rangle_{\sigma(t)}\; dt.$$
Thus
 \begin{equation}
\label{pull-back-map-2}
\phi(v_\cdot)=\phi\left(W_t\int_0^t W_s^{-1} \nabla _{v_s} X(\tilde{\parals_s} d\beta_s)
 \right)+\int_0^T  \left\langle  {\D \over dt}\phi^{\#}_t,
X(x_t)\dot h_t   \right\rangle \, dt 
\end{equation}
 For the second term on the right hand side, we have
\begin{eqnarray*}
\left|\int_0^T  \left\langle  {\D \over dt}\phi^{\#}_t,
X(x_t)\dot h_t   \right\rangle \, dt\right|
&\le&\|\phi^{\#}\|_{\HH_{x_\cdot}} \|h\|_H, \hskip 10pt \hbox{a.s.}
\end{eqnarray*}
Consequently the map
 $h\mapsto \phi_{x_\cdot}(\W_\cdot(X(x_\cdot)(\dot h_\cdot))$
is in $H^*$ almost surely.   For $p=0$ or $1\le p\le \infty$ it gives an
element of $L^p(\C_0\R^m; H^*)$ depending continuously on the restriction
of $\phi$ in $L^p\Gamma \HH^*$.

(2)  For the more interesting first term on the right hand side of (\ref{pull-back-map-2}),
\begin{equation}\Phi(h,\phi)\equiv \phi\left(W_\cdot \int_0^\cdot W_s^{-1} \nabla_{v_s}
X(\tilde{\parals_s} d\beta_s)\right),
\end{equation}
we assume that  $\phi$ is a geometric differential  form
 on $\C_{x_0}M$ which extends to give linear functionals on the $L^2$ tangent spaces
  $L^2T_\sigma \C_{x_0}M$, for some choice of a Riemannian metric on $TM$ extending that of $E$.  
Then there is a section $\alpha$ of the vector bundle
 $L^2T\C_{x_0}M\to \C_{x_0}M$ such that if $u\in L^2 T_\sigma \C_{x_0}M$ then
\begin{equation}
\label{form-L2-representation}
\phi(u)=\int_0^T\langle \alpha(\sigma)_t, u_t\rangle_{\sigma(t)}\;dt.
\end{equation}
Assume $\alpha$ is in $L^2$, i.e.
 $\int_{\C_{x_0}M} \|\alpha(\sigma)\|_{L^2}^2 \, d\mu_{x_0}(\sigma)<\infty$.
 We first  show that  for such $\phi$ the right hand side of 
 (\ref{I-phi}) makes sense and agrees with the pull back $\phi(v_\cdot)$.
It is easy to verify that
\begin{equation}
\label{TH-2.4-3}
\phi^{\#}_t
=\W_t\left(\Pi(W_\cdot^{-1})^*\int_\cdot^TW_s^*\alpha_s\,ds\right),
\end{equation}
where $W_r^*: T_{x_r}M\to T_{x_0}M$ is the adjoint of $W_r$.
 Let $e^1,\dots ,e^{n-p}$
be an orthonormal base for ker[$X(x_0)$].
 Set $e_s^j=\tilde {\parals_s} e^j$ and
 $\beta _s^j=\left\langle \beta _s,e^j\right\rangle _{\R ^m}$.
Then
\begin{eqnarray*}
\Phi (h,\phi ) &=&\sum_j\phi \left( W_{\cdot }\int_0^T\chi _{[0,\cdot
]}(s)W_s^{-1}\nabla_{v_s}X(e_s^j) \; d\beta _s^j\right) \\
&=&\sum_j\int_0^T\phi \left( \chi _{[0,\cdot ]}(s)W_{\cdot }^s \nabla
_{v_s}X(e_s^j)\right) d\beta _s^j
\end{eqnarray*}
where $W_t^s=W_t^{-1}W_s$, since $\phi_{x_\cdot} $ is $\F ^{x_0}=\mathcal{G}_0$
measurable. Using (\ref{form-L2-representation}) this gives
\begin{eqnarray*}
\Phi (h,\phi ) &=&\sum_j\int_0^T\left\{ \int_s^T\left\langle \alpha _t,W_t^s%
\nabla_{v_s}X(e_s^j)\,dt\right\rangle \right\} d\beta _s^j \\
&=&\int_0^T\left\langle \int_s^T\left( W_t^s\right) ^{*}\alpha
_t\,dt, \nabla_{v_s}X\left( \tilde{\parals_s} d\beta
_s\right) \,
\right\rangle _{x_s} \\
&=&\int_0^T\left\langle { \D \over d s}\phi _s^{\#},
\nabla_{v_s}X\left( \tilde{\parals_s} d\beta _s\right)
 \,\right\rangle _{x_s}.
\end{eqnarray*} 
This shows that $\phi(v_\cdot)=\I^*(\phi)$ agrees with the right hand side of (\ref{I-phi}), as expected.
\medskip

(3)  In general for $\psi $ a measurable $\HH $-form define
$$\tilde{\Phi}(\psi )_t:=\int_0^t\left\langle
{ \D \over d s}\psi
_s^{\#}, \nabla_{T\I  _s(-)}X\left(
\tilde{\parals_s} d\beta _s\right) \,\right\rangle _{x_s},\hskip
6pt 0\le t\le T.
$$

Note that $(\tilde\Phi(\psi)_t, 0\le t\le T) $  is a local $\mathcal{G}_{*}$-martingale with values in
$H^{*}$ since

$$\int_0^T\sup_{\left\| h\right\| _H\le 1}\left| \left( \nabla_{T%
\I  _s(h)}X\right) ^{*}\left( \frac{\D }{d s}\psi_s^{\#}\right) \right| _{x_s}^2\,ds<\infty ,
\quad \as
$$
using the fact that
$$
\sup_{0\le s\le T}\sup_{\left\| h\right\| _H\le 1}\left| \left( \nabla_{T\I  _s(h)}X\right) \right| _{\LL(\R ^m;E_{x_s})}^2<\infty
\quad \as .$$
 The usual stopping time argument, see Lemma \ref{stopping-times}, 
  shows that  $\tilde\Phi(\psi)(-)$ is a continuous map from $L^0\Gamma \HH^*$ to $L^0(\C_0\R^m;H^*)$.

To see that  the pull back map of a geometric differential 1-form $\phi$ evaluated at $h$ agrees with the
 right hand side of (\ref{I-phi}), we only need to show that $\tilde\Phi(\phi)(h)=\Phi(h,\phi)$. 
For this define a sequence of differential forms $\phi_n$, which extends over $L^2TC_{x_0}M$, by
$$\phi_n(\sigma)(u)=\phi\left(\parals_\cdot (\sigma)\int_0^T \lambda_n(s-\cdot)\parals_s^{-1}(\sigma) u_s ds\right)$$
for suitable $\lambda_n: \R^n \to \R^n$, $n=1,2,\dots$ so that $\phi_n$ converges to $\phi$ on $(T_\sigma C_{x_0}M)^*$ for almost all $\sigma$, and observe that $\tilde\Phi(\phi_n)(h)\to \tilde\Phi(\phi)(h)$ in particular from the convergence of $\phi_n$ in $\HH^*$ and Lemma \ref{stopping-times}.

Furthermore by the Burkholder-Davis-Gundy inequalities, \cite
{Daprato-Zabczyk2002}, for $0<p<\infty $%
\begin{eqnarray}
\E \Big| \tilde{\Phi}(\psi )\Big| ^p &\le &C_p\E \Big| \int_0^T\left| \left( \nabla_{T\I  _s(-)}X\right) ^{*}\left(\frac{\D }{d s}\psi _s^{\#}\right) \right| _{\LL_2(H; \R^m)}^2\,ds\Big| ^{p/2}  \nonumber \\
&\le &C_p\E  \Big\{ \sup_{0\le s\le T} 
\left| \left( \nabla_{T \I  _s(-)}X\right) ^{*}\right|
 _{\LL_2(H; \LL(E_{x_s}; \R ^m))}^p\,\left\| \psi \right\| _{\HH _{\cdot }^{*}}^p \Big\}.
\label{continuity-1}
\end{eqnarray}
For $1\le p<\infty $ this gives by H\"older's inequality, for
 $\epsilon>0$,
$$\left| \tilde{\Phi}(\psi )\right| _{L^{p}}
\le \hbox{const.}\,\n{\psi}_{L^{p+\epsilon}}\cdot \n{\sup_{0\le s\le T}\left| T\I  _s(-)\right|
_{\LL(H; T_{x_s}M)}} _{L^{\frac{p(p+\epsilon )}\epsilon }}.$$
Since $\sup_{0\le s\le T}\left| T\I  _s(-)\right| $ lies in $L^q$ for
all $1\le q<\infty $ we see that $\tilde{\Phi}$ gives a continuous linear
map
$$\tilde{\phi}:L^{p+\epsilon }\Gamma \HH ^{*}\longrightarrow L^p\left(
\C_0\R ^m;H^{*}\right)$$
for all $1\le p<\infty $ and $\epsilon >0$.

Combining the two terms in (\ref{pull-back-map-2})
we obtain a continuous linear map
$$
\psi \mapsto \int_0^T\left\langle \frac{\D }{d t}\psi _t^{\#},%
\nabla_{T\I  _t(-)}X\left( \tilde{\parals_t} d\beta
_t\right) +X(x_t)(\frac d{dt}-)\right\rangle _{x_t}\,dt
$$
from $L^{p+\epsilon }\Gamma \HH ^{*}$ to $L^p(\C_0\R ^m;H^{*})$ in
the relevant range of $1\le p<\infty ,\epsilon > 0$, which agrees with
$$
\psi \mapsto \I  ^{*}(\psi )
$$
when $\psi $ is an $L^\infty $ section of $(L^2T\C_{x_0}M)^{*}$.

Supposing furthermore that condition ($M_0$) holds, observe that (\ref{continuity-1}) gives
$$\n{ \tilde{\Phi}(\psi )} _{L^p}\le (C_p)^{1/p}\,\n{ \psi}_{L^p}  \cdot
\n{ \E \Big\{ \sup_{0\le s\le T}\big|\nabla_{T\I  _s(-)}X\big| ^p\,\Big|
 \F ^{x_0}   \Big\}} _{L^\infty}
$$
with
\begin{eqnarray*}
&&\E \Big\{ \sup_{0\le s\le T}\big| {\nabla}_{T\I  %
_s(-)}X\big| ^p\,\Big| \F ^{x_0} \Big\}\\
 &&\le  \sup_{0\le s\le T}\left| \nabla_{W_s(-)}X
    \right| _{\LL\left(T_{x_0M}; \LL(\R^m; E_{x_s})\right) }^p 
 \E \Big\{ \sup_{0\le s\le T}\big| W_s^{-1}T\I  _s(-)\big|
_{\LL_2(H; T_{x_0M)}}^p\,\Big| \F ^{x_0} \Big\}
\end{eqnarray*}
which is essentially bounded by condition ($M_0)$ and Proposition
\ref{proposition-3.2}. Thus in this case we can take $\epsilon =0$
and (a) is proved.

To see that $\I^*: L^0\to L^0$ is injective  note that by (\ref{I-phi})
and the independence of $\beta$ and $\F^{x_0}$,
$$
\E \left\{ \I  ^{*}(\phi(-) )\left| \F ^{x_0}\right.
\right\} =\int_0^T\left\langle \frac{\D }{d s}\phi
_s^{\#},\,X(x_s)(\frac d{ds}-)\right\rangle _{x_s}\,ds
=\phi(\overline{T\I}_{x_\cdot}-)
$$
(with a suitable interpretation of the conditional expectation if
 $\I  ^{*}(\phi )$ is not in $L^1$, e.g. see
 Elworthy-LeJan-Li \cite{Elworthy-LeJan-Li-book} p66).
Using the inverse of $\overline{T\I}$, see Property \ref{property-1} we obtain
\begin{equation}
 \label{continuity-2}
\phi _{x_\cdot}(-)=\E \left\{ \left.\I  ^{*}(\phi )\left( \int_0^{\cdot
}Y_{x_{s\cdot }}\left( \frac{\D }{d s}-\right) ds\right) \right|
\F ^{x_0}\right\}. 
\end{equation}
This proves injectivity, giving a left inverse for
$\I^*$.

Moreover since $v\mapsto \int_0^\cdot Y_{x_s}\left({\D\over ds}v_s\right)ds$
is an isometry of $\HH_{x_\cdot}\to H$ almost surely we see 
\begin{eqnarray*}
\|\phi\|_{L^p}&=&\| \E\left\{ \I^*(\phi)| \F^{x_0} \right\} \|_{L^p}\\
&\le & \|\I^*(\phi)\|_{L^p}
\end{eqnarray*}
completing the proof.
 \end{proof}

\begin{remark}
\label{re-3.6}
Although the term ${ \D \over ds}  \phi^{^{\#}}_s$ appearing  in (\ref{I-phi}) 
 may depend on the whole path $\{ x_s : 0 \le s\le T \}$,  the stochastic integral 
 there can and was considered as an It\^o integral by regarding $\beta_\cdot$ as a martingale with respect to $\F^{x_0}\cup \F_t$. We can also 
treat it as a Skorohod integral, c.f. X-D Li \cite{X-D-Li-2003} for pull backs by the stochastic development map. In fact if $\phi\in L^{p+\epsilon}\Gamma\HH$ the stochastic integral can be written as
  $$\int_0^T \left\langle\nabla_{T\I_s ( - )}Y\left(
{\D \over ds}  \phi^{\#}_s \circ \I\right), dB_s\right \rangle$$
 and interpreted as a Skorohod integral on $\C_0 \R^m$, i.e. as 
  $$ ( d^q )^{*}\left  ( \alpha \mapsto \int_0^T   \chi_{[0,\cdot]}(s) 
  \left\langle \nabla_{T\mathcal{I}_s ( - )} Y (
  {\D \over  ds}  \phi^{\#}_s \circ \I), {d\over ds} \alpha_s\right\rangle_{\R^m} ds\right) $$
   for $\frac{1}{p} + \frac{1}{q} = 1$ when $1 < p < \infty$. Here 
   $$d^q\equiv d^q_{H^*} : \Dom ( d^q ) \subset L^q (\C_0 \R^m ; H^{\ast} ) \longrightarrow 
   L^q (\C_0 \R^m ; \LL_2 ( H; H^{\ast} ) ).$$
   \end{remark}
   
   \begin{proof}
Set $G=\LL_2\left(H; L_0^{2,1}\left([0,T]; T_{x_0}M\right)\right)$. Define
  $\theta: \Omega\to  \LL_2(H;   G)$
 by
$$\theta(\alpha)(h)_\tau
=\int_0^T\chi_{[0,\tau]}(s) W_s^{-1} \nabla_{T\I_s(h)}X(\dot \alpha_s) ds, \quad \alpha, h \in H, 0 \le \tau \le T.$$
Then $$(d^q_G)^*(\theta)(-)_\tau
=\int_0^T\chi_{[0,\tau]}(s) W_s^{-1} \nabla_{T\I_s(-)}X(dB_s)
\in \LL(H;T_{x_0}M).$$

 Suppose first that $\phi$ is smooth and cylindrical. Then 
$$\Phi(-, \phi)= \phi_x \left( W_\cdot \int^T_0 \chi_{[ 0, \cdot]} ( s )
     W_s^{- 1} \nabla_{T\I_s ( - )} X ( dB_s ) \right)
          = \phi_x \left( W_\cdot ( d^q_G )^* ( \theta) (-)_\cdot\right).$$

   Now, if $\{ E^k \}_{k = 1}^{\infty}$ is an orthonormal base for $H$, $g: \C_0 \R^m \to \R$ is $C^{\infty}$ smooth cylindrical, and $h\in H$,
  \begin{eqnarray*}
&&  \E \left[ g \phi_{x_\cdot}\left( W_\cdot  (d^q_G )^* \theta (-)_\cdot \right)(h)\right]
= \E \left[ g \phi_{x_\cdot}\left ( W_\cdot  (d^q _{L_0^{2,1}([0,T];T_{x_0}M)})^* (\theta (-)(h)) \right)_\cdot \right]\\
&&  =\E \left[\sum_j d^q \left( g \phi_{x_\cdot} ( W_\cdot- ) \right) ( E^j ) ( \theta( E^j  )(h)_\cdot)\right] \\
    && =\E \left[\sum_j  d^qg ( E^j ) \phi_{x_\cdot} ( W_\cdot\theta( E^j )(h)_\cdot ) \right]
    +\E\left[ \sum_j g  d^q( \phi_{x_\cdot} ( W_\cdot- ) ) ( E^j) \left( \theta( E^j )(h) _\cdot \right) \right]\\
      && =\E \left\langle d^q g, \phi_{x_\cdot}(W_\cdot\theta(-)(h)_\cdot\right\rangle_{H^*}
    +\E\left[ \sum_j g  d^q( \phi_{x_\cdot} ( W_\cdot- ) ) ( E^j) \left( \theta( E^j ) (h)_\cdot \right) \right].
      \end{eqnarray*}
Furthermore
\begin{eqnarray*}
&&\sum_j d ( \phi_{x_\cdot} ( W_\cdot - ) ) ( E^j ) (\theta ( E^j )(h)_\cdot )\\
 &&= \sum_j d ( \phi_{x_\cdot} ( W_\cdot- ) )
\left(\int_0^\cdot K^\perp(x_s) \dot E^j _sds\right) ( \theta (\int_0^\cdot K(x_s)E_s^jds )(h) ), 
\end{eqnarray*}
since $T\mathcal{I}( E^j ) = T\I(\int_0^T K^\perp(x_s) \dot E^j _sds)$ and $\nabla_{T\I_s ( \cdot )} X ( \dot{E^j_s }) 
= \nabla_{T\I_s( \cdot)} X ( K ( x_s ) ( \dot E^j_s )  )$ by the defining property of $\nabla$, 
Proposition \ref{theorem-connection}.
 Note that the expression is independent of the choice of basis and  so vanishes giving
 \begin{equation}
 \Phi(h,\phi)
 =\phi_x \left( W_\cdot ( d^q _G)^*( \theta(-)_\cdot (h) )\right) 
= \left(d^q _{H^*}\right)^*      \left( \phi_{x_\cdot} \left( W_\cdot\theta  (-)(h)_\cdot\right)\right).
   \end{equation}
Note that for $h\in H$,
$$\phi_{x_\cdot} ( W_\cdot\theta  ( - )(h))_t^{\#}=\int_0^t
 \nabla_{T\I_s (h)}Y ( {\D \over ds} \phi^{\#}_s) ds\in H, $$
to obtain the desired result.

For general $\phi \in L^{p + \epsilon}\Gamma \HH$ we can take
$C^{\infty}$ cylindrical one-forms $\phi^j$, $j = 1$ to $\infty,$
converging to $\phi$ in $L^{p + \epsilon}\Gamma \HH$. By the Theorem,
$\I^* ( \phi^j ) \to \I^*( \phi)$ in $L^p$ and we see therefore that 
$\Phi ( -, \phi^j )$ is convergent in $L^p$. Thus 
$\nabla_{T\mathcal{I}_s ( - )} Y_{-} ({ \D\over ds} \phi^{\#}_s)$ is in
the domain of $( d^q )^*$ and the result holds.
\end{proof}

\begin{remark}
From  (\ref{I-phi}) we see that $\I^*(\lambda \phi)=\lambda \I^*(\phi)$ for all $\phi\in L^0\Gamma\HH^*$, 
if $\lambda\in L^0(\C_{x_0}M;  \R)$. This is because $\lambda\circ \I\in \G_0$ for all
such $\lambda$. From the Skorohod integral representation, 
$$\I^*(\phi)=\int_0^T \left \langle \nabla_{T\I_s ( - )}Y\left(
{\D \over ds}  \phi^{\#}_s \circ \I\right), dB_s\right \rangle+\int_0^T  \left\langle  {\D \over ds}\phi^{\#}_s,
X(x_s)\dot h_s   \right\rangle \, ds, 
$$
this is less obvious. However for sufficiently regular $\lambda$,
\begin{eqnarray*}
&&\int_0^T \left \langle \nabla_{T\I_s ( - )}Y\left(
{\D \over ds} \left(\lambda \phi^{\#}_s\right) \circ \I\right), dB_s\right \rangle\\
&=&\lambda \int_0^T \left \langle \nabla_{T\I_s ( - )}Y\left(
{\D \over ds}  \phi^{\#}_s \circ \I\right), dB_s\right \rangle\\
&&-\int_0^T \left \langle \nabla_{T\I_s ( - )}Y\left(
{\D \over ds}  \phi^{\#}_s \circ \I\right), {d\over ds} \nabla(\lambda\circ \I)\right \rangle_{x_s},
\end{eqnarray*}
and as for the proof of Remark  \ref{re-3.6} the second term vanishes. 
{\hspace*{\fill}
  \begin{math}\Box\end{math}\medskip}
\end{remark}

For suitable $h:\C_0\R ^m\rightarrow H$ define a measurable vector
 field $\overline{T\I  (h)}$ on $\C_{x_0}M$ by
$$
\overline{T\I  (h)}(\sigma )=\E \left\{ T\I  
(h(\cdot ))\left| \I  (\cdot )=\sigma \right. \right\}
$$
for $\mu _{x_0}$-almost all $\sigma $ in $\C_{x_0}M$. 
Note that if $h$ is $\F^{x_0}$ measurable with $h=\bar h\circ \I$ then
$\overline{T\I(h)}(\sigma)=\overline{T\I}_\sigma(\bar h(\sigma))$ for $\overline{T\I}$
as in Property \ref{property-1}. For completeness we
give the following extension of Theorem 2.2 of Elworthy-Li \cite{Elworthy-Li-Hodge-1}:

\begin{corollary}
\label{corollary-pull-back}
For $\delta > 0$ and $1<q<\infty $ the map
 $h\mapsto \overline{T\I  (h)}$ gives a continuous linear map
$$\overline{T\I  (-)}\colon
  L^q(\C_0\R ^m;H)\rightarrow L^{q-\delta}\Gamma \HH .$$
(This map is the co-joint of $\I  ^{*}$ in the sense that
\begin{equation}
\int_{\C_0\R ^m}\I  ^{*}(\phi )(h)dP
=\int_{\C_0\R ^m}\phi \left(
\overline{T\I  (h)}\right)dP   \label{continuity-3}
\end{equation}
for $\phi \in L^{p+\epsilon }\Gamma\HH^*$ and $h\in L^q(\C_0\R^m;H)$,
 taking $\epsilon >0$ and $\frac 1p+\frac 1q=1$.)

If Condition $(M_0)$ holds we can allow $\ \delta =0$ and also
$q=\infty $, and
$\overline{T\I  (-)}:L^q(\C_0\R ^m;H)\rightarrow
L^q\Gamma \HH $  is surjective, $1<q\le \infty $.
\end{corollary}

\begin{proof} 
First note that when $\phi $ is a true form
$$\int_{\C_0\R ^m}\I  ^{*}(\phi )(h)dP
=\E \left\{ \phi \left( \E \left\{ T\I  (h)
 \left| \F ^{x_0}\right. \right\} \right) \right\}
 =\E \phi \left( \overline{T\I  (h)}\right)$$
so (\ref{continuity-3}) holds and it holds for 
 $\phi \in L^{p+\epsilon }\Gamma \HH ^{*}$,
 $h\in L^q(\C_0\R ^m;H)$ for $\frac 1p+\frac 1q=1$  by continuity, 
from Theorem \ref{th:I-star}. Consequently for such
$\phi$ and $h$,
$$
\Big| \int_{\C_0\R^m}\phi \left( \overline{T\I  (h)}
            \right)dP \Big|
 \le \n{\I^{*}(\phi )} _{L^p}\cdot \n{h}_{L^q}
\le  \mathrm{const.}\cdot\n{\phi}_{L^{p+\epsilon }}\,\n{h}_{L^q}.
$$
 Since this holds for all $\phi $ in
$\displaystyle{L^{p+\epsilon }}$ we see
that $\displaystyle{\overline{T\I  (h)}
                 \in L^{\frac{p+\epsilon }{p+\epsilon -1}}}$
and is continuous linear in $h$ into
 $\displaystyle{L^{\frac{p+\epsilon }{p+\epsilon -1}}}$.
 Thus if $h\in L^q$ then $\displaystyle{\overline{T\I  (h)}
 \in L^{\frac{p+\epsilon }{p+\epsilon -1}}}$
 for all $\epsilon>0$ so that $\displaystyle{\overline{T\I  (-)}}$
 is continuous  linear from $L^q$ to $L^{q-\delta }$ for any $\delta >0$,
 with $\delta =0$ allowed if condition $(M_0)$ holds.

Surjectivity comes from the fact that $\overline{T\I  (-)}$ has a right inverse
\begin{equation}
v\mapsto \int_0^{\cdot }Y_{x_s}\left( \frac{\D }{ds}v_s\circ \I\right) ds
\label{continuity-4}
\end{equation}
mapping $L^q\Gamma \HH $ to $L^q\left( \C_0\R ^m;H\right) $ 
as in Property \ref{property-1}. 
\end{proof}

Analogously to (\ref{TH-2.4-3}) we have also:
 \begin{remark}
Let $N: [0,T]\times \C_{x_0}M\to \R^p$ be a continuous semi-martingale on the filtered 
probability space $\{\C_{x_0}M, \mu_{x_0}, \F^{x_0}_*\}$ where 
$\F^{x_0}_t=\sigma\{x_s: 0\le s\le t\}$. Suppose $\alpha_\cdot$ is a locally bounded
 section of $\LL(\R^p;L^2T\C_{x_0}M)$ which is $\F^{x_0}_*$ adapted. Then the mapping
 \begin{eqnarray*}
 L^2\HH&\to& L^0(\C_{x_0}M;\R)\\
 U&\mapsto&\int_0^T\left\langle U(\sigma), \alpha(\sigma)_s dN_s(\sigma)\right\rangle_{\sigma(s)}
 \end{eqnarray*}
 can be considered as the map $U\mapsto \varphi(U(\cdot))$ for $\varphi$ the $H$-one-form determined by the $\HH$-vector field $\varphi^{\#}$ with
 $$\varphi_t^{\#}
 =\W_t \left(\Pi(W_\cdot^{-1})^*\int_\cdot^T W_s^* \alpha_s dN_s\right), \qquad 0\le t\le T.$$
 If the martingale part of $N_\cdot$ is zero then no adapteness is required.
 \end{remark}

\begin{remark} 
The proof of Theorem 2.2 in Elworthy-Li \cite{Elworthy-Li-Hodge-1} as it stands has a trivial
 mistake in the last line and only gives $\overline{T\I(h)}\in L^1$ when
 $h\in L^2$. However it is easily modified by not taking expectations in
the proof of Lemma 3.3 of \cite{Elworthy-Li-Hodge-1}: this shows that $\overline{T\I(h)}\in L^\infty$ if
$h\in L^2(\C_0\R^m,\sigma\{\beta_s: 0\le s\le T\};H)$, not just $L^2$,
using the observation that
$\sup_{0\le r\le T}\E\{|T\I_r|^2 \big| \F^{x_0}\}<\infty$
from Aida-Elworthy \cite{Aida-Elworthy}, i.e. a special case of
 Proposition \ref{proposition-3.2} above.
\end{remark}

\section{Sobolev calculus on $\C_{x_0}M$ and its intertwining by It\^o maps}
\label{section-4}

If $M$ is given a Riemannian metric then $\C_{x_0}M$ gets a Finsler structure
defined by
$$\|v\|_{\sigma}:=\sup_{0\le t\le T} |v_t|_{\sigma(t)},
\hskip 20pt v\in T_\sigma \C_{x_0}M.$$
By the compactness of $M$ different metrics produce uniformly equivalent
Finsler norms. From Elworthy-LeJan-Li \cite{Elworthy-LeJan-Li-book} and 
Elworthy-Ma \cite{Elworthy-Ma}
and the compactness of $M$, for almost all $\sigma$ the inclusion
$\HH_\sigma\hookrightarrow T_\sigma \C_{x_0}M$ is continuous and its norm
is in $L^p$ as a function of $\sigma$ for all $p\in [1,\infty)$ and is essentially bounded for any of these Finsler norms if condition (M) holds. 

 A function $f: \C_{x_0}M\to \R$ will be said to be $BC^1$
if it is Fr\'echet differentiable and is bounded together with its
 differential $df$ considered as a section of $T^*\C_{x_0}M$, again
using any of these Finsler structures. Note that then $df$ restricts to
$\HH$ to give an element of $L^p\Gamma\HH^*$ for all $1< p<\infty$ and lies in $L^\infty\Gamma\HH^*$ if condition $M$ holds, by
the definition of $\HH$ and Proposition \ref{norm-TI}.

Denote by $\Cyl$ the space of smooth 
cylindrical functions on $\C_{x_0}M$. 
Let $\Dom(d_\HH)$ be a linear subspace of $L^\infty(\C_{x_0}M;\R)$ with
\begin{equation}\label{domain-1}
\Cyl\subset \Dom(d_\HH)\subset BC^1.\end{equation}
 Define
$$d=d_\HH:  \Dom(d_\HH)\to \cap_{1\le p<\infty} L^p\Gamma \HH^*$$
 by restriction:
$$(d_\HH f)_\sigma=\left. df_\sigma \right\vert_{ \HH_\sigma}.$$

Let $\dom(d_\HH)$ be the space of equivalence classes of
 $\Dom(d_\HH)$ under equality up to sets of $\mu_{x_0}$-measure zero. However
  after this section we will not
 distinguish between $\Dom(d)$ and $\dom(d)$.
 By a standard result the set of smooth cylindrical functions is dense in $L^p(\C_{x_0}M; \R)$, $1\le p<\infty$,
 and therefore so is $\dom(d_\HH)$.

We next give the proof of closability of $d_\HH$ restricted to
 $\dom(d_\HH)$ in our context. For classical
Wiener space this is one of the basic results of Malliavin calculus
and gives for each $1\le p<\infty$ a closed linear operator
$$d\equiv d^p: \dom(d^p)\subset L^p(\C_0\R^m;\R)\to L^p(\C_0\R^m; H^*).$$
This is proved by the standard method of integration by parts in
 Nualart \cite{Nualart-book}. (Strictly speaking his basic domain does not consist of $BC^1$ functions but it is easy to see, and well
known, that this gives the same closures $d^p$.)  On the path space $d_{\HH}$ is closable,
from Elworthy-LeJan-Li \cite{Elworthy-LeJan-Li-book},  in $L^2(\C_{x_0}M; \R)$. For $p>1$ the
following theorem follows in the same way from the integration by
parts results in Elworthy-LeJan-Li \cite{Elworthy-LeJan-Li-CR},
\cite{Elworthy-LeJan-Li-book}.   Related results in this context and a
 detailed discussion of the Dirichlet forms which arise  can be found
in Elworthy-Ma \cite{Elworthy-Ma} and \cite{Elworthy-Ma-Kiev}.
See also Driver \cite{Driver92} and  Hsu \cite{Hsu-Quasi-invariance}.

\begin{theorem}
\label{theorem-closability}
For $1\le p<\infty$ the operator $d_\HH$ can be considered as a linear
operator
$$d_\HH: \dom(d_\HH)\subset L^p(\C_{x_0}M;\R)\to L^p\Gamma\HH^*.$$
It is closable for each $1<p<\infty$ and for $p=1$ if condition
$(M_0)$ holds.
\end{theorem}

\begin{proof}
It is not immediately obvious that $\displaystyle{d_\HH}$ is defined on
$\dom(d_\HH)$: we need to know that  if $\displaystyle{f\in \Dom(d_\HH)}$ is
 $\displaystyle{\mu_{x_0}}$-almost surely $0$
then $\displaystyle{d_\HH f=0}$ almost surely (c.f. Proposition 3.5 in
Elworthy-Ma  \cite{Elworthy-Ma}).
This comes together with closability from the following proof
that if
$\displaystyle{f_j \in \dom(d_\HH)}$  has  $\displaystyle{f_j\to 0}$ in
 $L^p$ and  $\displaystyle{d_\HH f_j\to \theta}$ in
 $\displaystyle{L^p\Gamma\HH^*}$   then $\theta=0$.
For this note that $f_j\circ \I\to 0$ in $L^p(\C_0\R^m;\R)$ and
$f_j\circ \I$ belongs to the domain of $d$, as is well known, for example
by Wong-Zakai approximations. Moreover
$$d^p(f_j\circ \I)=d_\HH f_j\circ T\I=\I^*(d_\HH f_j).$$
By Theorem \ref{th:I-star}, $\I^*(d_\HH f_j)$ converges to $\I^*(\theta)$
in $L^{p^\prime}$ for any $p^\prime \in [1,p)$ and for
 $p^\prime=p=1$ if condition
$(M_0)$ holds. Since $d^{p^\prime}$ for $\C_0\R^m$ is closed this shows
 that $\I^*(\theta)=0$. But Theorem  \ref{th:I-star} shows $\I^*$
is injective and so $\theta=0$ as required.
\end{proof}

\bigskip

Let $\displaystyle{d^p: \dom(d^p)\subset L^p(\C_{x_0}M;\R)
\to L^p\Gamma\HH^*}$, $1\le p<\infty$, be the closures of
 $d_\HH$ given by Theorem \ref{theorem-closability}.

We now comes to our main result on the chain rule, or intertwining.
For $1\le p^\prime\le p<\infty$ let
$$\I^*\colon  L^p(\C_{x_0}M;\R) \to L^p(\C_0\R^m;\R)
\hookrightarrow L^{p^\prime}(\C_0\R^m;\R)$$
 denote the map $\I^*f=f\circ \I$ as well as the map from $L^p\Gamma\HH^*$ to $L^{p'}(\C_0\R^m;H^*)$,
 as defined by Theorem \ref{th:I-star}.  We will see  below
that the It\^o map for our s.d.e.  can be used in some sense
as a substitute for a chart for the `differential structure'
given by the calculus. See also Aida-Elworthy \cite{Aida-Elworthy}, Aida \cite{Aida97},
Elworthy-Li- \cite{Elworthy-Li-Hodge-1} for the gradient case,
Elworthy-LeJan-Li \cite{Elworthy-LeJan-Li-book}, Elworthy-Li \cite{Elworthy-Li-vector-fields}
more generally, and for related work see Fang-Franchi \cite{Fang-Franchi-97}, \cite{X-D-Li-2003} and Cruzeiro-Malliavin \cite{Cruzeiro-Malliavin-96}.

\begin{theorem}
\label{theorem-closability-2}
Suppose $1\le p^\prime\le p< \infty$, with $1<p^\prime< p$ unless
condition ($M_0$) holds.  The operators
$$\I^*d^p:  \dom(d^p)\subset L^p(\C_{x_0}M;\R)
 \longrightarrow L^{p^\prime}(\C_0\R^m;H)$$
and 
$$d^{p^\prime} \I^* :  \big\{f\in L^p \; \big\vert  \;
\I^*f\in \dom(d^{p^\prime})\big\} \subset  L^p(\C_{x_0}M; \R)
\longrightarrow  L^{p^\prime}(\C_0\R^m;H)$$
 are densely defined. Moreover
\begin{enumerate}
\item[(i)]
$\I^* d^p$ is closable in general and closed if condition $(M_0)$
holds and
 $p=p^\prime$;
\item[(ii)]
$d^{p^\prime }\I^*$ is closed.
\item[(iii)]
\begin{equation}
\label{domain-d}
\I^*d^p\subset d^{p^\prime}\I^*.
\end{equation}
\end{enumerate}
\end{theorem}
\begin{proof} 
The denseness of the domain of $d^{p^\prime}\I^*$ and the fact that it is closed 
are automatic by continuity of $\I^*$, giving (ii). It is clear that $\I^*d^p$ is densely defined. 

 For (i) first suppose that condition $(M_0)$ holds and $p=p^\prime$. Then
$\I^*$ on $L^p$ $\HH$-forms is continuous and has closed range, 
by Theorem \ref{th:I-star}. It follows that its composition $\I^*
d^p$ with the closed operator $d^p$ is closed. In the general case suppose
$1\le p^\prime <p$ and
 $\{f_j\}_{j=1}^\infty$ is a sequence in $\dom(d^p)$ with $f_j\to 0$ in
$L^p(\C_{x_0}M;\R)$ and $\I^*(d^p f_j)\to \theta$ in
$L^{p^\prime}(\C_0\R^m; H^*)$, some $\theta$.
By Theorem \ref{th:I-star}(b) we know $\theta=\I^*(\alpha)$ for
 some $\alpha \in L^{p'}\Gamma\HH$ with $d^pf_j\to \alpha$ in $L^{p'}$.
  Since
$f_j\to 0$ in $L^{p^\prime}$, $d^{p^\prime}$ is closed and
 $\dom(d^p)\subset  \dom (d^{p^\prime})$,  we have
$$\alpha=\lim_{j\to \infty} d^p f_j =\lim_j d^{p^\prime} f_j=0.$$
Thus $\theta=0$ 

As observed in the proof of closability, if $f\in \dom(d_\HH)$ then
$\I^*(f)\in \dom(d^q)$ all $1\le q<\infty$ and
$$\I^*(d_\HH f)=d^{q}(\I^*f), \qquad f\in \dom(d_\HH)$$
From this and using (i), (ii)
$$(\I^*d_\HH)^{p,p^\prime}
=\left(d^{p^\prime}\I^*\left\vert_{\dom(d_\HH)}\right.\right)^{p,p^\prime}
\subset d^{p^\prime}\I^*$$
where $({})^{p,p^\prime}$ indicates the closure as an operator from
$L^p(\C_{x_0}M;\R)$ to $L^{p^\prime}(\C_0\R^m;\R)$.

The inclusion (\ref{domain-d}) will follow if we show
$$(\I^*d_\HH)^{p,p^\prime}=(\I^*d^p)^{p,p^\prime}.$$
This is clear since if $f\in \dom(\I^*d^p)\subset\dom(d^p)$ there exist
$f_j\in \dom(d_\HH)$ with $f_j\to f$ in $L^p$ and $d_\HH f_j\to d^pf$
in $L^p$, but by continuity of $\I^*$ this implies that
$f\in \dom((\I^*d_\HH)^{p,p^\prime})$.
\end{proof}

Let $ \D^{p,1}(\C_{x_0}M;\R)$, $\D^{p,1}(\C_{x_0}M)$ or $\D^{p,1}$
denote the domain of $d^p$ with its graph norm
$$\|f\|_{\D^{p,1}}=\left(\n{d^pf}_{L^p}^p +\n{f}^p_{L^p}\right)^{1\over p}.$$
Note that these spaces depend on the choice of $\Cyl\subset \dom(d_\HH)\subset BC^2$.
But see \S\ref{se:weak-uniqueness}.

The boundedness of $\I^*$ on $\D^{2,1}$ in the next corollary
 was known from  Aida-Elworthy \cite{Aida-Elworthy} for gradient
Brownian stochastic differential equations.

\begin{corollary}
\label{co:closed-range}
The pull back $\I^*$ determines a continuous linear map
$$\I^*: \D^{p,1}(\C_{x_0}M;\R)\longrightarrow
 \D^{p^\prime,1}(\C_0\R^m;\R),$$
for $1<p^\prime<p<\infty$, with the property that 
\begin{equation}
\label{Ito-map-norm}
\|f\|_{\D^{p,1}}\le \|\I^*(f)\|_{\D^{p,1}}, \quad  \hbox{for } f\in \cup_{q>1} \D^{q,1}.
\end{equation}
 If condition $(M_0)$ holds we can take  
$1\le p^\prime=p<\infty$ and then the map $\I^*$ has closed range in the case of $p=p'$,
and (\ref{Ito-map-norm}) holds for $f\in \D^{1,1}$.
\end{corollary}

\begin{proof}  From the theorem
$\dom(d^p)\subset \dom(d^{p^\prime}\I^*)$. A comparison result of
H\"ormander,  (see Yosida \cite{Yosida-book}, Theorem 2, \S 6 of chapter II, p79)
therefore implies that there is a constant $C_{p,{p^\prime}}$
with
$$\n{d^{p^\prime}\I^*f}_{L^{p^\prime}}\le C_{p,{p^\prime}} \|f\|_{p,1}.$$
(Alternatively use Theorem \ref{th:I-star}.)
This, plus the continuity of $\I^*$ on $L^p$, gives the required
continuity of $\I^*$ on $ \D^{p,1}(\C_{x_0}M;\R)$ into
$ \D^{p^\prime,1}(\C_0\R^m;\R)$.

Inequality (\ref{Ito-map-norm}) holds by (\ref{I-star-domination}) and the intertwining (\ref{domain-d}),
and implies that $\I^*$ has closed range when $p=p'$ and Condition $(M_0)$ holds.
\end{proof}

\begin{remark}
\begin{enumerate}
\item[(1)]
As we see in Theorem \ref{th:Markov-unique-2} below, an outstanding question is 
whether $\I^*[\D^{p,1}(\C_{x_0}M;\R)]=\D^{p,1}_{\F^{x_0}}(\C_0\R^m;\R)$, the space of $\F^{x_0}$ 
measurable elements of $\D^{p,1}(\C_0\R^m;\R)$. This appear to be unknown even for gradient Brownian 
motion systems, with Levi-Civita connection, on any space with curvature (e.g. $S^n$ any $n$).
\item[(2)]
Even if Condition $(M_0)$ does not hold we see that for $1<q\le p<\infty$,
$$\I^*[\D^{q,1}\left(\C_{x_0}M;\R\right)]\cap \D^{p,1}\left(\C_0\R^m;\R\right)$$
is closed in $\D^{p,1}(\C_0\R^m;\R)$. However, even given Condition $(M_0)$, we do not know whether $\I^*(f)\in\D^{p,1}(\C_0\R^m;\R)$ and $f\in \D^{q,1}(\C_{x_0}M;\R)$, some $1<q<p<\infty$, imply that $f\in \D^{p,1}(\C_{x_0}M;\R)$. See section \ref{se: weak-differentiability} below.
\end{enumerate}
\end{remark}

\section{The divergence operator and the spaces $\D^{p,1}\HH$}
\label{section-divergence}
\subsection{The divergence operator  $\div$ }
 
From now on we shall take $\Dom(d_\HH)$ to be closed under multiplication by elements of $\Cyl$, the set of smooth cylindrical functions on the path space.  Assume $1<p<\infty$ and ${1\over p}+{1\over q}=1$. Define
$$\div\equiv {\div}^{p}: \Dom({\div}^p)\subset L^p\Gamma\HH\longrightarrow L^p(\C_{x_0}M;\R)$$
to be the cojoint of $-d^q$.  That is $V\in \Dom ({\div}^p)$ if and only $V^{\#}$ is in $\Dom(d^*)$
and ${\div}^p V :=-d^*(V^{\#})$.  So $V\in \Dom({\div}^p)$ if  and only if $V$ is in $L^p$ and
 there is a constant $C(V)$ such that,
 $$\left |\int df(V)d\mu_{x_0}\right |_{L^p}\le  C(V)\cdot |f|_{L^p}, \quad \forall f\in \D^{q,1}(\C_{x_0}M; \R).$$
 If Condition $(M_0)$ holds we can take $p=1$ and $q=\infty$.
If it is necessary to distinguish the underlying path spaces for the divergence operator, we  shall use $\Dom_{\C_{x_0}M}({\div}^p)$ or  $\Dom_\Omega({\div}^p)$.
 
 Define  $\,\Y:  \HH\to H$ whose restriction to $\HH_\sigma$, 
 $\sigma\in \C_{x_0}M$, is given by
\begin{equation}
  \label{Y}
  \Y_\sigma(h)(\cdot)=\int_0^\cdot Y_{\sigma_s}\left( {\D\over ds} h_s \right)\,ds, \qquad h\in \HH_\sigma^1.
  \end{equation}
Then $ \overline{T\I}_\sigma(\Y_\sigma): H\to H$ is the identity map.

For any $h: \C_{x_0}M\to H$ define $\K^\perp h: \C_{x_0}M\to H$ by
\begin{equation}\label{divergence-2}
(\K^\perp h)(\sigma)_t=\int_0^t K^\perp (\sigma_s)\dot h(\sigma)_s ds, \qquad 0\le t\le T
\end{equation}\label{divergence-3}
and for $h: \C_0\R^m\to H$ write
$$(\K^\perp h)(\omega)_t=\int_0^t \K^\perp (x_s(\omega))\dot h(\omega)_s ds, \qquad 0\le t\le T.$$
Note that  for all $h: \C_0\R^m\to H$,
  \begin{equation}\label{divergence-3-3}
 \overline{T\I}(\K^\perp h)=  \overline{T\I}( h), \quad  \overline{T\I(\K^\perp h)}=  \overline{T\I( h)}.
 \end{equation}
 
 \begin{proposition}
 \label{pr:divergence-2}
Let $1< p<\infty$. For any $h\in \Dom (\div^p )$ on the Wiener space,  and any $0<\delta<p-1$, $\overline{T\I(h)}\in \Dom (\div^{p-\delta} )$ on the path manifold. Furthermore
\begin{equation}\label{divergence-4}
\div (\overline{T\I(h)})=\overline{\div h} =\overline{-\int_0^T \langle  \dot h_s, dB_s\rangle}.
\end{equation}
If condition $(M_0)$ holds we take $\delta=0$ and also allow $p=\infty$.
\end{proposition}

\begin{proof} 
Take $h\in \Dom ({\div}^p) \in L^p(\C_0\R^m; H)$ then by Corollary \ref{co:closed-range}, for $\epsilon>0$,
$$\int d(\I^*f)(h) dP=-\int \I^*f(\div h) dP, \quad \forall \, f\in \D^{q+\epsilon,1}(\C_{x_0}M).$$
This implies that for all $f\in \D^{q+\epsilon,1}(\C_{x_0}M)$,
$$\int df\left (\overline{T\I(h)}\right) d\mu_{x_0}=\int \I^*(df)(h)dP=-\int \I^*f(\div h) dP
=-\int f(\overline{\div h}) d\mu_{x_0}$$
and so for all $ h\in \Dom ({\div}^p)$,
$$\div (\overline{T\I(h)})(\sigma)=\overline{\div h}(\sigma)$$
and (\ref{divergence-4}) holds since $\div h=-\int_0^T \langle  \dot h_s, dB_s\rangle$,  the Skorohod integral.
\end{proof}

\begin{corollary}
\label{divergence-skorohod}
Let $1<p'<p<\infty$. An $\HH$-vector field $V$ is in the domain of $\div^{p'}$ if it is in $L^{p'}$ and the vector field $h$ on $\C_0\R^m$ given by
 $$h:=(\Y V)\circ \I: \C_0\R^m\to H$$  is such that $\dot h$ is Skorohod integrable, i.e. $h\in \Dom(\div^p)$ for $\C_0\R^m$. If so 
 \begin{equation}\label{divergence-5}
 (\div V)(x_\cdot) =-
\E\left\{ \int_0^T  \left.\left \langle {\D \over ds} V(x_\cdot), X(x_s) dB_s \right \rangle_{x_s} \right| \F^{x_0}\right\}
 \end{equation}
 where the right hand side is interpreted as the Skorohod integral.  In this case
$$\div(V)=\overline{\div(\I^*(\Y(V)))}.$$
If condition $(M_0)$ holds we can take $1<p=p' \le \infty$.
\end{corollary}
 \medskip
 
 Let   $\D_{\F^{x_0}}^{p,1}$ be the closed subspace of  $\D^{p,1}(\Omega;\R)$ consisting of $\F^{x_0}$ measurable functions. 
 
 \begin{corollary}
\label{le_derivative}
 If $g\in \D^{p,1}_{\F^{x_0}}$ then $dg(-)=dg(\K^\perp -)$ almost surely.
\end{corollary}
 
\begin{proof} 
 Take $h\in H$. Then $\K^\perp h\in \D^{q,1}$, a subset of $\Dom(\div)$ by Kree-Kree \cite{Kree-Kree}.
 Since  by Proposition \ref{pr:divergence-2} 
 $\overline{\div h}=\div\left(\overline {T\I(h)}\right)=\div\left(\overline {T\I(\K^\perp h)}\right)$,
 \begin{eqnarray*}
\E dg(h) &=&-\E g \div h=-\E\left( g\E\{\div h |\F^{x_0}\}\right)\\
 &=&-\E \left(g \div\left (\overline {T\I(h)}\right)\circ \I\right)=-\E\left( g \div\left (\overline {T\I(\K^\perp h)}\right)\circ \I \right)\\
 &=&-\E \left(g \div(\K^\perp h)\right) =\E dg(\K^\perp h).
 \end{eqnarray*}
 Replace $h$ by $\lambda h$ where $\lambda \in \Cyl$ to conclude that almost surely
 $dg(h)=dg(\K^\perp h)$.
 \end{proof}
 
Finally we observe the following version of Corollary \ref{divergence-skorohod}:
\begin{proposition}
For an $L^q$ $\HH$-1-form $\phi$, if $\E \{\I^*\phi |\filt \}\in \Dom_\Omega(d^*)$ then $\phi$ belongs to 
$ \Dom_{\C_{x_0}M}(d^*)$ and 
$$(d^*\phi)_\sigma=\E\{d^*\E\{\I^*\phi|\filt\}|\I=\sigma\}=\overline{d^*[\phi_\I(\overline{T\I}_\I(-))]}(\sigma).$$
\end{proposition}
\begin{proof} 
Just note that  for $f\in \Cyl$, 
\begin{eqnarray*}\int\langle d(\I^*f), \E\{\I^*\phi|\filt\}\rangle dP
&=&\int \langle \I^*f,  \E\{d^*\E\{\I^*\phi|\filt\}|\filt \}\rangle dP\\
&=&\int  f(\sigma)  \E\{d^*\E\{\I^*\phi|\filt\}|\I=\sigma\} d\mu_{x_0}(\sigma)
\end{eqnarray*}
on one hand and, since $\E\{\I^*(-)|\F^{x_0}\}$ is an isometry on one forms,
\begin{eqnarray*}
\int\langle d(\I^*f), \E\{\I^*\phi|\filt\}\rangle d\mu &=&\int\langle \E\left\{\I^*(df)|\filt\right\}, \E\{\I^*\phi|\filt\})\rangle dP \\
&=&\int \langle df, \phi \rangle d\mu_{x_0}
\end{eqnarray*}
on the other hand.
\end{proof}

\subsection{Hilbert space valued $L^p$ functions}
\label{s_Hilbert}

Let $G$ be a  separable Hilbert space and $B=\HH$ or $\HH^*$, or a similar `tensor bundle'.
 Denote by $L^p\Gamma B$ and $L^p\Gamma(G\otimes B)$ respectively the $L^p$ sections of $B$ and
 those of the tensor product of the `bundle' $B$ with the trivial $G$ bundle over $\C_{x_0}M$. 
 We always use $\otimes$ to refer Hilbert space completions and $\otimes_0$ the incomplete algebraic tensor products.
 For each densely defined linear map $T$ from  $L^p(\C_{x_0}M; \R)$ to  $L^p\Gamma B$ there is a naturally defined  linear operator $T^G\equiv \id\otimes T $ from $L^p(\C_{x_0}M; G)$ to $ L^p\Gamma (G\otimes B)$ with domain $\Dom(T^G)\equiv \Dom(T)\otimes_0 G$, namely
$$\Dom(T^G)=\big\{F: \C_{x_0}M\to G \;\big|\;
     F(\sigma)=\sum_{j=1}^n f_j(\sigma) g_j,   g_j\in G, f_j\in \Dom(T), n\in{\mathbb N} \big\}$$
and such that $T^G(f\otimes g)=T(f)\otimes  g$ for $f\in Dom(T)$ and $g\in G$.

\begin{proposition}
\label{closable_G}
If $T$ with $\Dom(T)$ is a closable operator then so is $T^G\equiv \id\otimes T$ with 
$\Dom (T^G)=G\otimes_0 \dom(T)$.\end{proposition}
\begin{proof} Take $F_n \in \Dom(T^G)$ converging to  $0$ in $L^p(\C_{x_0}M; G)$ with $TF_n\to \alpha$.
 For an orthonomal basis  $\{e_j\}$ for $G$ 
write $F_n(\sigma)=\sum_{j=1}^\infty f_n^j(\sigma)e_j$ and  $\alpha=\sum_j e_j \otimes  \alpha_j$.  Then $Tf_n^j\to \alpha_j$ and $f_n^j\to 0$  in $L^p$ as $n\to \infty$.  Consequently  $\alpha_j=0$.
 \end{proof}
  
Now take $T$ to be $d_\HH$ with $\Dom(d_\HH)$ the set of smooth cylindrical functions.
 Define 
$$d\equiv  d^p\equiv d^{p,G}: \Dom(d^{p,G})\subset L^p(\C_{x_0}M;G)\to
L^p\Gamma(G\otimes \HH^*)\sim L^p\Gamma(\LL_2(\HH;G))$$
 to be the closure of ${(d_\HH)}^G$ which exists by Proposition \ref{closable_G}. Its domain shall be denoted by $\D ^{p,1}(\C_{x_0}M;G)$ and is the closure of $\Dom(d)\equiv  G\otimes_0 \Dom(d_\HH)$ under the graph norm.

The following elementary lemma is useful in section \ref{se: D-HH-spaces}.
   \begin{lemma}
   \label{le:multiplication}
 Suppose $1<p<\infty$, allowing $p=1$ if Condition $(M_0)$ holds. Let $G_1$ and $G_2$ be two Hilbert spaces and 
 $\theta: M\to \LL(G_1; G_2)$ a $C^1$ map. Suppose $f: C_{x_0}M\to L^2([0,T];G_1)$ is in
  $\D^{p,1}\left(\C_{x_0}M;L^2([0,T];G_1)\right)$. Then
 $$ \Theta(f): C_{x_0}M\to L^2([0,T]; G_2)$$
 given by $\sigma\mapsto [t\mapsto \theta(\sigma_t)(f(\sigma)_t)$ is in $\D^{p,1}\left(\C_{x_0}M;L^2([0,T];G_2)\right)$ with
 \begin{equation} 
 \label{multiplication-rule}
[d^p\left(\Theta(f)\right)_\sigma(v)]_t=(d\theta)_{\sigma_t}(v_t)(f(\sigma)_t)+\theta(\sigma_t)(d^pf_\sigma(v))_t.
 \end{equation}
  \end{lemma}
  
  \begin{proof}
  Take $f_n\to f$ in $\D^{p,1}\left(\C_{x_0}M;L^2([0,T];G_1)\right)$ where 
  $$f_n(\sigma)=\sum_{j=1}^{k_n} \lambda_j^n(\sigma) h^{j}$$
  with $\lambda_j^n$ real valued smooth cylindrical functions and  $h^j\in L^2([0,T];G_1)$. Then 
  $$\Theta(f_n)(\sigma)(t)= \sum_{j=1}^{k_n} \lambda_j^n(\sigma) \theta(\sigma_t) h^{j}_t.$$
  Clearly $\Theta(f_n)\to \Theta(f)$ in $L^p\left(\C_{x_0}M; L^2([0,T];G_2)\right)$ and 
  $\Theta(f_n)$ is Fr\'echet differentiable as a map into $L^2([0,T];G_2)$. For $v\in \HH_\sigma$,
  \begin{eqnarray*}
  [d_\HH\left(\Theta(f_n)\right)_\sigma(v)]_t
  &=&\sum_{j=1}^{k_n} \left(d_\HH \lambda_j^n\right)_\sigma (v) \theta(\sigma_t)h^{j}_t +\sum_{j=1}^{k_n} \lambda_j^n(\sigma) (d\theta)_{\sigma_t}(v_t)h^{j}_t \\
&=& (d\theta)_{\sigma_t}(v_t)(f_n(\sigma)_t)+\theta(\sigma_t)(df_n)_\sigma(v)_t) 
  \end{eqnarray*}
 and we see that $d_\HH\left(\Theta(f_n)\right)$ converges in
  $L^p\left(\C_{x_0}M; L^2([0,T];G_2)\right)$ with limit given by the right hand side of (\ref{multiplication-rule}).
   Finally just  approximate $\theta(\sigma_t)$ by a sequence of functions $\theta_n\in \D^{p,1}\left( L^2([0,T];G_2)\right)$, for example
   set $$\theta_n(\sigma)_t={t-t_j^n\over t_{j+1}^n-t^n_j} \theta\left(\sigma(t^n_{j+1})\right)+\theta\left(\sigma(t_j^n)\right),
   t_j^n\le t<t^n_{j+1}$$
   for suitable partitions $0\le t_1^n\le \dots \le t_j^n\le \dots \le T$.   \end{proof}

 \subsection{Pull back of Hilbert space  valued functions and $H$-forms}
 \label{Pull back of Hilbert space valued forms}
 
 We follow the notation of section \ref{s_Hilbert}.
 For $p=2$ there is the canonical isometry of $L^2(\C_{x_0}M;\R)\otimes G$ with $L^2(\C_{x_0}M;G)$ mapping $\Dom(d) \otimes G$ onto $\Dom( d^G)$.

By a $G$-valued $\mathcal{H}$-1-form $\phi $ on $\C_{x_0}M$ we mean a
measurable section of the bundle $\LL_2\left( \mathcal{H}; G\right) $,
(or equivalence class of such sections under almost sure equality). There is the
standard identification of $\LL_2\left( \mathcal{H}_\sigma ; G\right) $
with $G\otimes \mathcal{H}_\sigma $. It is given by 
 $(g\otimes h)(v)=\left\langle h,v\right\rangle _{\mathcal{H}_\sigma }g$. We shall use $\phi ^{\#}$ to denote the
section of $G\otimes \mathcal{H}$ corresponding to a $G$-valued $\mathcal{H}$-1-form $\phi $.
 Note that we can differentiate such $\phi
^{\#}$ to obtain $\left( 1_G\otimes
{{\D }\over ds}\right) \phi _{\sigma ,s}^{\#}\in
G\otimes T_{\sigma (s)}M$, $0\le s\le T$, $\sigma \in \C_{x_0}M$.
 
   If $f: \C_{x_0}M\to G$ is in $\D^{p,1}$,  its differential $(df)_\sigma\in \LL_2(\HH; G)$ then determines the
 gradient $\nabla f\in G\otimes \HH$ by $\nabla f(\sigma)=(df)_\sigma^{\#}$. In the Nualart-Pardoux
 notation we obtain $(s\mapsto D_sf)\in \Gamma(G\otimes L^2\Epsilon)$ by the isometries 
 $$\LL_2(\HH;G)\sim G\otimes \HH\stackrel{1\otimes {\D\over ds} }{\sim}G\otimes L^2\Epsilon.$$
 
\begin{theorem}
 \label{th:pull-back-forms}
Theorem \ref{th:I-star} and Theorem \ref{theorem-closability-2} and Corollary \ref{co:closed-range}
 hold for $G$-valued functions. For a measurable $G$-valued $\mathcal{H}$-1-form $\phi
$ in $L^p$ the pull back $\I ^{*}(\phi ): \C_0\R ^m\rightarrow \LL_2\left( H;G\right) $
 is given by, for $h\in H$,
\begin{equation}
\I ^{*}(\phi )(h)=\int_0^T\left[ \left( 1_G\otimes { \D  \over ds}\right) \phi _s^{\#}\right] ^{\#}\left( \nabla_{T\I _s(h)}X(\tilde{\parals_s}d\beta _s)+X(x_s)(\dot{h}_s)ds\right).
 \label{pull-back-G-form}
\end{equation}
\end{theorem}

\begin{proof} 
 The proofs are essentially the same as those for the real valued case. Here we assume that  $p=2$ 
 to demonstrate the proof. 
  If the $G$-valued $\phi $ has $\phi ^{\#}=\sum_{j=1}^kg_j\otimes \psi _j^{\#}$ where $\psi _j^{\#}\in \Gamma \mathcal{H}$ 
  and $g_j\in G$ for $j=1$ to $k$, we immediately obtain (\ref{pull-back-G-form}) from (\ref{I-phi}). The
general case and the rest of Theorem \ref{th:I-star}  follow by taking limits, in
particular by observing that (\ref{continuity-2}) remains true in the $G$-valued
case. The proof of Theorem \ref{theorem-closability-2} is based on (\ref{continuity-2}) and so is easily
seen to extend to the $G$-valued case. The crucial remark is that the stochastic integral in 
 (\ref{pull-back-G-form}) can be considered as an $\LL_2(H; G)$-valued integral applied to $h$. To see this and perform the necessary estimates we need to show
 $$\int_0^T \left\|\left[ \left( 1_G\otimes { \D  \over ds}\right) \phi _s^{\#}\right] ^{\#}
 \nabla _{T\I_s(-)} X^j\right\|^2 _{\LL_2(H,G)} ds$$
is finite almost surely for $j=1, .\dots, m$.  First noted that the assumption that  $\phi_\sigma: \HH\to G$ is in $L_2$ implies
that so is $V\mapsto \phi_\sigma(\W_\cdot V)$ from $L^2\Epsilon$ to $G$, where $L^2\Epsilon$ is the Hilbert bundle of $E$ valued $L^2$ tangent vector fields on $\C_{x_0}M$.  Moreover 
  $$\phi_\sigma(\W_\cdot V)=\int_0^T\left \langle  \left( 1_G\otimes { \D  \over ds}\right)  \phi_s^{\#}, 
 V_s \right\rangle_{\sigma(s)} ds.$$ 
 The norm of $\phi_\sigma(\W_\cdot -)$ equals $\|\phi_\sigma\|^2_{\LL_2(\HH_\sigma; G)}
$ and  is given by
      $$	\|\phi_\sigma\circ \W\|^2_{\LL_2(\LL^2\Epsilon; G)}
       =  \int_0^T\left \| \left( 1_G\otimes { \D  \over ds}\right)  \phi_s^{\#}\right\|^2_{\LL_2(E_{\sigma(s)}; G)} ds, $$
which is finite. Now for $\sigma=\I(\omega)$,
  \begin{eqnarray*}
&&\int_0^T \left\|\left[ \left( 1_G\otimes { \D  \over ds}\right) \phi _s^{\#}\right] ^{\#}
 \nabla _{T_\omega\I_s(-)} X^j\right\|^2 _{\LL_2(H; G)} ds\\
 &\le&\int_0^T \left\|\left[ \left( 1_G\otimes { \D  \over ds}\right) \phi _s^{\#}\right] ^{\#}\right\|_{\LL_2(E_{\sigma(s)}; G)}^2
\left\| \nabla _{T_\omega\I_s(-)} X^j\right\|^2 _{\LL(H; E_{\sigma(s)})}  ds\\
&\le&\|\phi_\sigma\circ \W\|^2_{\LL_2(L^2\Epsilon; G)} \cdot \sup_s \left\| \nabla _{T_\omega\I_s(-)} X^j\right\|^2 _{\LL(H;E_{\sigma(s)})} \\
&\le&\hbox{constant } \cdot   \|\phi_\sigma\circ \W\|^2_{\LL_2(L^2\Epsilon; G)} \sup_s\|T_\omega \I\|^2_{\LL(H; T_\sigma \C_{x_0}M)}<\infty.
 \end{eqnarray*}       
  There are no difficulties with Corollary \ref{co:closed-range}.

 \end{proof}

 \subsection{The space $\D^{p,1}\HH$}
\label{se: D-HH-spaces}

 For $1< p<\infty$, and allowing $p=1$ if Condition $(M_0)$ holds, define
\begin{equation}
\D^{p,1}\HH=\left\{V \in L^p \Gamma\HH \quad \big | \quad \Y_\cdot( V(\cdot)) \in \D^{p,1}(\C_{x_0}M; H)\right\},
\end{equation}
 Equip this space with the obvious norm:
 $$\|V\|_{\D^{p,1}\HH}
 =\| \Y( V(\cdot))\|_{\D ^{p,1}(\C_{x_0}M;H)}=\left(\E|\Y(V)|^p  +\E |d(\Y(V))|^p\right)^{1\over p}.$$
 This  depends only on the connection on $E$ not on the specific stochastic differential equation (\ref{sde}), or equivalently 
 not on the particular choice of $X$ as can be easily seen by Lemma \ref{le:multiplication}.
 In fact $\D^{p,1}\HH$ can also be described by covariant differentiation, see \S\ref{se:High-Sobolev} below. 
  Similarly we can define $\D^{p,1}\HH^*$:
\begin{equation}
\D^{p,1}\HH^*=\left\{\phi \in L^p \Gamma\HH^* \quad \big | \quad \phi(\X(-)) \in \D^{p,1}(\C_{x_0}M; H^*)\right\}.
\end{equation}

We have the following  analogue of a fundamental result of Kree-Kree \cite{Kree-Kree} for $\C_0 \R^m$:

\begin{theorem}\label{th:divergence-21}
For $1< p'<p<\infty$, the set $\D^{p,1}\HH$ is contained in $\Dom(\div^{p'})$ and  $\div^{p'}: \D^{p,1}\HH\to L^{p'}(\C_{x_0}M; \R)$  is continuous. If Condition $(M_0)$ holds we may take $p=p'$.
\end{theorem}

\begin{proof}  Take $r$ with $p'<r<p$.
	 If $V\in \D^{p,1}\HH$ then $\Y(V) \circ \I\in \D^{r,1}(\Omega; H)\subset \Dom_\Omega(\div^r)$, by Theorem \ref{th:pull-back-forms} and the corresponding
result for Wiener space of Kree-Kree \cite{Kree-Kree}. Then $V=\overline{T\I(\Y(V))} \in \Dom_{\C_{x_0}M}(\div^{p'})$ by Proposition \ref{pr:divergence-2}.
Finally note that
$$\D^{p,1}\HH \stackrel{\I^*(\Y-)}{\longrightarrow}\D^{r,1}(\Omega; H)
\stackrel{\div-}{\longrightarrow}  L^{p'}(\Omega;\R) 
\stackrel{\hbox{conditional expectation}} {\longrightarrow} L^{p'}(\C_{x_0}M;\R)$$
is a continuous map and by Corollary \ref{divergence-skorohod}, agrees with  the divergence operator restricted to $\D^{p,1}\HH$.
 Note that  the continuity also follows from the closed graph theorem.
\end{proof}

The following is a compliment of Proposition \ref{pr:divergence-2}:
\begin{proposition}
\label{pr-Domain-divergence-1}
For $1<p<\infty$, set 
\begin{equation}
\label{total_set_of_H}
U=\left\{h\in \D^{p,1}(\Omega, H)  \quad \left |  \quad \overline{T\I(h)}\in \D^{p,1}\HH  \right.\right\}.
\end{equation} 
Then $U$ is total in $\D^{p,1}(\Omega; H)$ and thus  total in the domain, 
$\Dom(\div^p)$, of the divergence  on Wiener space.  
\end{proposition}

\begin{proof} Consider the family of  functions
$$U_1\equiv \left\{k' _\cdot\exp{\left(\int_0^T \langle \dot  k_s, dB_s \rangle-{1\over 2}\int_0^T|\dot k_s|^2 ds\right)}\;\Big |\; k', k\in H \right\}.$$

Since the exponential martingales are total in $\D^{p,1}(\Omega;\R)$ it is clear from the definition that $U_1$ is total in $\D^{p,1}(\Omega; H)$ and so it is sufficient to show that 
$$\Image [\overline{T\I(h)}: h\in U_1] \subset \D^{p,1}(\HH),$$
or equivalently that
$\Y(\overline{T\I(h)})$ belongs to $\D^{p,1}(\C_{x_0}M; H)$. In fact for 
$$h_\cdot=k'_\cdot\exp{\left(\int_0^T  \langle \dot k_s, dB_s \rangle-{1\over 2}\int_0^T |\dot k_s|^2 ds\right)},$$
we can write
\begin{eqnarray*}
h_\cdot&=&k'_\cdot \exp\left({\int_0^T\langle \dot k_s, {\parals_s} d\breve B_s \rangle-{1\over 2}\int_0^T|K^\perp(x_s) \dot k_s|^2 ds}\right)\\
&&\qquad \cdot \exp{\left(\int_0^T\int_0^T\langle \dot  k_s, K(x_s)dB_s \rangle-{1\over 2}\int_0^T |K(x_s)\dot  k_s|^2 ds\right)}.
\end{eqnarray*}
\begin{eqnarray*}
\overline{T\I_t(h_\cdot)}&=&\exp{\left(\int_0^T \langle \dot k_s, {\parals_s} d\breve B_s \rangle-{1\over 2}\int_0^T | K^{\perp}(x_s)\dot k_s|^2 ds\right)} \\
&& \times \;\E\left\{ \left. T\I_t \left(k'_\cdot
\exp{\int_0^T \langle \dot k_s, K(x_s)dB_s \rangle-{1\over 2}\int_0^T | K(x_s) \dot k_s|^2 ds}\right)  \right| \filt\right\}.
\end{eqnarray*}
Set $$f_t=\E\left\{\exp{\left(\int_0^T \langle \dot  k_s, K(x_s)dB_s \rangle-{1\over 2}\int_0^T |K(x_s) \dot k_s|^2 ds\right)} | \F^{x_0}\vee \F_t\right\}.$$
Then $(f_t, 0\le t\le T)$ is a martingale with respect to $\{\F^{x_0}\vee \F_t\}$ and so
$$f_t =\exp{\left(\int_0^t \langle \dot  k_s, K(x_s)dB_s \rangle-{1\over 2}\int_0^t | K(x_s) \dot k_s|^2 ds\right)},$$
giving
 \begin{eqnarray*}
&&\E\left\{ \left. T\I_t(k'_\cdot)
\exp{\left(\int_0^T \langle \dot k_s, K(x_s)dB_s \rangle-{1\over 2}\int_0^T |  K(x_s) \dot k_s|^2 ds\right)}\right| \F^{x_0} \right\}\\
&=&\E\left\{ \left.T\I_t(k'_\cdot)
\exp{\left(\int_0^t\langle \dot k_s, K(x_s)dB_s \rangle-{1\over 2}\int_0^t|  K(x_s) \dot k_s|^2 ds\right)} 
\right| \F^{x_0} \right\}.
\end{eqnarray*}
On the other hand,  if  we set
$$V_t=T\I_t(k'_\cdot)
\exp{\left(\int_0^t\langle \dot k_s, K(x_s)dB_s \rangle-{1\over 2}\int_0^t|K(x_s)\dot  k_s|^2 ds\right)},$$
then as in Elworthy-LeJan-Li \cite{Elworthy-LeJan-Li-book}, Elworthy-Li 
\cite{Elworthy-Li-Hodge-1}
\begin{eqnarray*}
\D V_t&=&\nabla X(V_t)dB_t+{1\over 2}\nabla X(V_t)(K(x_t)\dot k_t)dt\\
&&+T\I_t(k')
\exp{\left(\int_0^t\langle\dot   k_s, K(x_s)dB_s \rangle-{1\over 2}\int_0^t|  K(x_s) \dot  k_s|^2 ds\right)}
 \langle \dot   k_t, K(x_t)dB_t \rangle\\
 &=&\nabla X(V_t)dB_t+{1\over 2}\nabla X(V_t)(K(x_t)\dot k_t)dt+\langle \dot  k_t, K(x_t)dB_t \rangle V_t.\\
\end{eqnarray*}
Consequently, for $\bar V_t(\sigma)=\E\{V_t |x_\cdot =\sigma\}$, 
\begin{eqnarray*}
{\D\overline{V_t} \over dt}\circ \I
&=& {1\over 2}\nabla_{\overline{V_t}(x_\cdot)} X(K(x_t)\dot k_t)dt.
\end{eqnarray*}
Thus
\begin{eqnarray*}
{d\over dt} \Y(\overline{V_t}(x_\cdot))=Y\left({\D\over dt} \bar V_t(x_\cdot)\right)
&=&{1\over 2}Y\left(\nabla X(\overline{T\I_t}(\Y(\overline{V_\cdot(x_\cdot)})(K(x_t)\dot k_t)\right).
\end{eqnarray*}
Since solutions of such a stochastic differential equation are in $\D^{p,1}$ for all $1<p<\infty$ by standard results,  $\Y(\overline V_\cdot)\in \D^{p,1}(\C_{x_0}M; H)$. See Lemma \ref{le:cylindricals-5} below.
Finally note that
$$\overline{T\I_t(h_\cdot)}_\sigma=\exp{\left(\int_0^T \langle \dot k_s, {\parals_s} d\breve B_s \rangle-{1\over 2}\int_0^T | K^\perp(\sigma_s) \dot k_s|^2 ds\right)} \overline{V_t}(\sigma).$$
Since such exponential martingales belong to $\D^{q,1}$ for all finite $q$ we see that 
$$\sigma\mapsto\Y(\overline{T\I_t(h_\cdot)}_\sigma )=\Y(\bar V_t(\sigma))  \exp{\left(\int_0^T \langle \dot k_s, {\parals_s} d\breve B_s \rangle-{1\over 2}\int_0^T |K^\perp(\sigma_s) \dot k_s|^2 ds\right)} $$
belongs to $\D^{p,1}(\C_{x_0}M;H)$ for all $1\le p<\infty$.
\end{proof}

 \section{On the Markov Uniqueness of $d$} 
 \label{section:Markov-Unique}
 Throughout section \ref{section:Markov-Unique} we take $\Dom(d_\HH)=\Cyl$. To define the weak derivatives we shall need to assume $\D^{q,1}\HH^*\subset \Dom(d^p)^*)$, which is guaranteed by Theorem \ref{th:divergence-21} if Condition $(M_0)$ holds or if $X$ is injective.
\subsection{Weak Differentiability }
\label{se: weak-differentiability}
For ${1\over p}+{1\over q}=1$ and $1<p<\infty$,  the weak Sobolev space $ W^{p,1}(\C_{x_0}M)$, abbreviated as $W^{p,1}$,  is the domain of the adjoint of the restriction of $(d^p)^*$ to $\D^{q,1}\HH^*$: 
  \begin{equation}W^{p, 1} = \Dom \left(\big( (d^p)^*\left|_{\D^{q, 1} \HH^*} \right.\big)^*\right)
  \end{equation}
   furnished with its graph norm.
More precisely a function $f$ belongs to $W^{p,1}$ if and only if it is in $L^p$ and there is a constant $C(f)$ such that 
$$\left |\int_{\C_{x_0}M} f d^*\phi \;d\mu_{x_0} \right|\le C(f)|\phi|_{L^q}, \; \forall \phi\in \D^{q,1}\HH^*.$$
  Equivalently,
$$W^{p,1}=\left\{f\in L^p: \big|\int \div(V) f d\mu_{x_0} \big|\le C|V|_{L^q},  \hbox { for all }V\in \D^{q,1} \HH, \hbox{ some } C \right\}.$$
 If $f\in W^{p,1}$ it has a ``weak derivative'' $\tilde df\in L^p\Gamma \HH^*$ defined by 
\begin{equation}
\label{weak-derivative}
 \int \tilde d f(V(\sigma) )  d\mu_{x_0}(\sigma) =-\int f (\sigma) \div V(\sigma) d\mu_{x_0}(\sigma) , \; \forall V\in \D^{q,1}\HH,
 \end{equation}
and  $$\tilde d=\left({(d^p)}^*|_{\D^{q,1}\HH^*}\right)^*$$
as a closed densely defined operator on $L^p$. Denote by   $\tilde \nabla$  the corresponding weak gradient with values in $L^p\Gamma\HH$. Note that $\tilde d$ is an extension of $(d, \D^{p,1})$.

Recall that $\Cyl$ denotes the space of smooth cylindrical functions on $\C_{x_0}M$. Set
\begin{equation} 
\label{cylindrical H-forms}
\hbox{Cyl}^0 \HH^{\ast} =\hbox{ linear span }  \{ g\, dk\; \vert \;g, k \in
  \Cyl \} 
  \end{equation}
  and define
\begin{equation} {}^0 W^{p, 1}= \Dom ( d^{\ast}\; |\; {\Cyl}^0 \HH^*)^*.
\end{equation}
From Proposition \ref{dense-cylindrical} below,
 $ {\Cyl}^0 \HH^*\subset \D^{q, 1} \HH^*\subset \Dom(d^*)$ and so
  $$\D^{p, 1} \subseteq\;   W ^{p, 1} \subseteq  {}^0 W^{p, 1}.$$

In Theorem \ref{th-Sobolev-sp}  below we show that $W ^{p, 1} =  {}^0 W^{p, 1}$.
 \begin{theorem}
   \label{th:weak_Sob}
 Suppose Condition ($M_0$) holds. For $1<p<\infty$, the following are equivalent:
 \begin{enumerate}
 \item[(i)]  $f\in W^{p,1}(\C_{x_0}M;\R)$
 \item[(ii)] $ \I^*(f)\in \D^{p,1}(\C_0\R^m; \R)$
 \item [(iii)] $f\in W^{r,1}(\C_{x_0}M;\R)\cap L^p(\C_{x_0}M;\R)$ some $r\in(1,p)$ and the weak derivative 
 $\tilde d f$ is in $L^p$. 
 \end{enumerate}
 Moreover 
  \begin{enumerate}
 \item[(iv)] there is the following intertwining of $\tilde d$ and $\I^*$: 
 $$d(\I^*f)=\I^*(\tilde d f),  \hbox{ for all  }f\in W^{p,1}(\C_{x_0}M;\R).$$
 \item[(v)] if $f \in W^{p, 1} (\C_{x_0}M;\R)$, 
  \begin{equation}
  \label{weak-22}
    (\tilde d f)_{\sigma} =\mathbb{E} \{ d (\mathcal{I}^{\ast}  f 
    )_{\omega} | x_. ( \omega ) = \sigma \} \Y_{\sigma}.
  \end{equation}
   \end{enumerate}
 \end{theorem}

\begin{proof} 
 Let $V \in \mathbb{D}^{q, 1} \HH$. Suppose $\I^*(f) \in \D^{p, 1}$ then 
 by Corollary \ref{divergence-skorohod} and Kree-Kree \cite{Kree-Kree},
    \begin{eqnarray*}
    \int_{\C_{x_0} M} f \div ( V ) d \mu_{x_0} &=& \int_{\C_0\R^m}
    \mathcal{I}^{\ast} ( f )  \div ( V ) \circ \mathcal{I}d\mathbb{P}  \\
    & =& \int_{\C_0
 \R^m} \mathcal{I}^{\ast} ( f ) \div \mathcal{I}^{\ast}(\Y_- ( V ( - ))
  ) d\mathbb{P}\\
  &=& - \int_{\C_0
 \R^m} d (\mathcal{I}^{\ast} ( f ) )_{\omega} (\Y_{x_{\cdot}(\omega)}
  ( V (x_{\cdot}( \omega )) ) d\mathbb{P}( \omega )\\
  &=&- \int_{\C_{x_0} M}
  \mathbb{E} \{ d (\mathcal{I}^{\ast} ( f ) )_{\omega} | x_\cdot ( \omega ) =
  \sigma \} \Y_{\sigma} ( V ( \sigma ) ) d \mu_{x_0}(\sigma).
  \end{eqnarray*}
Thus (ii) implies (i) and (\ref{weak-22}) holds.

 To show (i) implies (ii), suppose $f\in W^{p,1}(\C_0\R^m)$ and take 
$V= \overline{T\I(h)}$ where $h\in U$, as defined in Proposition \ref{pr-Domain-divergence-1}. 
By definition $V \in \D^{q,1}\HH$ and by Proposition \ref{pr:divergence-2} $\div (V)=\overline{\div (h)}$. So
 \begin{eqnarray*}
 \int_\Omega \div(h) \I^*(f)  dP&=&\int_{\C_{x_0}M} \div(V) f d\mu_{x_0}\\
 &\le& C(f) \|\overline{T\I(-)}\| |h|_{L^q},
\end{eqnarray*}
by Corollary \ref{corollary-pull-back}.
Since $U$ is total in $\D^{q,1}(\Omega, H)$  the inequality holds for all $h\in\D^{q,1}(\Omega; H)$.
Consequently $\I^*(f)\in W^{p,1}(\Omega)=\D^{p,1}(\Omega; \R)$, using Sugita \cite{Sugita85}.

Next observe for $f\in W^{p,1}(\C_{x_0}M)$, $h\in \D^{q,1}$,
\begin{eqnarray*}
&&\int_\Omega d(\I^*f)(h)dP = -\int_\Omega \div(h) \I^*(f)  dP=-\int_{\C_{x_0}M} \div(V) f d\mu_{x_0}\\
&=&\int _{\C_{x_0}M} (\tilde df)(V) d\mu_{x_0}=\int_\Omega(\I^*\tilde df)(h) dP,
\end{eqnarray*}
which gives $d(\I^*f)=\I^*(\tilde d f)$.

To see the equivalence of (i) and (iii), take ${1\over r}+{1\over s}=1$. That (i) implies (iii) is trivial.
To obtain (i) from (iii), take $f\in W^{r,1}(\C_{x_0}M;\R)$, with $\tilde d f\in L^p$.  For any $U\in \D^{s,1}\HH$,
$$\left|\int \div(U)(\sigma) f (\sigma)\mu(d\sigma) \right|   
=\left|\int \tilde d f (U)\mu(d\sigma) \right|
\le  \|\tilde d f\|_{L^p}|U|_{L^q}.$$
which, by continuity using Theorem \ref{th:divergence-21}, holds for all $U\in  \D^{q,1}\HH$ if $f\in L^p$, giving (i).
\end{proof}

\begin{corollary}
\label{co-Dirichlet-form}
 Suppose Condition ($M_0$) holds. The symmetric form 
$$\varepsilon(f,g)=\int_{\C_{x_0}M} \langle \tilde df, \tilde dg \rangle d\mu_{x_0}$$
with domain  $W^{2,1}(\C_{x_0}M; \R)$  is a Dirichlet form.\end{corollary}
\begin{proof}
Just observe that,  by Theorem \ref{th:weak_Sob}, the usual chain rule holds for composition on the left by $BC^1$ functions on $\R$.
\end{proof}

Note that if $\I^*[\D^{2,1}(\C_{x_0}M;\R)]=\D^{2,1}_{\F^{x_0}}(\Omega;\R)$ then, by Theorem \ref{th:weak_Sob},  $\D^{2,1}(\C_{x_0}M;\R)=W^{2,1}(\C_{x_0}M;\R)$ and $\I^*d=d\I^*$. In particular we have the Markov uniqueness.  Furthermore there is equality of the  following two Dirichlet forms:
$$\int_{\C_{x_0}M}  |\tilde d f|^2\; d\mu_{x_0}
=\int_\Omega \big| \E\{d \I^*(f) \big | \F^{x_0}\}\big|^2  dP $$
and there is a constant $c$ with
$$\int_{\C_{x_0}M}|\tilde d f|^2 d\mu_{x_0} \le \int_\Omega  |d\I^*f|^2dP
\le c\int_{\C_{x_0}M} |\tilde d f|^2 d\mu_{x_0}, \quad f \in W^{2,1}(\C_{x_0}M;\R),$$
c.f. Driver \cite{Driver-non-equivalence}, Shigekawa \cite{Shigekawa95}.

\begin{corollary}
 Suppose Condition ($M_0$) holds.
If $\I^*[\D^{p,1}]=\D^{p,1}_{\F^{x_0}}$ for the It\^o map $\I$ of one stochastic differential equation which induces $(\mu_{x_0}, \nabla)$, then it holds for all such It\^o maps.
\end{corollary}
\begin{proof}
This is immediate from Theorem \ref{th:weak_Sob} since $W^{p,1}$ depend only on $\mu_{x_0}$ and $ \nabla$.
\end{proof}

\begin{proposition}
For $1<p<\infty$,
 \begin{enumerate}
 \item[(1)]
 $\D^{p,1}(\C_{x_0}M;\R)$ is a closed subspace of $W^{p,1}(\C_{x_0}M;\R)$.
 \item[(2)]
 Set $W_0^{p,1}=\{f\in W^{p,1} | \int f d\mu_{x_0}=0\}$.  Suppose Condition ($M_0$) holds. Then
 $\tilde d: W^{p,1}_0(\C_{x_0}M;\R) \to L^p\Gamma \HH^*$ has closed range and  is a linear isomorphism onto its image.
 \end{enumerate}
\end{proposition}

\begin{proof}  Part (i) is automatic since $d\subset \tilde d$. 
For (2) the continuity holds by definition of the graph norm while injectivity comes from the result for $\C_0\R^m$ using 
Theorem \ref{th:weak_Sob}. 

To show $\tilde d$ has closed range, take $\phi\in L^p\Gamma\HH^*$ with $\phi=\lim_j \tilde d  f_j$ in $L^p$, for $f_j\in W^{p,1}_0(\C_{x_0}M;\R)$.
 Then $\I^*(\tilde d f_j)=d(\I^* f_j)$ by  Theorem \ref{th:weak_Sob}  and
 $$\I^*\phi=\lim_{j\to \infty}\I^*(\tilde d f_j) =\lim_{j\to \infty} d(\I^* f_j)$$
 in $L^p$ by Theorem \ref{th:I-star}. Since $d$ on the Wiener space has closed range, e.g. from Shigekawa \cite{Shigekawa86}, we see $\I^*f_j\to g$ in $\D^{p,1}(\Omega;\R)$ for some $g\in\D^{p,1}_{\F^{x_0}}$. By Theorem \ref{th:weak_Sob},  $g=\I^*(f)$ for some $f\in W^{p,1}_0$ and $f=\lim f_j$ in $L^p$. Since $\tilde d$ is closed,
   $f\in W^{p,1}_0$ with $\tilde d f=\phi$. Thus $\tilde d$ has closed range.
 \end{proof}

\subsection{Markov uniqueness}
\label{section-Markov-unique} 

Let $\underline{t}= \{ t_1, \ldots, t_k \}$ with $0 \le t_1 < \ldots < t_k
  \le T$, and write $M^{\underline{t}} \cong M \times \ldots \times M$. 
  A cylindrical $q$-form on $\C_{x_0} M$ is of the form
  $(\ev_{\underline{t}} )^{\ast} \varphi$ where $\varphi$ is a smooth $q$-form on $M^{\underline{t}}$
  and 
$${\ev}_{\underline{t}} : \C_{x_0} M
  \to M^{\underline{t}}$$
  is the evaluation given by
$$  {\ev}_{\underline{t}}(\sigma)=\left(\sigma(t_1), \dots, \sigma(t_k)\right).$$ 
  \begin{theorem}
  \label{theorem-cylindrical}
  The space $\Cyl^0 \HH^{\ast}$,  defined by (\ref{cylindrical H-forms}), is total
  in the space $\D ^{q  , 1} \HH^*$ for $1 \le q < \infty$.
\end{theorem}
\begin{proof}
 Let   $(\ev_{\underline{t}} )^{\ast} \varphi$ be a typical cylindrical one form on $\C_{x_0} M$.
   In local co-ordinates $\varphi$ can be represented by an
  expression such as $\Sigma_{j = 1}^N \varphi^j dx^j$, some $N$. It
  follows using a partition of unity that $\varphi = \Sigma_{j = 1}^N g^j
  df ^j$, for a finite set of smooth cylindrical functions
  $g^{j}, f^j$ on $M^{\underline{t}}$. Thus $\Cyl^0 \HH^*$ spans  the space of all smooth cylindrical 1-forms. The conclusion can be drawn from Proposition \ref{dense-cylindrical} in \S\ref{section-density} below, which states there exists a set of smooth
   cylindrical 1-forms  which is dense in $\D^{q,1}\HH^*$.  \end{proof}

\begin{remark}
  If our It\^o map is induced by a gradient stochastic differential equation then the last part of the proof
  above is unnecessary since the forms $\varphi^{\tau, e}$ defined in (\ref{phi-tau-e})  would be in $\Cyl^0\HH^*$ as
  $Y=dj$ for $j$ the immersion of $ M$ in $ \R^m$, and so  Proposition \ref{dense-cylindrical} 
 shows directly that the set 
  $$\{ f \varphi^{\tau, e} : f \hbox { is } C^{\infty}  \hbox{and cylindrical } , e \in \R ^m, 0 \leqslant
  \tau \le T \} $$ of elements of $\Cyl^0\HH^*$ is total in $ \D ^{q  , 1}  \HH^*$. 
\end{remark}

Combining Theorem \ref{theorem-cylindrical} with  Theorem \ref{th:divergence-21} on the continuity of ${\div}^q $ on $ \D^{q,1} \HH^*$, we see if Condition ($M_0$) holds,
$$\overline{(d^p)^*  |_{ {\Cyl^0}\HH^*}}=\overline{(d^p)^*  |_{\D^{q,1}\HH^*}}$$
 as operators on $L^q$ for $1\le q<\infty$. From this,
\begin{theorem}
\label{th-Sobolev-sp}  
 Suppose Condition ($M_0$) holds. Let $1<p<\infty$, then 
     
                                    $$W^{p, 1}=  {}^0 W^{p, 1}$$
 \end{theorem}
 Following Eberle \cite{Eberle-book}, consider  the space of bounded functions in 
 ${}^0 W^{2, 1}$ closed under the weak Sobolev norm, which we shall denote by ${}^0 W^{2,1}_\infty$.
 \begin{proposition}
  \label {Pro-Markov-uniqueness}   Suppose Condition ($M_0$) holds. Then   ${}^0 W^{2,1}_\infty=W^{2, 1}$. 
\end{proposition}
\begin{proof}
Take $f\in {}^0 W^{2,1}$ so $f\in W^{2,1}$ by Theorem  \ref{th-Sobolev-sp}. Then $g\equiv f\circ \I\in \D^{2,1}(\Omega)$ by Theorem \ref{th:weak_Sob}. Set 
$f_n={f\over 1+{1\over n} f^2}$ and $g_n=\I^*(f_n)$.
Then $|f_n|\le{1\over 2}$ is bounded. Now
 $$dg_n(\omega)(h)={d g(\omega)h\over 1+{1\over n} g^2} \cdot{n-g^2\over n+g^2}$$
is bounded and $g_n\in \D^{2,1}(\Omega)$  converges to $g$  in $\D^{2,1}$.
Consequently $f_n\in {}^{0}W^{2,1}_\infty$ and converges to $f$ in $W^{2,1}$ and so in ${}^0 W^{2, 1}$.  
Thus $ {}^0 W^{2,1}= {}^0 W^{2,1}_\infty$.
\end{proof}

Consider the symmetric diffusion operator $L=-d^*d$ on $L^2(\C_{x_0}M; \R)$ with domain $\Cyl$, the set of 
smooth cylindrical functions. It is called {\it Markov unique} if and only if there is only one symmetric 
sub-Markovian $C^0$ contraction semi-group $(P_t)$ on $L^2(\C_{x_0}M; \R)$ whose generator extends $L$.  
 Equivalently there is a unique extension of the corresponding Dirichlet form
   among the family of quasi-regular semi-Dirichlet forms. Markov uniqueness implies that there is at 
   most one reversible diffusion solving the corresponding martingale problem, c.f. Eberle \cite{Eberle-book}, though such results 
   go back to Takeda \cite{Takeda-uniqueness}.  We can apply a result of Eberle to our situation to obtain:
    
\begin{theorem} 
\label{th:Markov-unique-2}
 Suppose Condition ($M_0$) holds. The following are equivalent:
\begin{enumerate}
\item Markov uniqueness for $-d^*d$;
\item     
                                $W^{2, 1}=  \D^{2,1}$;
                                    \item $ \I^*[\D^{2,1}]=\D^{2,1}_{\filt}$.
                                    \end{enumerate}
 \end{theorem}
 
 \begin{proof}
It follows from Corollary \ref{co-Dirichlet-form} and  a result of Eberle (page 115, \cite{Eberle-book})
  that  Markov uniqueness is equivalent to $\D^{2, 1} = {}^0 W^{2, 1}_\infty$.
   Proposition \ref{Pro-Markov-uniqueness} then shows that  (1) and (2) are equivalent.
The equivalence of (2) and (3) is immediate from Theorem \ref{th:weak_Sob}.
\end{proof}

For discussion on when $ \I^*[\D^{2,1}]=\D^{2,1}_{\filt}$ holds see Remark \ref{remark-unique-5}.

\subsection{Cylindrical 1-forms are dense in $\D^{q,1}\HH^*$ }
\label{section-density}

The following lemma will be useful technically, leading to the proof of  Theorem \ref{theorem-cylindrical},  in view of the fact that we do not know if $\W^{2, 1} =\D ^{2, 1}$ on our path spaces. 
\begin{lemma}\label{le:density1}
  Let $G$ be a separable Hilbert space and $L$ a dense family of linear
  functionals on $G$. Suppose $f : \C_{x_0} M \to G$ satisfies
  \begin{enumerateroman}
    \item[(i)] $\mathcal{I}^{^{\ast}} ( f ) \in \D ^{p, 1}(\C_0\R^m; G)$  
    \item[(ii)] $l \circ f \in \D ^{p, 1}(\C_{x_0}M;\R)$ for all $l \in L$.
  \end{enumerateroman}
  Then $f \in \D ^{p, 1}(\C_{x_0}M;G)$.
\end{lemma}

\begin{proof}
For $f$ satisfying (i) and (ii) and $l\in G^*$ take $l_n\in L\to l$. Then
$$\I^*(l_nf)=l_n(\I^*f)\stackrel{\D^{p,1}}{\longrightarrow} l(\I^*f)= \I^*(l f).$$
By (\ref{Ito-map-norm}) $l_nf\stackrel{\D^{p,1}}{\longrightarrow}l f$ and so
we have
  $$(iii)  \qquad l \circ f \in \D ^{p, 1}(\C_{x_0}M;\R) \hbox{ for all  $l \in G^{\ast}$}.$$
    Take an orthonormal base $\{ g_n \}_n$ for $G$ and let $\Pi_n$ be the
  orthogonal projection onto the subspace spanned by its first $n$ elements, $n
  = 1, 2, \dots$. By (iii),  $\Pi_n \circ f \in \D ^{p, 1}(\C_{x_0}M;\R)$.

  Now as $n \to \infty$ we see $\Pi_n \circ f \to f$
  almost surely and also, if $\tilde d f$ is as in (\ref{weak-22}),
   for almost all $\sigma$, 
$$ \|d(\Pi_n \circ f )_{\sigma} - \tilde df_{\sigma} \|^2_{L _2} = \|
  (\tilde d f)_{\sigma}^{\ast} \Pi_n - (\tilde df)_{\sigma}^{\ast}
  \|_{L _2 ( G ;\HH_{\sigma} )}^2 
  = \sum_{p > n} \| \tilde  df_{\sigma} \circ (\tilde df_{\sigma} )^{\ast} ( g_p ) \|^2$$       
which converges to zero since $\tilde df_\sigma \in G\otimes \HH$.    
  From this we can apply the monotone convergence theorem to see that $\Pi_n
  \circ f \to f$ in the  $L^p$ graph norm and the result follows. 
\end{proof}

The next two lemmas are essentially `well known':

\begin{lemma}
\label{le:density2}
  Let $\tilde{\nabla}$ be a connection on the trivial $\R ^k$-bundle over $M$. Suppose $J : M \to \LL (\R^k; \R^k)$ is $C^{\infty}$ and define the
  `$J$-damped' parallel translation 
  $Z_t(\sigma): \{x_0\}\times \R ^k \to \{\sigma_t\}\times\R ^k$
  along $\sigma$ in  $\C_{x_0} M$ by
  \begin{equation}
  \label{cylindricals-1}\left\{
  \begin{array}{llll}
   \frac{\tilde{D}}{dt  } Z_t ( e )&=&J_{\sigma ( t )} \left(Z_t ( e ) \right), \quad &0
     \le t \le T \\
   Z_0 ( e ) &=&e, \quad &e \in \R ^k . 
   \end{array}\right.
   \end{equation}
    Then the principal part $ Z_t^P$ of $Z_t$ as a map from  $\C_{x_0}M \to \LL (\R ^k ;\R ^k)$ is in $\D ^{p, 1}, 
    1 \le p < \infty$, for each $t$.  Its H-derivative is an $\LL (\R ^k ;\R ^k )$-valued
  $H$-1-form: $h \mapsto dZ_t^P ( h )$ with
  \begin{equation}
  \label{density1}
    d Z_t^P ( h ) = A_t ( W_t^{-1}h_t ) + \int_0^t B_s (W_s^{-1} h_s ) ds  +
    \int_0^t C_s (W_s^{-1} \frac{\D}{ds } h ) ds, \quad h\in\HH 
  \end{equation}
where $$A_t, B_t, C_t \in L^p \left(\C_{x_0}M; \; \LL (T_{x_0}M ;\LL(\R^k;\R ^k ))\right), \; 1 \le p < \infty,\;
0 \le t \le T.$$ 
Furthermore $A_t ( \sigma )$, $B_t ( \sigma )$, and  
$C_t ( \sigma) : T_{\sigma ( t )} M \to \LL (\R ^k;\R ^k )$  are almost surely continuous in $t$ with 
$$\E \left( \sup_{0 \le t \le T} |A_t |_{\LL ( T_{x_0} M ; \LL (\R ^k ;\R ^k ) )}  \right)^p< \infty,$$ 
      $$\E \left( \sup_{0 \le t \le T} |B_t |_{\LL ( T_{x_0} M ; \LL (\R ^k ;\R ^k ) )}  \right)^p< \infty,$$ $$\E \left( \sup_{0 \le t \le T} |C_t |_{\LL ( T_{x_0} M ;\LL (\R ^k ;\R ^k ) )}  \right)^p< \infty.$$  
   \end{lemma}

\begin{proof}
  That $ Z_t \in \D ^{p, 1}$ for all $p$ and each $t$ is standard when
  we are using Brownian motion measure and the Levi-Civita connection, e.g. see
  the Appendix in Aida \cite{Aida97}, or L\'eandre \cite{Leandre93} and the proofs go
  over to our situation. The computation leading to (\ref{density1}) is also
  standard, going back, at least, to Bismut. See also Driver \cite{Driver92},  Cruzeiro-Malliavin \cite{Cruzeiro-Malliavin-96}. In particular if $h$ is an adapted
  $L^q$-$H$-vector field on $\C_{x_0} M$, from (\ref{cylindricals-1}), the covariant derivative
  $\tilde{\nabla} _h  Z_t$ satisfies the covariant Stratonovich equation
  \[ \tilde{D} \tilde{\nabla}_h  Z_t ( e ) = \left[ (\tilde{\nabla}_{h_t} J) (
      Z_t ( e ) ) + J_{\sigma ( t )} ( \tilde{\nabla}_{h_{}}  Z_t ( e ) )\right]
     dt   -\tilde{R} ( h_t, \circ d \sigma_t )  Z_t ( e ) .
     \]
 Thus
  \[ \tilde{\nabla}_h  Z_t ( e ) =  Z_t \int_0^t  Z_s^{- 1}\left \{
    ( \tilde{\nabla}_{h_s} J) (  Z_s ( e ) ) ds  - \tilde{R} ( h_s, \circ d
     \sigma_s )  Z_s ( e )\right\} . \]
  In fact by integration by parts, treating $\int_0^t Z _s^{- 1}
  \tilde{R} (W _s -, \circ d \sigma_s )  Z_s ( e )$ as an
  $L ( T_{x_0} M ;\R ^k )$-valued integral we see

  \begin{eqnarray*}
&& \int_0^t  Z_s^{- 1} \tilde{R} ( h_s, \circ d \sigma_s )  Z_s ( e ) =
     \left( \int_0^t  Z_s^{- 1} \tilde{R} (W _s -, \circ d \sigma_s )  Z_s
     ( e ) \right) W _t^{- 1} h_t\\
  &&\hskip 90pt  - \int^t_0 \left( \int_0^{\tau}  Z_s^{- 1} \tilde{R}
     (W _s -, \circ d \sigma_s )  Z_s ( e ) \right)\left(
     W _{\tau}^{- 1} \frac{\D }{d \tau} h_{\tau} \right) d \tau.
      \end{eqnarray*}   
      So  our expression for $\tilde{\nabla}_h  Z_t ( e )$ is `tensorial' in $\HH$
      and so holds for arbitrary elements in $\HH$.
  
  To obtain (\ref{density1}) we observe that the principal part of $\tilde{\nabla}_h  Z_t ( - ) $
  and $dZ_t^P ( h )$ only differ by  $ \tilde{\Gamma} ( \sigma   ( t ) ) ( h_t ) (  Z_t^{p} )$
 where $\tilde \Gamma$ is the Christoffel symbol $\tilde{\Gamma}$.
    Finally the required estimates come from Proposition \ref{norm-TI} and the analogous  observation that $|\sup_t Z_t|$ and $|\sup_t Z_t^{-1}|$  
    are in $L^p$ for all $p$.
\end{proof}

\begin{lemma}
\label{le:cylindricals-5}
  Let $Z$ be the `$J$-damped' parallel translation from Lemma \ref{le:density2}. Then
  for $p > 1$, and $\epsilon\in (0, p-1)$, the map $f \longmapsto Z( f )$
  gives a continuous linear map
  \[ \D ^{p, 1} \left( \C_{x_0} M ; L^2 \left( [ 0, T ] ;\R ^k \right)\right)
     \to \D ^{p - \varepsilon, 1} \left( \C_{x_0} M ; L^2 \left( [ 0,
     T ] ;\R ^k \right)\right) . \]
\end{lemma}

\begin{proof}
  Take $f$ in $ \D ^{p, 1} \left( \C_{x_0} M ; L^2 ( [ 0, T ] ; \R ^k )\right)$. It is standard and easy to
   see using Lemma \ref{le:density2} and weak
  differentiability, \eg Sugita \cite{Sugita85} Corollary 2.1, that
  $\I^* ( Z ( f) ) \in \D ^{p - \varepsilon, 1} \left(\C_0 \R ^m ; L^2 ( [ 0, T ] ;\R ^k ) \right)$. 
  Let $E_n$ be the subspace of $L^2 ( [ 0, T ] ; \R ^k )$ spanned by the polynomials of
  degree less than or equal to $n$, $n = 1, 2, \ldots$ , with $\pi_n (Z(f))$ the
  orthogonal projection of $Z(f)$ into $E_n$. Since the evaluations $\{\ev_t, 0\le t\le T\}$
    span a dense linear subspace of linear functionals on each $E_{n}$, we
  can apply the previous lemmas to see that
   $\pi_n(Z( f ))  \in \D ^{p - \varepsilon, 1} ( \C_{x_0}M ; L^2 ( [ 0, T ] ;\R ^k ) )$
    for each $n = 1, 2, \ldots$. However as in the proof
  of Lemma \ref{le:density1} we see that $\I^* (  \pi_n(Z( f)) )\to \I^* ( Z( f) )$ in the $L^{p - \varepsilon}$ graph norm. Therefore  
  $Z( f) \in \D ^{p - \varepsilon, 1} \left( \C_{x_0} M ; L^2 ( [ 0, T ] ;\R ^k ) \right) .$
  Continuity follows as usual from continuity into $L^{p - \varepsilon}$ and
  the closed graph theorem.
\end{proof}

For cases including non-elliptic diffusion measures take a Riemannian metric $\langle, \rangle^{\prime}$
on $TM$ extending that of $E$ with a metric connection $\nabla^1$ and adjoint ${\nabla^1}^\prime$
extending $\nabla$ and $\nabla^\prime$ as in \S\ref{subsection-2.1}.  Take a surjective vector bundle map $\tilde X: \underline \R^{\tilde m}\to TM$ for some $\tilde m$,  inducing the metric $\langle, \rangle^{\prime}$, which extends $X$
if $\R^m$ is considered as a subspace of $\R^{\tilde m}$.
For $x\in M$ let $\tilde Y_x: T_xM\to \R^{\tilde m}$ be the usual right inverse of $\tilde X(x)$. Take a connection $\tilde \nabla$ on $\R^{\tilde m}$ conjugate to ${\nabla^1}^\prime$ on $[\ker \tilde X]^\perp$ and arbitrary on $\ker \tilde X$. In the elliptic case we could take $\tilde X=X$, $\tilde Y=Y$ and if moreover $\nabla$ is the Levi-Civita connection then $\tilde \nabla$ could be the standard metric connection on $\R^m$ induced by $K$, as in Elworthy-LeJan-Li \cite{Elworthy-LeJan-Li-book}.

In order to prove the density of cylindrical forms in $\D ^{q, 1}\HH^*$ fix $\tau \in [ 0, T ]$ and $e \in \R ^{\tilde m}$. Let $\varphi^{\tau, e}$ be the
cylindrical one-form given by
\begin{equation}
\label{phi-tau-e}
 \varphi^{\tau, e} ( V ) 
= \langle \tilde Y_{\sigma ( \tau )} ( V_{\tau} ), e \rangle, \qquad V \in T_{\sigma} \C_{x_0} M.
\end{equation}
In fact we will show that the set $\{f\phi^{\tau,e}:f \in \Cyl, e \in \R ^{\tilde m}, 0 \le \tau \le T \}$
  is total in $\D ^{q, 1}\HH^*$.
Set        
\begin{equation}U^{\tau, e} = ( \varphi^{\tau, e} )^{^{\#}}.
\end{equation}
Using the fact that
\[ W_{\tau} \int_0^{\tau} W_s^{- 1}  \frac{\D }{ds } V_s ds  = V \tau \]
we see that
\begin{equation}
\label{density2-2}
  \frac{\D }{ds } U^{\tau, e}_s =\chi_{[ 0, \tau ]} ( s ) \Pi_s(W_s^{-1})^*(W_\tau)^*\tilde X(\sigma(\tau))(e),
  \end{equation}
  where $\Pi_s=\Pi_s(\sigma):  T_{\sigma(s)}M  \to  E_{\sigma(s)}$ is the orthogonal projection and $W_s^*$ and ${W_s^{-1}}^*$  the adjoints using the extended metric $\langle, \rangle ^\prime$.

Now define $Z_t \equiv Z_t ( \sigma ) :\R ^{\tilde m}  \to\R ^{\tilde m}$ to be the damped parallel translation on the trivial $\R ^{\tilde m}$-bundle over $M$ along  $\sigma \in \C_{x_0} M$ given by
\[ Z_t ( a ) = \left\{
\begin{array}{ll}
      \tilde Y_{\sigma ( t )} W_t \tilde X ( x_0 ) a,  &\hbox{if } a \in  [\ker (\tilde X (x_0 ))]^{\perp}  \\
      \tilde{\parals_t} ( a ), &\hbox{if } a \in \ker (\tilde X ( x_0 ) ).
\end{array}
\right.
\]

By (\ref{covariant-3}), $Z_t$ solves (\ref{cylindricals-1}) for
$$J(x)(a)=\tilde Y(x)\left(-{1\over 2}Ric^{\#}_x +\nabla^1_-A\right) \tilde X(x)a.$$
Let $Z_t^*: \R ^{\tilde m}\to \R ^{\tilde m}$ be the usual adjoint of $Z_t$ and similarly for
$( Z^{- 1}_t )^{\ast}$. Set                      
$$\Lambda^{\tau, e}_s  = \chi_{[ 0, \tau ]} ( s ) ( Z_s^{- 1} )^{\ast}
Z^{\ast}_{\tau} e,$$                                           
then
\begin{equation}
\label{density2}
  \frac{\D }{ds } U^{\tau, e}_s = \Pi_s\tilde X ( \sigma ( s ) )
  \Lambda^{\tau, e}_s.
\end{equation}

\begin{lemma}
\label{le:density4}
Set \begin{equation}
\Xi = \{ f \Lambda^{\tau, e}_\cdot : f \in \Cyl, e \in \R ^{\tilde m}, 0 \le \tau \le T \}.
\end{equation}
  The sets $Z^*_\cdot [ \Xi ]$ and $\Xi$ are both total in 
  $\D ^{q, 1}\left( \C_{x_0} M ; L^2 ( [ 0, T ] ;\R ^{\tilde m} )\right)$.
  \end{lemma}

  \begin{proof}
  Since every element of the Haar
    basis, $\{ E^j \}_j$ say, of $L^2 ( [ 0, T ] ;\R^{\tilde m} )$ has the form 
    $\chi_{[ 0, \tau_2 ]} ( \cdot) e  - \chi_{[ 0, \tau_1 ]} (\cdot) e $ some $0
    \le \tau_1 \le \tau_2 \le T, e \in \R ^{\tilde m}$, using the definition of $\D ^{q, 1}$ we have
    $$\mathcal{T} := \left\{ f\chi_{[0,\tau]}(\cdot)e:  0\le \tau\le T, e\in \R^{\tilde m}, f \in  \Cyl\right\}$$
     is total in  $\D ^{q,1} \left( \C_{x_0} M ; L^2 ( [ 0, T ] ;\R^{\tilde m} )\right)$ for 
    $1 \le q < \infty$.
    On the other hand, by definition,
    \begin{equation}
    \label{density3}
      Z^{\ast}_\cdot [ \Xi ] = \left\{  f \chi_{[ 0, \tau ]} ( \cdot) Z^{\ast}_{\tau}    e_\cdot : f \in \Cyl, \;  e \in
      \R^{\tilde m},  \; 0 \le \tau \le T\right \},
    \end{equation}
    a subset of $\D ^{q,1} ( \C_{x_0} M ; L^2 ( [ 0, T ] ;\R^{\tilde m} )$ from Lemma \ref{le:density1}.
 
      For each $e\in \R^{\tilde m}$ and  $\tau \in [ 0, T ]$ define   $e^{\tau} : \C_{x_0} M \to \R^{\tilde m}$ by
    \[ e^{\tau} ( \sigma ) = \left(Z^{\ast}_{\tau} ( \sigma )\right)^{- 1} e. \]
  Note that $e^{\tau} \in \D ^{p, 1}, 1 \le p < \infty$, by
   Lemma \ref{le:density1} (applied to $( Z^{-1}_\cdot)^{\ast}$). Consequently given
    $\varepsilon > 0$ and $q_0 > 1$, there exist $g^{\tau} \in   \Cyl( C_{x_0} M ;\R^{\tilde m} )$ with
   $$	 \| g^{\tau} - e^{\tau} \|_{\D ^{q, 1}} < \varepsilon, \quad  1\le q \le q_0 .$$
    Then, for $1\le q<q_0$,
    \begin{eqnarray*}
    &&\left\| \chi_{[ 0, \tau ]} (\cdot) e  - \chi_{[ 0, \tau ]} (\cdot )
       Z^{\ast}_{\tau} g^{\tau} \right\|_{\D ^{q, 1} ( \C_{x_0} M ; L^2 (  [ 0, T ] ; \R^{\tilde m} ) )}\\
       &=& \| \chi_{[ 0, \tau ]} (\cdot) (
       Z^{\ast}_{\tau} e^{\tau} - Z^{\ast}_{\tau} g^{\tau})
       \|_{ \D ^{q, 1} ( \C_{x_0} M ; L^2 ( [ 0, T ] ; \R^{\tilde m} ) )} 	
      \le C \varepsilon \| Z_{\tau} \|_{\D ^{r, 1}} 
       	\end{eqnarray*}
   for sufficiently large $r$ and some constant $C$.
    
    Thus each
     $\chi_{[ 0, \tau ]} ( \cdot) e   \in \overline{\Span}^{q, 1} Z^*_\cdot [ \Xi ]$,  $ 1 \le q <\infty$
     and ${\mathcal T}\subset  \overline{\Span}^{q, 1} Z^*_\cdot [ \Xi ]$. Consequently for each $1 \le q < \infty$,
    $$ \overline{\Span}^{q, 1} Z^{\ast}_. [ \Xi ] =\D ^{q, 1}\left (\C_{x_0} M ; L^2 ( [ 0, T ] ; \R^{\tilde m} )\right),  $$
    as required for the first assertion.     For the second,  take $fE^j \in \Span\mathcal{T}$ and $q<p<\infty$. By 
   Lemma \ref{le:cylindricals-5} $Z^{\ast}_\cdot( fE^j ) \in \D ^{p, 1}$. Since by the first assertion $\D^{p,1}=\Span Z^*[\Xi]$, there is a sequence $\{ S^n \}_n$ in $\Span[\Xi]$ with
    $Z^{\ast}_\cdot ( S^n_\cdot) \to Z^{\ast} ( fE^j )$ in
    $\D ^{p, 1}$. Using Lemma \ref{le:cylindricals-5}, applied to ($Z^{\ast})^{- 1}$, we 
    see $S^n \to fE^j$ in $\D ^{q, 1}$, which implies the result by the totality of $\mathcal{T}.$
  \end{proof}

\begin{proposition}
\label{dense-cylindrical}
Smooth cylindrical 1-forms form a dense subspace of $\D^{q,1}\HH^*$.
\end{proposition}
\begin{proof}
By construction of $\tilde X$, if $\pi: \R^{\tilde m}\to \R^{\tilde m}$ is the projection map then $\Pi\tilde X=X\circ \pi$ and so
$$\D^{q,1}\HH= \left\{\W\Pi(\tilde X(h) )\; |\; h\in \D^{q,1}(\C_{x_0}M;L^2([0,T];\R^{\tilde m})\right\}.$$
By Lemma \ref{le:density4},  $\{\W_\cdot \Pi\tilde X(f\Lambda_\cdot^{\tau, e}): f \in \Cyl, e \in \R^{\tilde m}, 0 \le \tau \le T \}$ is total in $\D^{q,1}\HH$. Finally note that by (\ref{density2}), $\W_\cdot \Pi\tilde X(f\Lambda_\cdot^{\tau, e})=U_s^{\tau,e}$ and so the
set  $$\{f\phi^{\tau,e}:f \in \Cyl, e \in \R^{\tilde m}, 0 \le \tau \le T \}$$  is total in $\D^{q,1}\HH^*$.

To see that every cylindrical one-form gives an element of $\cap_{q<\infty}\D^{q,1}\HH$ suppose $\phi$ is one, given by
$$\phi_\sigma(v)=\Phi_{(\sigma(t_1), \dots, \sigma(t_k))}\left(v_{t_1}, \dots, v_{t_k}\right), \qquad v\in \HH_\sigma$$
where $\Phi$ is a smooth 1-form on $\stackrel {k} {\overbrace{M\times\dots \times M}}$.
We must show $$\phi\circ \overline{T\I}\in \D^{q,1}(\C_{x_0}M; H^*), \quad 1<q<\infty.$$
 For this write
\begin{eqnarray*}
\phi\circ T\I_\sigma(h)
&=&\Phi\left(\tilde X(\sigma(t_1))-, \dots, \tilde X(\sigma(t_k))-\right)\\
&& \left(\int_0^{t_1}\tilde Y_{\sigma(t_1)} W^s_{t_1} \tilde X(\sigma(s))\dot h_s ds, \dots,
\int_0^{t_k}\tilde Y_{\sigma(t_k)} W^s_{t_k} \tilde X(\sigma(s))\dot h_s ds  \right)
\end{eqnarray*}
for $h\in H\subset L_0^{2,1}\left([0,T]; \R^{\tilde m} \right)$. The first part of the right hand side is in
$$\cap_{q<\infty} \D^{q,1} \left(\C_{x_0}M, \left(\R^{\tilde m}\times \dots \times \R^{\tilde m}\right)^*\right),$$
being cylindrical. The second can be written as
$$h\mapsto \left (\int_0^{t_1} Z_{t_1}Z_s^{-1}K^\perp (\sigma(s))\dot h_s ds,
\dots, \int_0^{t_k} Z_{t_k}Z_s^{-1}K^\perp (\sigma(s))\dot h_s ds\right)$$
which by Lemma \ref{le:cylindricals-5} and \ref{le:multiplication} gives a continuous linear map 
$$H\to \D^{q,1} \left(\C_{x_0}M, (\R^{\tilde m}\times \dots \times \R^{\tilde m}\right).$$
 Thus $\phi\in \D^{q,1}$.
\end{proof}

\section{On uniqueness of $d$}
\label{d-uniqueness}
\subsection{A weak uniqueness result on the Gross derivative operator $d$}
\label{se:weak-uniqueness}

 If Condition $(M_0)$ holds the map $\mathcal{I}$* sends $\D^{2, 1} ( \C_{x_0}  M ;\R)$ to $\D^{2, 1}_{\F^{x_0}}(\Omega;\R)$ with closed range, by Corollary \ref{co:closed-range}. We investigate the question whether $\I^*[\D^{2, 1} ( \C_{x_0}  M ;\R)]=\D^{2, 1}_{\F^{x_0}}(\Omega;\R)$. We proceed using chaos expansions.  An $L^2$ real valued function $f$  on the Wiener space has a chaos expansion
$$f=\sum_{k=0}^\infty I_k(\alpha_k),   $$
for $\I_0=\id$ on constants, $\alpha_0=\E (f)$, 
$I_1(\alpha_1)=\int_0^T\langle \alpha_1(t_1), dB_{t_1}\rangle$, and for $k>1$,
\begin{equation}
I_k(\alpha _k)=k!\int_0^T\int_0^{t_k}\dots \int_0^{t_2}\left\langle \alpha
_k\left( t_1,\dots ,t_k\right) ,dB_{t_1}\otimes \dots \otimes
dB_{t_k}\right\rangle _{\otimes \R ^m},
\end{equation}
an  iterated It\^{o} integral. Here $\alpha_k$ is considered to be an element of
$L^2\left(\Delta^k; \R^{m}\otimes\dots\otimes \R^{m}\right)$ for $\Delta^k=\{(t_1, \dots, t_k), 0\le t_1\le t_2\le \dots t_k\le T\}$.

Let $R_n$ be the remainder term such that
$$f=\sum_{k=0}^n I_k(\alpha_k)+R_n .$$
Then $f\in \D^{2,1}(\Omega;\R)$ if and only if $R_n\to 0$ in $L^2$ and $\|d R_n\|\to 0$ in $L^2$ (\eg see Nualart \cite{Nualart-book}, Proposition 1.2.2).
\medskip

Let  $L^2_{\F^{x_0}}\equiv L^2_{\F^{x_0}}(\Omega, \R)$ be the closed subspace of  $L^2(\Omega;\R)$  whose elements are $\F^{x_0}$ measurable. 
\begin{lemma}
\label{le:conditional_chaos1}
If $f\in\D^{2,1}_{\F^{x_0}}(\Omega, \R)$ then the remainder term of its $L^2$ chaos expansion has the following form
\begin{equation}\label{chaos-1}
R_k=\int_0^T\int_0^{s_{k+1}}\dots \int_0^{s_2} \left\langle a_{k+1}(s_1, \dots, s_{k+1}), dB_{s_1}\otimes \dots \otimes dB_{s_{k+1}}\right\rangle
\end{equation}
for $a_{k+1}\in L^2\left(\C_0\R^m\times \Delta^{k+1}; \R^m\otimes \dots \otimes \R^m\right)$ such that
\begin{enumerate}
\item[(i)] 
Each $a_{k+1}(s_1, \dots, s_{k+1})$ is $\F^{x_0}_{s_1}$ measurable;
\item[(ii)]
$a_{k+1}(s_1,\dots, s_{k+1})=\left(K^\perp(x_{s_1}) \otimes \1\dots\otimes \1\right)\alpha^{k+1}(s_1, \dots s_{k+1})$, where $K^\perp(x): \R^m\to [\ker X(x)]^\perp$ is the  orthogonal projection. 
\item[(iii)]
Furthermore 
$a_{k+1}\in \D^{2,1}(\C_0\R^m; L^2(\Delta^{k+1}; \R^m\otimes\dots\otimes \R^m))$ and
\begin{equation}\label{chaos-2}
\|dR_{k}\|^2_{L^2}=(k+1)\|a_{k+1}\|_{L^2}^2+\|da_{k+1}\|^2_{L^2}.
\end{equation}
\end{enumerate}
\end{lemma}

\begin{proof} 
By the integral representation theorem there exists an $L^2$, $\F_*$-adapted process $a_1:  [0, T] \times \C_0\R^m \to \R^m$ with 
\begin{equation}\label{conditional_chaos}
f=\E f+\int_0^T \langle a_1(s), dB_s\rangle.
\end{equation}
Since $f=\E\{f|\F^{x_0}\}$ we see, using (\ref{noise-decomposition}), 
$$f=\E f+\int_0^T  \langle K^\perp (x_s)\E\{a_1(s)|\F^{x_0}\}, dB_s\rangle.$$
By the uniqueness of the integral representation (\ref{conditional_chaos})
$$K^\perp(x_s)\E\{a_1(s)|\F^{x_0}\}=a_1(s).$$
Thus $a_1(s)$ is $\F^{x_0}_s$ measurable and (i) and (ii) hold for $k=0$. Suppose for induction (i) and (ii) hold for $k-1$ some $k\ge 1$. Apply the integral representation theorem to the $\F^{x_0}_{s_2} $ measurable function $a_k(s_2,\dots, s_{k+1})$ to see
$$a_k(s_2,\dots, s_k)=\E \left(a_k(s_2,\dots, s_k)\right)+\int_0^{s_2} \left\langle
a_{k+1}(s_1, \dots, s_{k+1}), dB_{s_1}\right\rangle,$$
where $a_{k+1}(s_1, \dots, s_{k+1})$ is $\F^{x_0}_{s_1}$ measurable with  
$$\left(K^\perp(x_{s_1}) \otimes \1 \dots\otimes \1\right)
a_{k+1}(s_1, \dots, s_{k+1})=a_{k+1}(s_1, \dots, s_{k+1}).$$
By the uniqueness of the chaos expansion with remainder terms, we see that $R_{k}$ together with $a_{k+1}$ satisfies (i) and (ii).

To prove part (iii), we apply a result of Pardoux-Peng \cite{Pardoux-Peng90}, see Lemma 1.3.4 in Nualart \cite{Nualart-book}: if 
$$X=\int_0^T \langle u_s, dB_s\rangle$$
where $u$ is adapted and square integrable then $X\in \D^{2,1}$ implies $u_\cdot\in \D^{2,1}$. Moreover \begin{equation} \|d X\|^2_{L^2}=\int_0^T \E|u_s|^2 ds+\int_0^T \|d u_s\|^2_{L^2} \; ds.
\end{equation}
Apply this to $R_k$ iteratively,  to deduce that
$a_{k+1}(s_1, \dots, s_{k+1})$ belongs to  $\D^{2,1}$ and 

\begin{eqnarray*}
&&\|dR_k\|^2_{L^2}\\
&& =\int_0^T \E \left|
\int_0^{s_{k+1}}\dots \int_0^{s_2}     \left\langle a_{k+1}(s_1, \dots, s_{k+1}), dB_{s_1}\otimes \dots \otimes dB_{s_{k}} \right\rangle \right|^2 ds_{k+1}\\
&&+\int_0^T  \left \| d  \int_0^{s_{k+1}}\int_0^{s_k}\dots \int_0^{s_2}     \left\langle a_{k+1}(s_1, \dots, s_{k+1}), dB_{s_1}\otimes \dots \otimes dB_{s_{k}} \right\rangle  \right\|^2_{L^2} ds_{k+1}\\
&&=\|a_{k+1}\|_{L^2}^2+\int_0^T\int_0^{s_{k+1}} \E \left|
\int_0^{s_k}\dots \int_0^{s_2}     \left\langle a_{k+1}(s_1, \dots, s_{k+1}), dB_{s_1}\otimes \dots \otimes dB_{s_{k-1}} \right\rangle \right|^2 ds_k ds_{k+1}\\
&&+\int_0^T\int_0^{s_{k+1}}  \left\|d \int_0^{s_k}  \int_0^{s_{k-1}} \dots \int_0^{s_2}     \left\langle a_{k+1}(s_1, \dots, s_{k+1}), dB_{s_1}\otimes \dots \otimes dB_{s_{k-1}} \right\rangle  \right \|^2  _{L^2}ds_{k}ds_{k+1}\\
&&=\dots,
\end{eqnarray*}
giving (\ref{chaos-2}).

\end{proof}

Note that  Lemma \ref{le:conditional_chaos1} (iii) holds if $\F^{x_0}$ is replace by $\F_T$ everywhere in the statement. This shows that a function $f\in \D^{2,1}$ if and only if
\begin{equation}\label{chaos-4}
k\|a_k\|_{L^2}^2\to 0 \hbox{ and } \|da_{k}\|_{L^2}^2\to 0.
\end{equation}

\begin{lemma}
		\label{le:conditional_chaos2}
Let $f\in L^2(\Omega;\R)$ with chaos expansion $f=\sum_{k=0}^\infty I_k(\alpha_k)$.
		Set $\bar J_k(\alpha_k)(\sigma)=\E\left\{\I_k(\alpha_k)\big| x_\cdot=\sigma\right\}$.
		Then
$\bar J_k(\alpha_k)$ is in $\D^{2,1}(\C_{x_0}M)$. Consequently if $f\in\D^{2,1}_{\F^{x_0}}$ 
then each $f-\E\{R_k|\F^{x_0}\}$ is in $\I^*[\D^{2,1}(\C_{x_0}M;\R)]$. 

	\end{lemma}

\begin{proof} 
		For $0\le t\le T$ define \begin{equation}
I_{k,t}(\alpha _k):=k!\int_0^t\int_0^{t_k}\dots \int_0^{t_2}\left\langle
\alpha _k\left( t_1,\dots ,t_k\right) ,dB_{t_1}\otimes \dots \otimes
dB_{t_k}\right\rangle _{\otimes \R ^m}. \label{multiple-integrals}
\end{equation}
Thus
\begin{equation}
I_{k,t}(\alpha _k)=k\int_0^t\Big\langle I_{k-1,t_k}(\alpha _k\left( \dots
,t_k\right)) ,dB_{t_k}\Big\rangle _{\R ^m}  \label{Iterated-integrals}
\end{equation}
considering $s\mapsto I_{k-1,t_k}(\alpha _k(\dots ,s))$ as a random element in $%
L^2\left( [0,T]; \R ^m\right) $, or more precisely
\begin{equation}
I_{k,t}(\alpha _k)=k\int_0^t\left\langle \left( I_{k-1,t_k}\otimes id\right)
\left(\alpha _k ( \dots , t_k)\right) ,dB_{t_k}\right\rangle _{\R ^m}
\end{equation}
if we identified $L^2\left( \Omega ;\R \right)
\otimes \R ^m$ with $L^2\left( \Omega ;\R ^m\right) $. 
Set
\begin{eqnarray*}
J_{k,t}\left( \alpha _k\right)  &=&\E \left\{ I_{k,t}(\alpha _k)\;
|\quad \mathcal{F}^{x_0}\right\} ,\; k=0,1,2,\dots  \\
J_k\left( \alpha _k\right)  &=&J_{k,T}\left( \alpha _k\right) .
\end{eqnarray*}
Then inductively
$$J_{k,t}\left( \alpha _k\right) =k\int_0^t\left\langle J_{k-1,t_k}\left(
\alpha _k\left( \dots ,t_k\right) \right) , K^{\perp
}(x_{t_k})dB_{t_k}\right\rangle _{\R ^m}$$

Set
$\breve B_t=\int_0^t {\parals_s}^{-1}X(x_s)dB_s$, the martingale part of the  anti-development of $x_\cdot$ to see
 \begin{eqnarray}\nonumber
K^\perp(x_t)dB_t&=&Y(x_t)X(x_t)dB_t=Y(x_t){\parals_t}(x_\cdot)d\breve B_t;\\
 \bar J_{k,t}\left( \alpha _k\right) (\sigma)&=&k\int_0^t\left\langle \bar J_{k-1,t_k}\left(
\alpha _k\left( \dots ,t_k\right) \right) (\sigma), Y(x_{t_k}){\parals_{t_k}}(\sigma_\cdot)d\breve B_{t_k}
\right\rangle _{\R ^m}. 
 \label{bar-J}
\end{eqnarray}
The fact that $\bar J_{k,t}(\alpha_k)$ is in $\Dom(d)$ is essentially standard \eg\
see Cruzeiro-Malliavin  \cite{Cruzeiro-Malliavin-96} or the Appendix in Aida
\cite{Aida97}. For a gradient stochastic differential equation
(\ref{sde}) determined by an isometric $j: M \to \R^m$ it is especially clear since then 
$K^\perp (x_t)dB_t$ can be replaced by $d\tilde x_t-{1\over  2}\Delta j(x_t)dt$ for $\tilde x_t=j(x_t)\in \R^m$. 

Finally just observe that for $f\in \D^{2,1}_{\F^{x_0}}$,
$$f-\E\{R_k| \F^{x_0}\}=\sum_{j=1}^m \I^*\left(\bar J_k(\alpha_k)\right)$$
which is in $\I^*[\D^{2,1}(\C_{x_0}M)]$.   

\end{proof}

\begin{proposition}
\label{pr:chaos-expansion}
Suppose Condition $(M_0)$ holds.
Then $f\in Dom(\Delta)\cap L^2_{\F^{x_0}}(\Omega, \R)$ implies that $f\in \I^*[\D^{2,1}(\C_{x_0}M;\R)]$. 
\end{proposition}

\begin{proof}
Take $f\in \D^{2,1}(\Omega, \R)$. By Corollary \ref{co:closed-range}, $\I^*[\D^{2,1}(\C_{x_0}M)]$ is closed in $\D^{2,1}(\C_0\R^m;\R)$ 
and so by Lemma \ref{le:conditional_chaos2} to show $f\in \I^*[\D^{2,1}(\C_{x_0}M;\R)$ we only need to demonstrate that $\|d\left(\E\{R_k|\F^{x_0}\}\right)\|_{L^2}\to 0$.
Observe that, as iterated integrals,
\begin{eqnarray*}
\E \{R_k|\F^{x_0}\}
&=&\int_0^T\int_0^{s_{k+1}}\dots \int_0^{s_2}\\
&& \left\langle a_{k+1}(s_1, \dots, s_{k+1}), K^\perp(x_{s_1})dB_{s_1}\otimes \dots \otimes K^\perp (x_{s_{k+1}})dB_{s_{k+1}}     \right\rangle
\end{eqnarray*}
and for $h\in H$,
\begin{eqnarray*}
&&d\left(\E\{R_k|\F^{x_0}\}\right)(h)\\
&& =\sum_{l=1}^{k+1} 
\int_0^T \int_0^{s_{k+1}}\dots \int_0^{s_2}
\left\langle \, a_{k+1}(s_1, \dots, s_{k+1}), \right.\\
&&\quad\qquad\left. K^\perp(x_{s_1}) dB_{s_1}\otimes \dots \otimes d(K^\perp(\I_{s_l})) (h)dB_{s_l} \otimes\dots K^\perp (x_{s_{k+1}})dB_{s_{k+1}} \; \right\rangle\\
&&+\sum_{l=1}^{k+1}\int_0^T \int_0^{s_{k+1}}\dots \int_0^{s_2}
\left \langle  a_{k+1}(s_1, \dots, s_{k+1}), \right.\\
&&\quad\qquad\qquad\qquad \left. K^\perp(x_{s_1})dB_{s_1}\otimes \dots  K^\perp(x_{s_l})\dot h_{s_l}\otimes\dots K^\perp (x_{s_{k+1}}) dB_{s_{k+1}}\right\rangle\\
&&+\int_0^T \int_0^{s_{k+1}} \dots\int_0^{s_2}
\langle  d(a_{k+1}(s_1, \dots, s_{k+1}))(h), K^\perp(x_{s_1})dB_{s_1}\otimes \dots K^\perp (x_{s_{k+1}})dB_{s_{k+1}}  \rangle\\
&=&A_1(h)+A_2(h).
\end{eqnarray*}
where
\begin{eqnarray*}
A_1(h)&=&\sum_{l=1}^{k+1} 
\int_0^T \int_0^{s_{k+1}}\dots \int_0^{s_2} \left\langle a_{k+1}(s_1, \dots, s_{k+1}), \right.\\
&& \left. K^\perp(x_{s_1}) dB_{s_1}\otimes \dots \otimes d(K^\perp\circ \I_{s_l}) (h)dB_{s_l} \otimes\dots K^\perp (x_{s_{k+1}})dB_{s_{k+1}}\right\rangle
\end{eqnarray*}
On the other hand 
 \begin{eqnarray*}
 \E\{d R_k(h) |\F^{x_0}\}
&&= \sum_{l=1}^{k+1}\int_0^T \int_0^{s_{k+1}}\dots \int_0^{s_2}
\left \langle  a_{k+1}(s_1, \dots, s_{k+1}), \right.\\
&&\qquad \qquad \left. K^\perp(x_{s_1})dB_{s_1}\otimes \dots  \dot h_{s_l}\otimes\dots K^\perp (x_{s_{k+1}}) dB_{s_{k+1}}\right\rangle\\
&&+\int_0^T\int_0^{s_{k+1}}\dots \int_0^{s_2}\left\langle  \E \{d(a_{k+1}(s_1, \dots, s_{k+1}))(h) |\F^{x_0}\}, 
\right.\\
&& \qquad \left. K^\perp(x_{s_1})dB_{s_1}\otimes \dots K^\perp (x_{s_{k+1}})dB_{s_{k+1}}  \right\rangle.
\end{eqnarray*}
Thus
 $$d\left(\E\{R_k|\F^{x_0}\}\right)(h)=A_1(h)+C_2(h)+\E\{d R_k({\mathcal K}^\perp h) |\F^{x_0}\},$$ 
 where ${\mathcal K}$ is as defined by (\ref{divergence-3})
 and
\begin{eqnarray*}
C_2(h)&=&\int_0^T \int_0^{s_{k+1}}\dots \int_0^{s_2} 
\left\langle \; d(a_{k+1}(s_1, \dots, s_{k+1}))(h) \right.\\
&&  \left.  - \E \{d(a_{k+1}(s_1, \dots, s_{k+1}))({\mathcal K}^\perp h) |\F^{x_0}\}, \;K^\perp(x_{s_1})dB_{s_1}\otimes \dots K^\perp (x_{s_{k+1}})dB_{s_{k+1}} \right \rangle.
\end{eqnarray*}
However  
$$\|C_2\|^2_{L^2}\le 4 \int_0^T\int_0^{s_{k+1}}\dots \int_0^{s_2}\|d a_{k+1}(s_1,\dots, s_{k+1})\|^2 _{L^2}ds_1\dots ds_{k+1}=4\|da_{k+1}\|^2_{L^2}$$ 
and $\|\E\{dR_k(-)|\F^{x_0}\}\|_{L^2}\le \|dR_k(-)\|_{L^2}$ 
giving
$$ \|d\left(\E\{R_k|\F^{x_0}\}\right)\|^2_{L^2} \le 6\|A_1\|^2_{L^2} +24\|da_{k+1}\|^2_{L^2}+6|dR_k(-)\|^2_{L^2}.$$
For $f\in\D^{2,1}_{\F^{x_0}}$,  $\|dR_k(-)\|^2_{L^2}\to 0$ and $\|da_{k+1}\|^2_{L^2}\to 0$ by (\ref{chaos-2}). Thus we only need to show that $\|A_1\|_{L^2}\to 0$.
Now
$$\|A_1\|^2_{L^2} \le C k^2 \E\int_0^T\int_0^{s_{k+1}}\dots \int_0^{s_2}\|a^{k+1}(s_1,\dots, s_{k+1})\|^2 ds_1\dots ds_{k+1}\le C k^2\|R_k\|_{L^2}^2,$$
by Proposition \ref{proposition-3.2},
which converges to zero if $f$ belongs to $\Dom(\Delta)$.
\end{proof}

\begin{remark}
\label{remark-unique-5}
If we can show that $\|A_1\|^2_{L^2} \to 0$ without the condition $f\in \Dom(\Delta)$ we would have shown that $\I^*[\D^{2,1}(\C_{x_0}M;\R)]=\D^{2,1}_{\F^{x_0}}(\Omega;\R)$. This convergence should hold, though we do not have a proof,  as the similar term 
\begin{eqnarray*}
\|A_2(h)\|_{L^2}&\equiv& \Big\|\sum_{l=1}^{k+1}\int_0^T \int_0^{s_{k+1}}\dots \int_0^{s_2}
\big \langle  a_{k+1}(s_1, \dots, s_{k+1}), \\
&& K^\perp(x_{s_1})dB_{s_1}\otimes \dots  K^\perp(x_{s_l})\dot h_{s_l}\otimes\dots K^\perp (x_{s_{k+1}}) dB_{s_{k+1}}\big\rangle\Big\|_{L^2}\to 0,
\end{eqnarray*}
following from  $\E\{d R_k(h) |\F^{x_0}\}\to 0$.
\end{remark}

Observe that the pull back by $\I$ of $BC^2$ functions belong to $\D^{2,2}_{\F^{x_0}}(\Omega, \R)$ and  $\D^{2,2}_{\F^{x_0}}(\Omega, \R)\subset \I^*[\D^{2,1}(\C_{x_0}M;\R)]$ by  Proposition \ref{pr:chaos-expansion}. We are lead to the following :
\begin{corollary}
\label{co:weak-unique}
If Condition $(M_0)$ holds then
$BC^2$ functions on $\C_{x_0}M$ are in  $\D^{2,1}(\C_{x_0}M; \R)$.
\end{corollary}

Corollary \ref{co:weak-unique} corrects Theorem 2.1 of \cite{Elworthy-Li-CR-03} which stated that $BC^1$ functions are  in $\D^{2,1}(\C_{x_0}M; \R)$.

\subsection{On the uniqueness of $d$}

 \begin{definition} We say uniqueness holds for $d^p$ if the closure of $d^p$ is independent of the choice of its initial domain $d_\HH$ satisfying  $\Cyl\subset \Dom(d_\HH)\subset BC^1$.
 \end{definition}

\begin{remark}
Since $BC^1\subset W^{p,1}$ uniqueness for $d^p$ would be implied by Markov uniqueness. It would also follow more generally if the image $\I^*[\D^{p,1}]$ is independent of the choice of $\Dom(d_\HH)$.
\end{remark}

\begin{theorem}\label{th-strong}
Assume Condition $(M_0)$.
Suppose $T$ is a closed operator, 
$$T: Dom(T)\subset  L^2(\C_{x_0}M;\R) \subset L^2\Gamma \HH^* $$
with the properties that \begin{enumerate}
\item [(i)] $T$ agrees with $d_\HH$ on smooth cylindrical functions. 
\item [(ii)] $Dom(T^*)$ contains all smooth cylindrical one forms. 
\end{enumerate}
Then $d\subset T\subset \tilde d$,
where $\tilde d$ is the weak derivative.

 \end{theorem}

\begin{proof}
It is enough to show that $(\tilde d)^*\subset T^*$, which follows from
$${\tilde d}^*=\overline{d^*|_{\D^{2,1}}}=\overline{d^*|_{\Cyl}} \subset T^*.$$
\end{proof}

\begin{remark}
\label{Remark-7.9}
If $\E\{f|\F^{x_0}\}$ is in $\D^{2,1}(\Omega, \R)$ whenever $f$ is in 
$\D^{2,1}(\Omega, \R)$ then $\I^*[ \D^{2,1}]=\D^{2,1}_{\F^{x_0}}$ and so Markov uniqueness, and hence uniqueness of $d$, holds. To see this take $f$ in $\D^{2,1}_{\F^{x_0}}$ so $f=\lim_n f^n$ 
in $\D^{2,1}$ when $f^n$ is the sum of the finite terms in the chaos expansion of $f$. We saw above that each $\E\{f^n|\F^{x_0}\}$ lies in $\I^*[\D^{2,1}]$. Since the latter is closed $f$ itself must lie in it. See the Appendix, section \ref{Appen}.
\end{remark}

 \section{Covariant Differentiation}
 \label{section-covariant}
 
 Our main aim in this section is to define higher order Sobolev and weak Sobolev spaces and to prove the pull back theorem,  Theorem \ref{higher-order-pull}. 
 
 For $G$ a separable Hilbert space we can define $W^{p,1}G\equiv W^{p,1}(\C_{x_0}M;G)$ by
 Proposition \ref{closable_G} to be the domain, by graph norm, of the closure of $\tilde d^G$ whose domain $\Dom(\tilde d^G)$ consists of finite sums $\sum_{j} f_jg_j$ for
 $f_j\in W^{p,1}(\C_{x_0}M;\R)$, $g_j\in G$. Let 
 $$ \tilde d \equiv\tilde d^G : W^{p,1}G\subset L^p(\C_{x_0}M;G)\to L^p\Gamma(\LL_2(\HH;G))$$
 also denote this closure.
 \begin{proposition}
 \label{pr:weak-strong-sobolev}
 Assume Condition $(M_0)$.
 For $f: \C_{x_0}M\to \R$ the following are equivalent:
 \begin{enumerate}
 \item[(i)] $f\in W^{p,1}(\C_{x_0}M;G)$
 \item[(ii)] $\I^*f\in \D^{p,1}(\C_0\R^m;G)$;
 \item[(iii)] There is a constant $c_f$ such that if ${1\over p}+{1\over q}=1$ then for all $G$-valued 1-forms 
 $\phi\in \D^{q,1} \left(\C_{x_0}M; \LL_2(\HH;G) \right)$,
 \begin{equation}
 \left| \int_{\C_{x_0}M} \left\langle f, (d^G)^*\phi\right\rangle_G d\mu_{x_0}\right | \le
 c_f \|\phi\|_{L^q\Gamma\LL_2(\HH;G)}.
 \end{equation}
 If $f\in W^{p,1}(\C_{x_0}M;G)$ the intertwining formula of Theorem \ref{th:weak_Sob} and (\ref{weak-22}) extend to $G$-valued functions.
 \end{enumerate}
 \end{proposition}
 
 \begin{proof}
 Suppose $f\in W^{p,1}G$. Then $f=\lim_{l\to \infty} f_l$ in graph norm for some sequence $\{f_l\}_{l=1}^\infty$ in
 $G\otimes_0\Dom(\tilde d)$. By Theorem \ref{th:weak_Sob}, $\I^*(f_l)\in\D^{p,1}$ for each $l$ and
 $$d\I^*(f_l)=\I^*(\tilde d^Gf_l).$$
  By assumption $\tilde d^Gf_l\to \tilde d^Gf$ in $L^p$ and so by continuity of $\I^*$, Theorem \ref{th:pull-back-forms},
 $d\I^*(f_l)$ converges in $L^p$, showing that
 $\I^*(f)\in \D^{p,1}$ and $d(\I^*f)=\I^*(\tilde d^G f)$.
 
  Conversely, if $\I^*(f)\in \D^{p,1}(\C_0\R^m;G)$, taking an orthogonal base $\{g_i\}_{i=1}^\infty$ for $G$
  $$\I^*(f)=\lim_{k\to \infty}\sum_{l=1}^k \alpha_l g_l$$
  in $\D^{p,1}(\C_0\R^m;G)$ for $\alpha_l=\langle \I^*(f), g_l\rangle_G\in \D^{p,1}_{\F^{x_0}}(\C_0\R^m;\R)$. By 
  Theorem \ref{th:weak_Sob} $\alpha_l=\I^*(f_l)$ some $f_l\in W^{p,1}(\C_{x_0}M;\R)$. 
  Then $f=\lim_{k\to \infty} \sum_{l=1}^k f_lg_l$ in $L^p(\C_{x_0}M;\R)$. Since 
  $$\tilde d^G(\sum_{l=1}^k f_lg_l)(-)=\sum_{l=1}^k \tilde d f_l (-)g_l$$
   and $\tilde d f_l$ is given in terms of $d\alpha_l$ by equation (\ref{weak-22}), we see the 
   convergence is in $W^{p,1}(\C_{x_0}M; G)$
   and so $f\in W^{p,1}(\C_{x_0}M;G)$. 
   Thus (i) is equivalent to (ii). That (ii) is equivalent to (iii) can be seen as for the scalar case 
   in Theorem \ref{th:weak_Sob}. The only additional point is to observe that if $U$ is defined by (\ref{total_set_of_H})
   then $G\otimes_0U$ is total in $\D^{p,1}(\Omega, G\otimes H)$ and 
   $\overline{(\1\otimes T\I)(h)}\in \D^{p,1}(G\otimes \HH)$ for each $h\in G\otimes_0U$. 
   Thus $G\otimes_0U$ can take over the role $U$ played in the proof of Theorem \ref{th:weak_Sob}.
    This proof also shows that the analogue of (\ref{weak-22}) holds for $G$-valued functions.
   \end{proof}

  \subsection{The pointwise and the damped Markovan connections}
   To define higher order Sobolev spaces we need to introduce covariant derivatives.

  Using the notation of \S\ref{se:L2_tangent_spaces} consider the principal $\C_{\id}O(p)$-bundle
  $C_{x_0}OE\to \C_{x_0}M$ and its associated Hilbert bundle $L^2\Epsilon\to \C_{x_0}M$. As described in Eliasson
\cite{Eliasson67} our metric connection $\nabla$ on $E$ induces a so called pointwise connection $\bar\nabla$
  on $L^2\Epsilon$.
  \begin{equation}\label{pointwise-connection}
\bar\nabla_vU={D\over \partial s}U\left(\exp_{\sigma_\cdot} s v_\cdot\right)\Big|_{s=0},
\end{equation}
for $U$ a smooth section of $L^2\Epsilon$ and $v\in T_\sigma
C_{x_0}M$ where ${D\over \partial s}$ and $\exp$ come from $\nabla$.
  
  The almost surely defined map ${\D\over d\cdot}:\HH\to L^2\Epsilon$,  is an 
  isometric isomorphism which we used to give $\HH$ a vector bundle
   structure (at least over the subset on which
 ${D\over \partial \cdot}$ is defined). We also use it to pull back $\bar\nabla$ to obtain a
 connection on $\HH$, which we shall denote by ${\pnabla}$.
   By definition $$\pnabla=\W \bar \nabla { \D\over dt},$$
   see (\ref{wt}).
As usual these connections induce connections on the relevant
tensor bundles and in particular on the dual bundles
$(L^2\Epsilon)^*$ and ${\HH}^*$ respectively. Since the
connections are metric the latter are compatible with the natural
isometries
$$(L^2\Epsilon)^*\to L^2\Epsilon, \qquad\qquad
(\HH)^*\to \HH.$$

For  $X:M\times \R^m\to TM$  as
 in (\ref{sde}) define:
$$\bar X : C_{x_0}M\times L^2([0,T];\R^m) \to L^2\Epsilon $$
by
 $$\left(\bar X(\sigma)h\right)_t   =X(\sigma(t))(h(t)),
 \qquad 0\le t\le T$$
and its  right inverse:
$$\bar Y :  L^2\Epsilon \to L^2([0,T];\R^m)$$
$$ \bar Y(v)_t     =Y_{\sigma(t)}(v(t)),  \qquad 0\le t\le T.$$
Also define
$$ \X : C_{x_0}M\times H \to \HH$$ by
$$ \X(\sigma)(h)=\overline{T\I}_\sigma(h)\equiv
\W \bar X(\sigma)(\dot h_\cdot),
 \qquad 0\le t\le T$$
with right inverse $\Y:\HH\to H$, as defined by  (\ref{Y}).

\begin{proposition}
\label{pr:connection-3}
The connections $\bar\nabla$, $\pnabla$ on $L^2\Epsilon$ and $\HH$ are the 
connections corresponding to $\bar X$ and $\X$ respectively in the sense that 
$\bar\nabla_v \bar X(\dot h)=0$ and
$\pnabla_v \X(h)=0$ if $\dot h\in [\ker \bar X(\sigma)]^\perp$ or equivalently 
$h\in [\ker \X(\sigma)]^\perp$ 
and $v\in T_\sigma \C_{x_0}M$.

If $U$, ${\D\over \partial t}V$ are $C^1$ sections of $L^2\Epsilon$ and $v\in T_\sigma\C_{x_0}M$ then
\begin{equation}
\label{pointwise-connection-2}
\bar\nabla_vU=\bar X(\sigma) d\left(\alpha\mapsto \bar Y_\alpha U(\alpha)\right)(v)
\end{equation}
and
\begin{equation}
\label{Markovian-connection}
\pnabla_vV=\X(\sigma) d\left(\alpha\mapsto \Y_\alpha V(\alpha)\right)(v).
\end{equation}
\end{proposition}

\begin{proof}
By definition 
\begin{eqnarray*}
\bar \nabla _v \bar X(\dot h)_t&=&
{D\over ds} X\left(\exp_{\sigma(t)} sv_t\right)(\dot h(t))|_{s=0}\\
&=&\nabla_{v_t} X(\dot h(t))\\
&=&0, \hbox{ if } \dot h(t)\perp \ker[X(\sigma(t))].
\end{eqnarray*}
However $\dot h\perp \ker[\bar X(\sigma)]$ in $L^2([0,T];\R^m)$ holds if and only if $\dot h(t)$ is 
orthogonal to $\ker[X(\sigma(t))]$ for almost all $t$, and we see the pointwise connection
 is the connection corresponding to $\bar X$. From this (\ref{pointwise-connection-2}) holds
 as in Elworthy-LeJan-Li \cite{Elworthy-LeJan-Li-book}. The results for $\HH$ follow by conjugation with ${\D\over\partial t}$.
\end{proof}

In the case $E=TM$ the connection $\pnabla$ is the {\it damped Markovian connection} defined
 in a different way by Cruzeiro-Fang \cite{Cruzeiro-Fang95} and we refer to it as that in our more general situation. For the case of $M$ a Lie group with left invariant metric and connection, as Example 2 of section \ref{sec:Examples}, see section \ref{se:special-case}.
  
  \subsection{Covariant Gross-Sobolev derivatives}
  \label{sec:covariant}
  
 Let  $G$ be a separable Hilbert space.  First consider a smooth separable Hilbert bundle $\G$ over $\C_{x_0}M$ with a metric connection 
 $\tilde\nabla$ determined by a smooth surjective vector bundle map
   $\tilde X: \C_{x_0}M\times G\to \G$, with isometric right inverse
   $$\tilde Y_\sigma=\tilde X(\sigma)^*: \G_\sigma \to G, \qquad  \sigma\in\C_{x_0}M.$$
   Define $\Dom(\tilde\nabla^p)$ to be $\D^{p,1}\G$ for
   $$\D^{p,1}\G=\left\{ V\in L^p\Gamma\G \; |\;
   \tilde Y V\in \D^{p,1}(\C_{x_0}M;G)\right\}.$$
  For $V\in \D^{p,1}\G$, $\phi\in \D^{p,1}\LL_2(\G;K)$, $K$ a separable Hilbert space set
  \begin{eqnarray}
  \label{covariant-21}
  \tilde\nabla_v^pV&=&\tilde X(\sigma) d^p\left(\tilde YV\right)_\sigma(v), \qquad v\in \HH_\sigma.\\
  \tilde \nabla^p_v \phi&=&d^p(\phi (\tilde X-))_\sigma (v)(\tilde Y_\sigma-).
  \label{covariant-22}
  \end{eqnarray}
  
   Let $\tilde K^\perp(\sigma)=\tilde Y_\sigma \tilde X(\sigma)$ be the orthogonal projection
   of $G$ onto $\ker[\tilde X(\sigma)]^\perp$, $\sigma\in \C_{x_0}M$. 
    Note that, c.f. \cite{Elworthy-LeJan-Li-book}:
    \begin{equation}
    \label{nabla-K}
   \tilde X \tilde \nabla \tilde Y=0, \quad \hbox{ and } \quad \tilde Xd_\HH \tilde K^\perp=0.
    \end{equation}

 To obtain a closed covariant differentiation operator we impose:
  \medskip
  
\noindent   
  {\bf Condition $K^p$.} If $f\in \D^{p,1}(\C_{x_0}M;G)$ then $\tilde K^\perp(\cdot) f(\cdot)\in \D^{p,1}(\C_{x_0}M;G)$. 
  \medskip
  
  Note that this implies , by the closed graph theorem, that $f\mapsto \tilde K^\perp f$ is continuous from $\D^{p,1}(\C_{x_0}M;G)$ to $\D^{p,1}(\C_{x_0}M;G)$.   Note that if $\tilde X=\bar X$ or $\X$,  then Condition $(M_0)$ implies 
  Condition $K^p$ holds for all $p$ by Lemma \ref{le:multiplication}.
   
   \begin{proposition}
   \label{dense-DG}
   Assume Condition $K^p$. Then $C^\infty\Gamma\G\cap \D^{p,1}\G$ is dense in $\D^{p,1}\G$.
   \end{proposition}
   \begin{proof}
   For $V\in \D^{p,1}\G$ set $f=\tilde Y V$ and take smooth functions $f_k\in \D^{p,1}(\C_{x_0}M;G)$, $k=1,2,\dots$, converging to $f$ in $\D^{p,1}(\C_{x_0}M;G)$. Then $\tilde X f_k\to V$ in $L^p$. Observe that 
   $$\tilde\nabla(\tilde X f_k)=\tilde X d_\HH(\tilde Y\tilde X f_k)\to \tilde X d(\tilde Y\tilde X f)=\tilde \nabla V$$
   in $L^p$ by Condition $K^p$.
   \end{proof}
   
   The proposition is essentially due to the fact that Condition $K^p$ implies the smoothness of $\tilde X$ and in particular 
   \begin{equation}
 U \in\D^{p,1}(\C_{x_0}M;G)\Longleftrightarrow
 \tilde X(U) \in \D^{p,1}\G  .
\end{equation}

   \begin{corollary}
   Assume Condition $K^p$ . There is the Leibniz formula
   \begin{equation}
\label{Leibniz}
d^p(\tilde Y V)=(\tilde\nabla_-\tilde Y)V+\tilde Y\tilde\nabla_-^p V, \qquad V\in \D^{p,1}\G.
\end{equation}
   \end{corollary}
   \begin{proof}
   The formula holds for smooth $V$ and so in general by Proposition \ref{dense-DG} since $d^p$ is closed.
   \end{proof}
  \begin{lemma}
 Assume Condition $K^p$. Then for $1<p<\infty$, and for $1\le p<\infty$ if condition $(M_0)$ holds,
$$\tilde{\nabla}^p: \Dom(\tilde\nabla^p)
\subset L^p\Gamma\G\to L^p\Gamma(\LL_2(\HH;\G))$$
is a closed operator. 
\end{lemma}
\begin{proof}
Let $\{ U^k\} _{k=0}^\infty $ be a sequence in $\D^{p,1}\G$ such that 
 $\{U^k,\tilde{\nabla}^pU^k\} _{k=0}^\infty $ converges to some $(U,Z)$. 
 We must show that $U\in \D^{p,1}\G $
and $Z=\tilde{\nabla}^pU$. By definition $\tilde Y(U^k)\in \D^{p,1}(\C_{x_0}M; G)$ and
$\tilde Y(U^k)\to \tilde Y(U)$.   Apply the Leibniz formula to see that
  $$d(\tilde Y U^k)= (d_\HH \tilde K^\perp)(\tilde Y U^k)+ \tilde Y \tilde\nabla^p U^k
  \to (d_\HH \tilde K^\perp)(\tilde Y U)+\tilde Y Z.$$ 
 The convergence is in $L^p\Gamma\LL_2(\HH;G)$ by Condition $K^p$. Since $d$ is closed this shows that
 $\tilde Y U\in \D^{p,1}(\C_{x_0}M; G)$ and $\lim_{k\to \infty}d(\tilde Y U^k)=d(\tilde Y U)$.
 Consequently $U\in \D^{p,1}\G$ and, using (\ref{covariant-21}), $\tilde\nabla^p U^k= \tilde X d(\tilde Y U^k)\to \tilde X d(\tilde Y U)=\tilde\nabla^p U$, giving $Z=\tilde\nabla^p U$.
 \end{proof}

  Note that if $\G^1$ and $\G^2$ are smooth Hilbert bundles over $\C_{x_0}M$ with metric 
  connections given respectively by $$\tilde X^j: \C_{x_0}M\times G^j\to \G^j, j=1,2$$
  then the natural induced metric connections on the Hilbert bundles $(\G^1)^*$, 
  $\G^1\otimes \G^2$, $\LL_2(\G^1;\G^2)$ are determined by
  \begin{eqnarray*}\C_{x_0}M\times (G^1)^*&\to& (\G^1)^*\\
   (\sigma, l)&\mapsto& l\circ Y_\sigma,
   \end{eqnarray*}
   $$\tilde X^1\otimes \tilde X^2: \C_{x_0}M\times (G^1\otimes G^2)\to \G^1\otimes \G^2$$
  and
  \begin{eqnarray*}\C_{x_0}M\times\LL_2(G^1;G^2)&\to& \LL_2(\G^1;\G^2)\\
  (\sigma, T)&\mapsto& \tilde X^2(\sigma)T\tilde Y^1_\sigma
   \end{eqnarray*}
  respectively. Using these, iteratively, we can obtain closed operators acting on sections of the
  tensor bundles constructed from an initial Hilbert bundle $\G$, given the relevant Condition $K^p$.
  For `bundles' of the form $\LL_2(\HH; \LL_2(\HH;\G))$ which we can more compactly
  write as $\G\otimes (\HH^{\otimes^2})^*$ we will use the isometry
  ${\D\over\partial t}: \HH\to L^2\Epsilon$ to pull back the covariant derivative operator
  from the one obtained as above on $\G\otimes \left((L^2\Epsilon)^{\otimes^2}\right)^*$,
  (in this case) using the pointwise connection on $L^2\Epsilon$. For example
  $\D^{p,1}\left(\G\otimes (\HH^{\otimes^2})\right)$ is
  $$\left\{ V\in L^p \left(\G\otimes (\HH^{\otimes^2})\right)\, |\,
  \left(\1\otimes{\D\over\partial t}\otimes{\D\over \partial t}\right)V
   \in \D^{p,1}\left(\G\otimes(L^2\Epsilon)^{\otimes 2}\right)\right\}.$$
  
  \begin{theorem}
  \label{th:the-mess}
  Assume Condition $K^p$ for $\bar X$ for all $1<p<\infty$.
  Then for $1<p<\infty$ the above construction yields closed covariant 
  derivative operators with domain $\Dom(\pnabla^p)$:
  \begin{eqnarray*}
  \pnabla^p:&& L^p\Gamma \left( \left(L^2\Epsilon\right)^{\otimes r}\otimes
   \left(L^2\Epsilon\right)^{\otimes s}\otimes \HH^{\otimes a}\otimes \left(\HH^*\right)^{\otimes b}\right)\\
 &&  \longrightarrow L^p\Gamma\left( \left(L^2\Epsilon\right)^{\otimes r}\otimes \left(L^2\Epsilon \right)^{\otimes s}\otimes \HH^{\otimes a}
 \otimes \left(\HH^*\right)^{\otimes b+1}\right)
 \end{eqnarray*}
  for any $r,s,a,b,\in \{0,1,2,\dots, \}$. Moreover $V\in \Dom(\pnabla^p)$
  if and only if 
\begin{eqnarray*}  X^{r,s,a,b} V &\equiv&
\left(  \left(\otimes^{r}\bar X\right)\otimes 
\left( \circ (\otimes^{s}  \bar Y)\right) \otimes
 \left( \otimes^{a}\X\right) \otimes \left( \circ (\otimes^{b} \Y)\right)
 \right)V:\\
  \C_{x_0}M&\longrightarrow& \left(L^2([0,T];\R^m\right)^{\otimes^r} \otimes 
  \left(L^2([0,T];\R^m\right)^{\otimes^s}
  \otimes H^{\otimes ^a}\otimes(H^*)^{\otimes^b}
  \end{eqnarray*}
  is in $\D^{p,1}$ and then $\pnabla^p=(X^{r,s,a,b})^*d^p(X^{r,s,a,b}V)$.
  If Condition $(M_0)$ holds we may take $p=1$.
  \end{theorem}
  \begin{proof}
  It is only necessary to observe that Condition $K^p$ for $X^{r,s,a,b} $ is implied by Condition $K^q$ for all $1<q<\infty$ for $\bar X$, and if Condition $(M_0)$ holds then $X^{r,s,a,b} $ satisfies $K^1$ by Lemma
   \ref{le:multiplication}.
  \end{proof}

\subsection{The higher order Sobolev spaces $\D^{p,k}$}
\label{se:High-Sobolev}

Suppose that $\G$ is a smooth Hilbert bundle over $\C_{x_0}M$ with connection
 $\tilde\nabla$ as given in \S\ref{sec:covariant} by some $\tilde X$ which together with all tensor products  
  $$\tilde X\otimes (\otimes^a \bar X)\otimes (\circ \otimes^b \bar Y)$$
satisfies Condition $K^p$.
 For $a,b \in\{0,1,2,\dots\}$ we can inductively define $\pnabla^{p, (k)}$ and $\D^{p,k}\left(\G\otimes\HH^{\otimes a} \otimes (\HH^*)^{\otimes b}\right)$, $k=1,2,\dots$ as follows: \\
 Set $\pnabla^{p,(1)}=\pnabla^p$, defined as in Theorem \ref{th:the-mess}, with
 $$\D^{p,1}\left(\G\otimes\HH^{\otimes a} \otimes (\HH^*)^{\otimes b}\right)=\Dom(\pnabla^p).$$
 For $k=\{2,3,\dots\}$ set
$$\D^{p,k}\left(\G\otimes\HH^{\otimes a} \otimes (\HH^*)^{\otimes b}\right)
 =\left\{V\in \D^{p,1} \; |\; \pnabla^{p}V \in \D^{p, k-1}
  \left(\G\otimes \HH^{\otimes a}\otimes (\HH^*)^{\otimes (b+1)}\right)\right\}$$
  and $\pnabla^{p,(k)}=\pnabla^{p, (k-1)}\circ \pnabla^p$. Here we have used
  our usual identification of $\LL_2(H_1;H_2)$ with $H_2\otimes H_1^*$ for
  Hilbert spaces $H_1, H_2$. As usual we give  $\D^{p,k}$ the graph norm
  $$\n{V}_{p,k}
  =\left(\n{V}_{L^p}^p+\n{\pnabla^p V}^p_{L^p}+\dots+
  \n{\pnabla^{p, (k)}V}_{L^p}^p
  \right)^{1\over p},$$
  \ie the graph of the closed operator
  $\id\otimes \pnabla^p\oplus\dots\oplus \pnabla^{p, (k)}$.

\subsubsection{The higher order weak Sobolev spaces $W^{p,k}$}

Continuing with the previous notation let
\begin{eqnarray*}
(\pnabla^p)^*: \Dom(\pnabla^p)^*
&\subset& L^q\Gamma\left(\G\otimes \HH^{\otimes r}\otimes {\HH^*}^{\otimes s+1}\right)\\
&\longrightarrow& L^q\Gamma\left(\G\otimes \HH^{\otimes r}\otimes {\HH^*}^{\otimes s}\right)
\end{eqnarray*}
be the adjoint of $\pnabla^p$, ${1\over p}+{1\over q}=1$, with $q=\infty$ if
$p=1$.

\begin{lemma}
\label{le:higher-order}
Assume Condition $K^p$ holds for $\tilde X^{r,s}$ and  $\tilde\X^{r,s}$ where
\begin{eqnarray*}
\tilde X^{r,s}& = &\tilde X\otimes \left( \otimes^r \bar X\right)\otimes \left(\circ(\otimes^s  \bar Y)\right), \\
\tilde\X^{r,s}&=&\tilde X\otimes \left(\otimes^r  \X\right)\otimes \left(\circ(\otimes^s \Y\right).
\end{eqnarray*} Then
\begin{enumerate}
\item[(i)]
$(\pnabla^p)^*=\tilde \X^{r,s}(d^p)^*\left(\tilde \X^{r,s+1}\right)^*$.
\item[(ii)]
$\D^{q,1}\left(\G\otimes \HH^{\otimes r}\otimes {\HH^*}^{\otimes s+1}\right)
\subset \Dom(\pnabla^p)^*$.
\end{enumerate}
\end{lemma}

\begin{proof}
After conjugation with ${\D\over dt}$ if necessary we can assume that $r=s=0$. Then by (\ref{covariant-21}), $\tilde \nabla^p (-)=\tilde X d^p \tilde Y(-)$ and it is easy to see that
$$(\tilde\nabla^p)^*\supset \tilde X (d^p)^* (\tilde Y\otimes 1)=\tilde X(d^p)^* (\tilde Y\circ -).$$
Suppose that 
$\phi\in \Dom(\tilde \nabla^p)^* \subset L^q\Gamma(\G\otimes \HH^*)\sim L^q\Gamma \LL_2(\HH; \G)$.
For (i) it suffices to show that $\tilde Y\circ \phi\in \Dom((d^p)^*)$. For this take
 $g\in \D^{p,1}(\C_{x_0}M; G)$ then
 \begin{eqnarray*}
 &&\left|\int_{\C_{x_0}M} \left\langle \tilde Y\circ \phi, d^p g\right\rangle d\mu_{x_0}\right|\\
 &=& \left|\int_{\C_{x_0}M} \left\langle \phi, \tilde X d^p (\tilde K^\perp g)\right\rangle d\mu_{x_0}\right|
 + \left|\int_{\C_{x_0}M} \left\langle  \phi, \tilde X d^p(\tilde K g)\right\rangle d\mu_{x_0}\right|.
 \end{eqnarray*}
 Now $\tilde K^\perp g=\tilde Y (\tilde X g)$ and $\tilde X g\in \D^{p,1}\G$. So the
 first of these two terms is
 $$ \left|\int_{\C_{x_0}M} \left\langle  \phi, \tilde \nabla^p (\tilde X g)\right\rangle d\mu_{x_0}\right|
 \le c_\phi|\tilde X g |_{L^p} \le \const \cdot c_\phi\, |g |_{L^p}.$$
The second term is
 \begin{eqnarray*}
 \left|\int_{\C_{x_0}M} \left\langle  \tilde Y \phi, d^p(\tilde K g)\right\rangle d\mu_{x_0}\right|
 &=&\left|\int_{\C_{x_0}M} \left\langle  \tilde Y \phi, d^p(\tilde K)(-) g+\tilde K d^p g\right\rangle d\mu_{x_0}\right|\\
&\le&|\const \phi|_{L^q} |g|_{L^p}.
 \end{eqnarray*}
For (ii) suppose $\phi\in \D^{q,1} \LL_2(\HH; \G)$. By definition, 
$$\tilde Y\circ \phi \in \D^{q,1} \left(\LL_2(\HH;G)\right)\subset \Dom(d^p)^*$$
just as in the scalar case,  Theorem  \ref{th:divergence-21}. However (the easy part of ) (i) then shows $\phi\in \Dom(\tilde\nabla^p)^*$ as required.
\end{proof}

Define $W^{p,1}\left(\G\otimes \HH^{\otimes r}\otimes {\HH^*}^{\otimes s}\right)$ to be domain of the adjoint of the restriction of $(\pnabla^p)^*$
to $\D^{q,1}\left(\G\otimes \HH^{\otimes r}\otimes {\HH^*}^{\otimes s+1}\right)$.
Let \begin{eqnarray*}
\tilde \pnabla^p : W^{p,1}\left(\G\otimes \HH^{\otimes r}\otimes {\HH^*}^{\otimes s}\right)
&\subset& L^p\Gamma\left(\G\otimes \HH^{\otimes r}\otimes {\HH^*}^{\otimes s}\right)\\
&\longrightarrow& L^p\Gamma\left(\G\otimes \HH^{\otimes r}\otimes {\HH^*}^{\otimes s+1}\right)
\end{eqnarray*}
be the adjoint , considered as a closed operator.  

As for $\D^{p,k}$ we can define $W^{p,k}$ and $\tilde \pnabla^{p,(k)}$ iteratively,
 giving $W^{p,k}$ the analogous graph norm. Since 
 $\pnabla^{p,(k)}\subset \tilde \pnabla^{p, (k)}$ we see $\D^{p,k}$ is always a closed
 subspace of $W^{p,k}$. When $\G$ is a trivial vector bundle $C_{x_0}M\times  G$ we write
them as  $W^{p,k}(\C_{x_0}M;G)$ and $\D^{p,k}(\C_{x_0}M;G)$ respectively. By Proposition
 \ref{pr:weak-strong-sobolev} this agrees with the previous definition when $k=1$.
 
We will consider the following possible conditions on $\tilde X$:

\medskip

\noindent {\bf Condition K(N).}
  For $1\le k\le N$, for $g\in G$, the $k^{th}$ $\HH$-derivative of 
   $\tilde K^\perp(-)(g):\C_{x_0}M\to G$ has a bound
   $$
   \left| (\pnabla^{\otimes (k-1)} d_\HH \tilde K^\perp(-)(g))
   \right|_{G\otimes (\HH^*)^{\otimes (k-1)}}<c_N(\sigma)|g|_G,$$
   where $\bar c_N=\esssup_{\sigma\in\C_{x_0}M}  c_N(\sigma)$ is finite.

\noindent {\bf Condition K(N)$\D$.}
For each $g\in G$ the map 
$\tilde K^\perp(-)(g)$ is in $\cap_{p>1}\cap_{k=1}^N \D^{p,k}$
and Condition $K(N)$ holds.
\begin{remark}
\label{re:condition}
If Condition $M$ holds then Condition $K(N)\D$ holds for $\bar X$ and $\X$. 
To see this observe that Condition $K(N)$ holds because Condition $(M)$ implies that the norm $\|i_\sigma\|_\sigma$ of the inclusion $i_\sigma: \HH_\sigma\to T_\sigma \C_{x_0}M$ is in $L^\infty(\C_{x_0}M;\R)$ while all the derivatives of $K^\perp: M\to \LL(\R^m;\R^m)$ are bounded. Then Condition $K(N)D$ follows as in the proof of Lemma \ref{le:multiplication}. Condition $(M_0)$ suffices for Condition $K(1)D$.
\end{remark}

\begin{remark}
\label{Re-8.9}
The necessity of imposing Condition $K(N)D$ in order to discuss $\D^{p,k}$ sections of $\G$ reflects the fact that we have not shown that Fr\'echet $BC^\infty$ functions are in $\cap_{p>1}\cap_{k=1}^\infty \D^{p,k}$, and $\G$ should be ``of class $\D^{p,k}$ '' in some sense.
\end{remark}

\begin{lemma} 
\label{le:K-continuity}
Assume Condition $(M_0)$.
Let $1\le p<\infty$. 
\begin{enumerate}
\item[(i)]
Under Condition K(N)$\D$ the map 
$$f\mapsto \tilde K^\perp f$$
gives a continuous linear map 
$$\D^{p,k}(\C_{x_0}M;G)\to \D^{p,k}(\C_{x_0}M;G) $$
for $k\in\{1,\dots, N\}$ and $f\mapsto d_\HH\tilde K^\perp(-)f$ 
a continuous linear map
$$\D^{p,k-1}(\C_{x_0}M;G)\to \D^{p,k-1}(\C_{x_0}M;\LL_2(\HH;G))$$
for $k\in\{1,2,\dots, N\}$.
\item[(ii)] Under Condition K(N) the corresponding  results hold for weak derivatives.
\end{enumerate}
\end{lemma}

\begin{proof}
Assume Condition $K(N)$. Consider $\I^*(\tilde K^\perp(\cdot)f(\cdot))$ to see 
$$\tilde K^\perp(\cdot)f(\cdot)
\in W^{p,1} (\C_{x_0}M;G)$$ and 
$$\tilde d^p(\tilde K^\perp f)=d_\HH(\tilde K^\perp)(-)f+\tilde K^\perp(\cdot)\tilde d^p f$$
by Proposition \ref{pr:weak-strong-sobolev}. Repeat this for $d_\HH(\tilde K^\perp)(\tilde X-)f$ and higher derivatives to prove (ii) (The continuity comes from the closed graph theorem). For (i) 
assume K(N)$\D$. We already have the result for the weak derivatives. If $f\in \Dom(d_\HH)^G$,
\ie $f\in \Cyl\otimes G$ we see $\tilde K^\perp f\in \D^{p,1}$, as do successive
derivatives of $\tilde K^\perp(-)f$. Since $\D^{p,1}$ is 
closed in $W^{p,1}$ we obtain the result for $k=1$, and by iterating this the result for $1\le k\le N$ as required.

\end{proof}

\begin{proposition}
\label{pr:higher-order-connection}
Suppose $1<p<\infty$ and $k\in\{1,2,\dots\}$. Assume Condition $(M_0)$
and that $\tilde \X^{r,s}$ defined in Lemma \ref{le:higher-order} satisfies Condition $K(k)$. 
Then
$$ V\in W^{p,k} \left(\G\otimes \HH^{\otimes r}\otimes {\HH^*}^{\otimes s}\right)$$  
if and only
if  $(\tilde \X^{r,s})^*V\in \W^{p,k}\left(\C_{x_0}M; 
G\otimes H^{\otimes r}\otimes {H^*}^{\otimes s}\right)$.
Furthermore  $\tilde \pnabla^{p,(k)}= \tilde \X^{r,s}{ \tilde d}^{(k)} (\tilde \X^{r,s})^*V$.
The corresponding result holds for $\D^{p,k}$ assuming Condition K(k)$\D$.
\end{proposition}
\begin{proof}
After conjugation by ${\D\over dt}$ we can assume that $r=s=0$.
Suppose $V\in W^{p,1}\G$.  Then there exists $c_V$ such that
\begin{equation}
\label{covariant-adjoint}
\left |\int_{\C_{x_0}M} \left\langle V, (\tilde\nabla^q)^*Z\right\rangle_\sigma d\mu_{x_0}\right|
\le c_V|Z|_{L^q}, \qquad \forall \;Z\in \D^{q,1}\left(\G\otimes \HH^*\right). 
\end{equation}
On the other hand suppose 
$\phi\in \D^{q,1}(G\otimes \HH^*)\sim \D^{q,1}\LL_2(\HH;G)$.
Then, by definition,
$\tilde X\circ \phi\in \D^{q,1}\LL(\HH;\G) \sim \D^{q,1}(\G\otimes \HH^*)$ and so
\begin{eqnarray*}
&&\left |\int_{\C_{x_0}M} \left\langle \tilde Y V, (d^q)^* \phi \right\rangle_\sigma d\mu_{x_0}(\sigma)\right|
=\left |\int_{\C_{x_0}M} \left\langle  V, \tilde X (d^q)^* \phi \right\rangle_\sigma d\mu_{x_0}(\sigma)\right|\\
&\le&\left |\int_{\C_{x_0}M} \left\langle V, \tilde X(d^q)^* \tilde K^\perp \phi \right\rangle_\sigma d\mu_{x_0}(\sigma)\right|
+\left |\int_{\C_{x_0}M} \left\langle V, \tilde X(d^q)^* \tilde K \phi \right\rangle_\sigma d\mu_{x_0}(\sigma) \right|.
\end{eqnarray*}

Take $Z=\tilde X\circ \phi$ and use Lemma \ref{le:K-continuity} with equation (\ref{covariant-adjoint}) to bound the first of these two 
terms by a constant times $|\phi|_{L^q}$. 
To obtain a similar bound for the second term
observe that there is a constant such that
$$\left |\int_{\C_{x_0}M} \left\langle U, \tilde X(d^G)^* \tilde K \phi \right\rangle d\mu_{x_0}\right|
\le c|\phi|_{L^q}, \qquad \forall \; U\in L^p \Gamma(\G),$$
 because if $U\in\D^{p,1}$, a dense subset of $L^p$, we have
$\tilde Y U\in \D^{p,1}$ and
\begin{eqnarray*}
\left |\int_{\C_{x_0}M} \left\langle d^q(\tilde Y U), \tilde K\phi \right\rangle d\mu_{x_0}\right|
&=&\left |\int_{\C_{x_0}M} \left\langle (\tilde \nabla_ - \tilde Y) U
+\tilde Y\tilde \nabla_ -U, \tilde K\phi \right\rangle d\mu_{x_0}\right|\\
&\le&\bar c_1 |U|_{L^p} |\phi|_{L^q}
\end{eqnarray*}
by Condition $K(k)$, and since $\tilde K\tilde Y=0$.  This proves that $V\in W^{p,1}\G$ implies $\tilde Y V\in W^{p,1}(\C_{x_0}M; G)$.

Conversely suppose $\tilde Y V\in W^{p,1}(\C_{x_0}M;G)$ and $Z\in \D^{q,1}(\LL_2(\HH; \G))$.
We shall show (\ref{covariant-adjoint}) holds. Observe $\tilde Y\circ Z\in \D^{q,1}(\LL_2(\HH;G))$. By Lemma \ref{le:higher-order},
 \begin{eqnarray*}
 \left |\int_{\C_{x_0}M} \left\langle V, (\tilde \nabla^p)^*Z \right\rangle_\sigma d\mu_{x_0}(\sigma)\right|
 &=& \left |\int_{\C_{x_0}M} \left\langle V, \tilde X(d^p)^*(\tilde Y\circ Z) \right\rangle_\sigma d\mu_{x_0}(\sigma)\right|\\
 &=& \left |\int_{\C_{x_0}M} \left\langle \tilde Y V, (d^p)^*(\tilde Y \circ Z) \right\rangle_\sigma d\mu_{x_0}(\sigma)\right|\\
& \le& C_{\tilde Y V} |\tilde Y\circ Z|_{L^q} \le C_{\tilde Y V} |\tilde  Z|_{L^q}.
 \end{eqnarray*}
 So $V\in W^{p,1}\G$ and the result holds for $k=1$.
 
 Suppose now that the proposition holds for $W^{p,k-1}$ some $k\in\{2, 3,\dots\}$. Take $V\in W^{p,k}\G$.
  Then $V\in W^{p,k-1}\G$ and $\tilde \pnabla^p V\in W^{p,k-1}\G$. Equivalently
  $\tilde Y V\in W^{p,1}(\C_{x_0}M; G)$ and 
  $\tilde Y\tilde \pnabla^p V\in W^{p,k-1}\LL_2(\HH;G)$. Now 
  $$\tilde d(\tilde YV)=\tilde d(\tilde K^\perp \tilde Y V)
  =(d_\HH \tilde K^\perp)(-)\tilde Y V+\tilde K^\perp \tilde d(\tilde Y V)
  =(d_\HH \tilde K^\perp)(-)\tilde Y V+\tilde Y \tilde \pnabla^p V,$$
  which belongs to $W^{p,k-1}(\C_{x_0}M; G)$ by the previous lemma. 
   Thus $V\in W^{p,k}\G$  implies $\tilde Y V\in W^{p,k}(\C_{x_0}M;G)$. 
    
    Conversely if $\tilde Y V\in W^{p,k}$ then 
    $$\tilde Y\tilde \pnabla^p V=\tilde d(\tilde Y V)-\tilde \pnabla^p_-\tilde Y(V)
    =\tilde d(\tilde YV)-(d_\HH \tilde K^\perp)(-)\tilde Y V$$
    is in $W^{p, k-1}$ by the previous lemma, and so by the induction hypothesis $V\in W^{p,k}\G$.
\end{proof}

Note that we have not discussed the independence of $\D^{p,k}\G$ from the choice of $\tilde X$. However our main interest is in sections of $L^2\Epsilon$ or $\HH$ and related tensor bundles. For these we used $\tilde X$ derived from a stochastic differential equation and for such $\tilde X$ this is clear, as in Lemma \ref{le:multiplication}. An extension of Proposition \ref{dense-DG} to $\D^{p,k}$ and of the discussion in section \ref{d-uniqueness} to prove that $BC^
\infty$ functions are in $\cap_{p>1}\cap_{k=1}^\infty \D^{p,k}$, \cf Remark \ref{Re-8.9}, would give the general case.

As for the proof of Theorem \ref{th:divergence-21}, we have, using Remark \ref{re:condition}
\begin{theorem}
\label{th:continuity-2}
Assume Condition $(M_0)$.
  For $1<p<\infty$, the set $W^{p,1}\HH$ is contained in $\Dom({\div}^p)$ and 
  ${\div}^p:W^{p,1}\HH\to L^p(\C_{x_0}M;\R)$ is continuous.
  \end{theorem}
 
\subsection{Intertwining of higher order derivatives}

First we define the Sobolev spaces over the Wiener space relative to the It\^{o} map.
\begin{definition}
For any sub $\sigma $-algebra ${\mathfrak a}$ of the Borel $\sigma $-algebra of $\R ^m$
and separable Hilbert space $G$, the space $\D _{\frak{a},G}^{p,1}$
consists of those $F\in L^p\left( C_0\R ^m;G\right) $ s.t.
\begin{enumerate}
\item  $F$ is $\frak{a}$-measurable,
\item   $F\in {\Dom}_G(d^p)$.
\end{enumerate}
Inductively $\D _{\frak{a},G}^{p,k}$ consists of those $F\in L^p\left( C_0\R^m;G\right) $ such that
\begin{enumerate}
\item  $F\in \D _{\frak{a},G}^{p,1}$
\item  $\E \left\{ d^pF|\frak{a}\right\} : C_0R^m\rightarrow L_2\left(
H;G\right) $ is in $\D _{\frak{a},\LL_2(H; G)}^{p,k-1}$.
\end{enumerate}
For $F\in \D _{\frak{a},G}^{p,k}$ define
$$
\left( \left\| F\right\| _{\D _{\frak{a},G}^{p,k}}\right) =\left(
\sum_{r=0}^k\left\| \left( d^p\E \left\{ -|\frak{a}\right\}\right)^r F
\right\| _{L^p}^p\right) ^{1\over p}.
$$
\end{definition}
These spaces are Banach spaces since $d^p\E \left\{ -|\frak{a}\right\} F$
is the composition of a closed operator following a bounded operator and thus
a closed operator. 
\medskip

 The following theorem corrects the version in \cite{Elworthy-Li-CR-03}.
\begin{theorem}
\label{higher-order-pull}
Assume Condition $(M)$. Let $G$ be a separable Hilbert space. Then $f:\C_{x_0}M\to G$ is in
$W^{p,k}(\C_{x_0}M; G)$ some $1\le p<\infty$, $k\in\{1,2,\dots\}$ if and only if 
$\I^*(f)\in\D^{p,k}_{\F,G}$. Moreover $$\I^*:W^{p,k}(\C_{x_0}M;G)\to \D_{\F,G}^{p,k}$$
is a continuous linear isomorphism.
\end{theorem}

\begin{proof}
For $k=1$ this is Proposition \ref{pr:weak-strong-sobolev}.  Suppose it holds for 
some $k\in \{1,2,\dots\}$. Then if $f\in W^{p,1}(\C_{x_0}M;G)$ 
\begin{eqnarray*}
\E\{d\I^*(f)|\F^{x_0}\}&=&\E\{\I^*(\tilde df)|\F^{x_0}\}
=\I^*(\tilde d f \circ \overline{T\I})=\I^*(\tilde d f\circ \X),
\end{eqnarray*}
which belongs to $\D^{p,k}$ if and only if $\tilde d f\circ \X\in W^{p,k}(\C_{x_0}M;H^*)$.
The last holds if and only if $f\in W^{p, k+1}(\C_{x_0}M;G)$ by Proposition \ref{pr:higher-order-connection}
and so the result holds for $k+1$.
\end{proof}

\begin{corollary}
\label{co:intertwining2}
Under the conditions of Theorem \ref{higher-order-pull},
 $\I^*$ maps $\D^{p,k}(\C_{x_0}M;G)$ onto a closed subspace of $\D^{p,k}_{\F,G}$.
\end{corollary}

Note that many authors, \eg L\'eandre \cite{Leandre93}, Li \cite{X-D-Li-2003},  have defined Sobolev spaces for $\C_{x_0}M$ using the flat connection on $\HH$ defined by the trivialisation of $\HH$ given by $V_\cdot\mapsto \parals_\cdot^{-1} V_\cdot$ (in the Levi-Civita case). This cannot be expected to agree with our definition because of the curvature term occurring in the covariant derivative of $\paral^{-1}$, \cf the proof of Lemma \ref{le:density2}. In particular Corollary \ref{co:intertwining2} should not in general hold with those definitions. For a covariant calculus using Markovian connections see Cruzeiro-Malliavin \cite{Cruzeiro-Malliavin00}.

\section{Special case: no redundant noise}
\label{se:special-case}

Suppose $X(x): \R^m\to T_xM$ is injective for each $x$, so $p=m$
and $X$ trivialises $E$. Or equivalently $\F=\F^{x_0}$. Examples include left and right invariant stochastic differential equations on Lie groups, 
Example 2 of section \ref{sec:Examples} and the canonical stochastic differential equation on the orthonormal frame bundle, Example 4 of the same section.  Condition $(M)$ may not hold for the injective case, see for example 2C on page 24 of Elworthy-LeJan-Li \cite{Elworthy-LeJan-Li-book}. However we still get complete intertwining and all the results in this article.  In fact condition $(M_0)$ can be removed in the proof of the key Theorems,  Th \ref{th:I-star}, Th \ref{theorem-closability-2}, Co \ref{co:closed-range} and their corresponding Hilbert space valued versions. Furthermore the Conditions $K$ needed to define $\D^{p,k}$ holds trivially. The following theorem extends some of the results of Fang-Franchi \cite{Fang-Franchi-97} who were concerned with Brownian motion measure on Lie groups.

\begin{theorem}
\label{th:injective-case}
Suppose $X$ is injective. Then for $1<p<\infty$
\begin{enumerate}
\item  
$\I^*$ maps $\D^{p,k}(\C_{x_0}M;\R)$ isometrically onto $\D^{p,k}(\Omega;\R)$, $k=1,2,\dots$, and $$\I^*d^p=d^p\I^*$$
$$\I^*\pnabla^{p,(k)}d^p=(d^{p})^{(k+1)}\I^*$$
where, for $\psi\in L^p\Gamma \LL\left(\otimes^{k+1}\HH;\R\right)$,
$$\I^*(\psi): \Omega\to \LL(\otimes^{k+1}\HH;\R)$$ is given by $\I^*(\psi)=\left(\psi\circ \otimes ^{(k+1)}\overline{T\I}\right)\circ \I$.
\item 
On $\HH$-1-forms $(d^p)^*\I^*=\I^*(d^p)^*$. An $\HH$-vector field $V$ lies in $\Dom(\div^p_{\C_{x_0}M})$ if and only if
$(\Y V)\circ \I$ is in $\Dom(\div^p_\Omega)$; if so
$$(\div V)\circ \I=-\int_0^T\left\langle {\D\over ds}V(x_s),
X(x_s)dB_s\right\rangle_{x_s}$$
where the integral is a Skorohod integral.
\item The Laplacians, or `Ornstein-Uhlenbeck operators',  are conjugated by $\I^*$, as operators on $L^p$:
$$\I^*((d^p)^*d^p=(d^p)^*d^p\I^*.$$
\end{enumerate}
\end{theorem}

\begin{proof}
Since $T \I=\overline{T\I}$, it acts isometrically on $H$ with inverse $\Y$. The proof of Proposition \ref{pr:divergence-2} gives that
$h\in \dom(\div_\Omega^p)$ if and only if $T\I(h)\in \Dom(\div^p)$. Part 2 follows by (\ref{divergence-4}).
 For $k=1$, part 1 follows from part 2 since, by (\ref{I-phi}), $\I^*$ acts isometrically on $\HH$-1-forms as well as on functions. For $k>1$ and higher order Sobolev spaces, observe that for $G$ a separable Hilbert space,  $\D^{p,k}_{\F^{x_0};G}=\D^{p,k}(\Omega, G)$. Finally part 3 follows from part 1 and part 2.

\end{proof}
\begin{remark}
From Theorem \ref{th:injective-case}  we see Markov uniqueness holds if $X$ is injective, without assuming condition $(M_0)$. In fact the proof in section \ref{section-density} shows that $\Cyl^0\HH^*$ is total in $\Dom(d^*)$ on $\C_{x_0}M$ and so $\D^{2,1}\HH^*$ is dense in $\Dom(d^*)$   and $\overline{(d^p)^*  |_{ {\Cyl^0}\HH^*}}=\overline{(d^p)^*  |_{\D^{q,1}\HH^*}}$.  The argument leading to Theorem \ref{th:Markov-unique-2} in section \ref{section-Markov-unique}  proves $\D^{2,1}=W^{2,1}$ and  Markov uniqueness.
 For the stronger result of the essential self-adjointness of $d^*d|_{\Cyl}$ on $L^2(\C_{x_0}M;\R)$ see the method of Aida \cite{Aida93}.
\end{remark}

For completeness and as an example we show that for Lie groups $G$ with left invariant connection our connection $\pnabla$ agrees with the `Levi-Civita' connection used previously, \cf  Freed \cite{Freed88}, Fang-Franchi \cite{Fang-Franchi-97}, Aida \cite{Aida95},  Driver-Lohrenz \cite{Driver-Lohrenz}. See also Shigekawa  \cite{Shigekawa97}, Shigekawa-Taniguchi \cite{Shigekawa-Taniguchi96}. For this let $M=G$ with left invariant stochastic differential equation as in Example 2 of section \ref{sec:Examples}. Since the adjoint connection $\hat \nabla$ of $\nabla$ is the flat right invariant connection we have $\HH$-vector fields $V^y(\sigma)$ given by
$V(\sigma)_t:=TR_{\sigma(t)}(y_t)$, $y\in L^{2,1}_0([0,T];\g)$, for 
$TR_g: \g= T_eG\to T_gM$ the derivative of right translation $R_g$ by $g\in G$.
\begin{proposition}
For $y,z$ in $L^{2,1}_0([0,T];\g)$,
$$\left(\pnabla_{V^z(\sigma)} V^y \right)_t=TR_{\sigma(t)}[\dot y(t), z(t)].$$
\end{proposition}
\begin{proof}
By definition of $X$ and $\pnabla$ and since the Ricci curvature is zero and there is no `drift',
\begin{eqnarray*}
\left(\pnabla_{V^z(\sigma)} V^y \right)_t&=&
TL_{\sigma(t)}\left[ 
d\left( \rho\mapsto TL_\rho^{-1} {\D\over d\cdot} V_\cdot ^y\right)_\sigma(V^z(\sigma))\right]_t\\
&=&TL_{\sigma(t)}\left[ 
d\left( g\mapsto TL_g^{-1} TR_g \dot y_t\right)_{\sigma(t)}  (TR_{\sigma(t)}(z_t))\right]\\
&=&TR_{\sigma(t)}[\dot y(t), z(t)]
\end{eqnarray*}
as expected.
\end{proof}

\section{Appendix: The conditional expectations of exponential martingales}
\label{Appen}
Let $\varepsilon(a)$ be the exponential martingale 
$$\varepsilon(a):=\exp{\left(\int_0^T \langle \dot  a_s, dB_s \rangle-{1\over 2}\int_0^T|\dot a_s|^2 ds\right)}$$
for $a\in H$ and set $\tilde \varepsilon(a)=\E\{ \varepsilon(a)|\F^{x_0}\}$. By Proposition \ref{pr:chaos-expansion}
we know $\tilde\varepsilon(a)$ is in $\D^{2,1}_{\F^{x_0}}$. As evidence that $\E\{-|\F^{x_0}\}$  maps $\D^{2,1}$ into itself we show, c.f. Remark \ref{Remark-7.9}:
\begin{proposition}

There is a constant $C$ such that for all $a\in H$ 
$$\|\tilde\varepsilon(a)\|_{\D^{2,1}}\le C\|\varepsilon(a)\|_{\D^{2,1}}.$$
\end{proposition}
\begin{proof}
Set  $$\varepsilon_t=\varepsilon_t(a)=\exp{\left(\int_0^t \langle \dot  a_s, dB_s \rangle-{1\over 2}\int_0^t|\dot a_s|^2 ds\right)}.$$
Denote by $d^H$ the $H$-derivative to distinguish it from the stochastic differential. 
Then $\E|\varepsilon_t|^2=\exp\left({\int_0^t|\dot a_s|^2 ds}\right)$ and 
$\|\varepsilon(a)\|^2_{\D^{2,1}}=(1+\|a\|^2_H) \exp\left(\|a\|^2_H\right)$.
Conditioning the stochastic differential equation
$$d\varepsilon_t= \left\langle \dot a_t, dB_t\right\rangle \varepsilon_t$$
on $\F^{x_0}$ shows that $\tilde\varepsilon_t\equiv\tilde\varepsilon_t(a)= \E\{\varepsilon_t|\F^{x_0}\}$ satisfies:
$$d\tilde \varepsilon_t= \tilde\varepsilon_t \left\langle X(x_t) \dot a_t, X(x_t)dB_t\right\rangle.$$
So
$$\tilde\varepsilon_t=\exp{\left(\int_0^t \langle X (x_s)\dot  a_s, X(x_s) dB_s \rangle-{1\over 2}\int_0^t|X(x_s) \dot a_s|^2 ds\right)}$$
and $$\E|\bar\varepsilon_t|^2=\E \exp\left({\int_0^t|K^\perp(x_s) \dot a_s|^2 ds}\right),$$
where $K^\perp: M\times \R^m\to [\ker(X)]^\perp$ is the projection map $Y(x)X(x)$. From this, for $h\in H$,
\begin{eqnarray*}
d^H(\tilde\varepsilon)(h)
&=&\tilde\varepsilon\left\{\int_0^T \left\langle
X(x_s)\dot a_s, X(x_s)(\dot h_s)  ds+ \nabla_{T\I_s(h)} X dB_s  \right \rangle 
\right.\\\\
&& \left. +\int_0^T\left\langle \nabla_{T\I_s(h)}X(\dot a_s), X(x_s)dB_s\right \rangle 
-\int_0^T \left\langle  \nabla_{T\I_s(h)} X(x_s)(\dot a_s), X(x_s)\dot a_s\right\rangle ds\right\}.
  \end{eqnarray*}
  By Theorem \ref{th:I-star}, setting
  $\phi^{\#}=\tilde\varepsilon \W_s(X(x_\cdot)\dot a_\cdot)$, then
  $$\tilde\varepsilon \int_0^T \left\langle
X(x_s)\dot a_s, X(x_s)(\dot h_s  ds)+ \nabla X(T\I_s(h)) dB_s  \right \rangle 
=\I^*(\phi)(h).$$
Note that $\|\phi\|^2=\E \tilde \varepsilon^2\int_0^T |X(x_s) \dot a_s|^2 ds\le |a|^2_H \E \tilde \varepsilon^2$.
  
 Take any Riemannian metric on $TM$ extending that of $E$ and for $ e\in \R^m$ define: 
  $$Z_x(e): E_x\to T_xM$$
to be the adjoint of $\nabla_ -X(e): T_xM\to E_x$ for each $x$ in $M$. Then
\begin{eqnarray*}
&&\tilde\varepsilon\left\{ \int_0^T\left\langle \nabla_{T\I_s(h)}X(\dot a_s), X(x_s)dB_s\right \rangle 
-\int_0^T \left\langle  \nabla_{T\I_s(h)} X(x_s)(\dot a_s), X(x_s)\dot a_s\right\rangle ds\right\}\\
&=&\tilde\varepsilon\left\{ \int_0^T\left\langle T\I_s(h), Z_{x_s}(\dot a_s)\left(X(x_s)dB_s\right) \right \rangle 
-\int_0^T \left\langle  T\I_s(h) ,Z_{x_s}(\dot a_s)\left(X(x_s)\dot a_s)\right) \right\rangle ds\right\},
 \end{eqnarray*}
which can be verified to equal to $\I^*(\psi)(h)$, defined by Theorem  \ref{th:I-star},  where $$\psi^{\#}= \tilde\varepsilon \W_t\left(\Pi(W_\cdot^{-1})^*\int_\cdot^T W_r^* Z_{\sigma_r}(\dot a_r)X(\sigma_r)\left(dB_r-\dot a_rdr\right)\right).$$
We have used the expression (\ref{derivative-Ito-map}) of $T\I$ for the verification. 
It remains to estimate 
\begin{eqnarray*}
 \E|\I^*(\psi)|^2&=& \int_{C_{x_0}M}\int_0^T  \left|{\D\over ds}\psi^{\#}_s \right|^2ds \mu_{x_0}(d\sigma)\\
&=&\E \tilde \varepsilon^2 \int_0^T  \left(\Pi(W^{-1}_t)^*\int_t^T W_r^* Z_{x_r}(\dot a_r)X(x_r)(dB_r-\dot a_r dr)\right)^2 dt.
\end{eqnarray*}
For this let $\{y_t: 0\le t\le T\}$ be the solution to the stochastic differential equation
\[\left\{\begin{array}{lll}
dy_t&=&X(y_t)\circ dB_t+A(y_t)dt+2X(y_t)(\dot a_t)dt\\
y_0&=&x_0.
\end{array}\right.\]
Then by the Girsanov-Maruyama theorem, 
\begin{eqnarray*}
&& \E|\I^*(\psi)|^2\\
&=& \E\exp\left(\int_0^T|K^\perp(y_s)\dot a_s|^2ds\right) \int_0^T
\left(\Pi(W^{-1}_t)^*\int_t^T W_r^* Z_{y_r}(\dot a_r)X(y_r)(dB_r+\dot a_r dr)\right)^2 dt\\
&\le& 2\E\exp\left(\int_0^T|K^\perp(y_s)\dot a_s|^2ds\right) \int_0^T
\left(\Pi(W^{-1}_t)^*\int_t^T W_r^* Z_{y_r}(\dot a_r)(X(y_r)dB_r)\right)^2 dt\\
&&+ 2\E\exp\left(\int_0^T|K^\perp(y_s)\dot a_s|^2ds\right) \int_0^T
\left(\Pi(W^{-1}_t)^*\int_t^T W_r^* Z_{y_r}(\dot a_r)(X(y_r)\dot a_r) dr\right)^2 dt\\
&\le& \constant  \|a\|_H^2\exp( \|a\|_H^2)\\
&&+ \exp( \|a\|_H^2)\E\exp\left(-\int_0^T|K(y_s)\dot a_s|^2ds\right)\cdot \\
&&\hskip 100pt \int_0^T
\left(\Pi(W^{-1}_t)^*\int_t^T W_r^* Z_{y_r} (K(y_r)\dot a_r )X(y_r)(\dot a_r) dr\right)^2  dt\\
  \end{eqnarray*}
  where $K(x)=1-K^\perp(x)$, using the fact that $\nabla X=\nabla X\circ K$.
  Now 
  \begin{eqnarray*}
  &&\E\exp\left(-\int_0^T|K(y_s)\dot a_s|^2ds\right)\cdot 
 \int_0^T
\left(\Pi(W^{-1}_t)^*\int_t^T W_r^* Z_{y_r} (K(y_r)\dot a_r )X(y_r)(\dot a_r) dr\right)^2  dt\\
&\le&\const \E\exp\left(-\int_0^T|K(y_s)\dot a_s|^2ds\right)\cdot 
\int_0^T\left\{\int_t^T |K(y_r)\dot a_r|^2 dr\int_t^T|K^\perp(y_r)\dot a_r|^2 dr\right\}dt\\
&\le&\const  |a|_H^2 \E\int_0^T \left(\int_t^T \exp{\left(-\int_r^T|K(y_s)\dot a_s d_s|^2ds\right)}
|K^\perp(y_r)\dot a_r|^2 dr \right)dt \\
&\le&\const  |a|_H^2 \E\int_0^T  \exp{\left(-\int_t^T|K(y_s)\dot a_s d_s|^2ds\right)} dt \\
&\le&\const  |a|_H^2.
   \end{eqnarray*}
   From this we see $\E|\I^*(\psi)|^2\le \const |a|^2_H e^{|a|_H^2}$ as required.
\end{proof}

\subsubsection*{Acknowledgement}
This research was partially supported by EPSRC research grant GR/H67263 and NSF  research grant DMS 0072387.  Xue-Mei Li would like to acknowledge the support from the Royal Society and the Leverhulme Trust for her Senior Research Fellowship. The research benefited from our contacts with many colleagues, especially S. Aida, S. Albeverio, A. Eberle, S. Fang, Y. LeJan, Z. Ma, P. Malliavin, and M. R\"ockner.

\bigskip

\noindent{\bf \small  Addresses:}\\
{ \small K.D. Elworthy,
 Mathematics Institute, University of Warwick, Coventry CV4 7AL, U.K.}\\
{\small Xue-Mei Li (Xue-Mei Hairer), Mathematical Sciences, 
Loughborough University, 
Loughborough, LE11 3TU
U.K. \email xue-mei.li@lboro.ac.uk. }

\end{document}